\documentclass[11pt,reqno,english]{article}

\usepackage{amsmath,amsfonts,amssymb,graphicx,amsthm,url}
\usepackage[noadjust]{cite}
\usepackage{stmaryrd}
\usepackage{mathrsfs,booktabs,tabularx}
\usepackage{microtype,enumerate,tikz,xcolor,fancyhdr,needspace}
\usetikzlibrary{arrows.meta}
\usepackage{xifthen,xcolor,tikz,setspace}
\usetikzlibrary{decorations.pathmorphing,patterns,shapes,calc,decorations}
\usetikzlibrary{decorations.pathreplacing}
\usepackage{mathtools,upgreek,paralist}

\usepackage[final]{showkeys} %add in 'final' into parameter to remove showkeys

\definecolor{refkey}{gray}{.75}
\definecolor{labelkey}{gray}{.5}
\usepackage[algo2e,boxed,vlined,algoruled]{algorithm2e}
\usepackage{comment}
\usepackage[shortlabels]{enumitem}

\usepackage{graphbox}

\usepackage{extarrows}

\definecolor{inkblue}{HTML}{254C70}
\definecolor{softblue}{HTML}{DCEAF4}
\definecolor{plusgreen}{HTML}{276B53}
\definecolor{softgreen}{HTML}{E8F2EC}
\definecolor{frozengray}{HTML}{E7E8EA}

\colorlet{DarkGreen}{green!50!black}

\newcommand{\Osc}{\operatorname{Osc}}

\usepackage[letterpaper]{geometry}
\usepackage[colorinlistoftodos]{todonotes}
\presetkeys{todonotes}{inline, color=green}{}

\usepackage{accents}
\usepackage[algo2e,boxed,vlined,algoruled]{algorithm2e}

\usepackage{authblk}

\usepackage[small]{caption}
\usepackage[colorlinks=true, citecolor=blue]{hyperref}
\numberwithin{equation}{section}

\newcommand{\Var}{\operatorname{Var}}

\newcommand{\E}{\mathbb E}

\renewcommand{\epsilon}{\varepsilon}

\usepackage{color}
 \definecolor{refkey}{gray}{.5}
 \definecolor{labelkey}{gray}{.5}
\definecolor{light}{gray}{.9}

\usepackage{soul}

\usepackage[capitalise]{cleveref}

\newtheorem{theorem}{Theorem}[section]
\newtheorem*{theorem*}{Theorem}
\newtheorem{lemma}[theorem]{Lemma}

\crefname{claim}{Claim}{Claims}
\newtheorem{proposition}[theorem]{Proposition}

\newtheorem{observation}[theorem]{Observation}

\newtheorem{corollary}[theorem]{Corollary}

\theoremstyle{definition}

\newtheorem{definition}[theorem]{Definition}

\newtheorem*{definition*}{Definition}

\newtheorem{remark}[theorem]{Remark}
\newtheorem*{remark*}{Remark}

\newcommand{\tmix}{t_{\textsc{mix}}}

 \renewcommand{\epsilon}{\varepsilon}

\DeclareMathOperator{\var}{Var}
\DeclareMathOperator{\dist}{dist}

\DeclareMathOperator{\cov}{Cov}

\newcommand{\tv}{{\textsc{tv}}}

\newcommand{\eps}{\varepsilon}

\title{Polynomial mixing for 2D Ising with plus boundary conditions}

\author{Shirshendu Ganguly\thanks{Department of Statistics, UC Berkeley. Berkeley, CA. \url{sganguly@berkeley.edu}} { and }Reza Gheissari\thanks{Department of Mathematics, Northwestern University. Evanston, IL. \url{gheissari@northwestern.edu}}   }

\date{November 23, 2026}

\begin{document}

\maketitle

\vspace{-0.5cm}
\begin{abstract}
    We study the Glauber dynamics of the 2D Ising model at low temperatures. 
    In the absence of boundary conditions, the mixing time on an $n\times n$ box is exponential in $n$ due to a bottleneck between metastable plus and minus phases. When the symmetry is broken via all-plus boundary conditions, it has been predicted in the physics literature since work of Lifshitz (1962) that from an all-minus initialization, the minus droplet shrinks via mean-curvature flow on timescales that are quadratic in $n$, and thus, the system equilibrates in polynomial time. On the rigorous side, the strongest result was a quasi-polynomial mixing time upper bound due to Lubetzky, Martinelli, Sly and Toninelli (2013); a proof of polynomial mixing remained one of the major challenges in the field. 
    We prove this polynomial mixing time bound
    at all low temperatures. 
    
    The key step is to establish fast mixing of the interface separating the plus and minus phases on a box with Dobrushin boundary conditions. 
    In simpler models of height functions, polynomial bounds tend to rely on exact harmonic properties and connections to a discrete heat equation. 
    However, the Ising interface exhibits microscopic defects that necessitate the development of a more robust method. 
     We address this via a new multiscale approach to bound mixing times of general disordered random-walk bridges that homogenize to a Brownian bridge at mesoscopic scales. Another notable ingredient is a novel application of censoring to reduce the mixing time problem to proving only a weak Poincar\'e inequality, or warm-start mixing bound, allowing us to assume suitable regularity conditions on the interface and bubbles.
\end{abstract}

\tableofcontents

\section{Introduction}\label{sec:introduction}

The Ising model is one of the most studied models of phase transitions in statistical physics and probability theory. On a graph $G = (V,E)$ at inverse temperature $\beta$, it is the following probability distribution over $\{-1,+1\}$ spin assignments to the vertices of $V$: 
\begin{equation}\label{eq:Ising-distribution}
    \pi_{G,\beta}(\sigma) \propto \exp\Big(  \beta \sum_{v\sim w} \sigma_v \sigma_w\Big)\,,
\end{equation}
where $u\sim v$ indicates that $\{u,v\}$ is an edge in $E$. 
We will focus on the model in $\mathbb Z^2$, which is of particular interest~\cite{McCoyWu}. Early work of Peierls and Kramers and Wannier showed that the Ising model on $\mathbb Z^2$ undergoes a sharp phase transition at a critical point $\beta_c = \frac{1}{2}\log (1+\sqrt{2})$. This critical point marks a transition between a high-temperature regime $\beta<\beta_c$ where correlations between the spins at a pair of vertices decay exponentially in their distance, and the infinite-volume Gibbs measure is unique, and a low-temperature regime $\beta>\beta_c$ where long-range order emerges and distinct plus and minus phases, in which one of the two spins dominates, coexist.

The study of dynamics of the Ising model dates back to Glauber~\cite{Glauber}, who introduced the following Markov chain to model the microscopic evolution of an Ising system out of equilibrium. 
The continuous-time Ising Glauber dynamics on (finite subsets $\Lambda$ of) $\mathbb Z^2$ initialized with $X_0 = x_0$, assigns each vertex a rate-one Poisson clock, and if the clock at vertex $v$ rings at time $t$, it updates 
\begin{align}\label{eq:Glauber-update-rule}
    X_t(v)  = \begin{cases}
        +1 & \text{w.\ prob.\ }\propto \exp( \beta \sum_{w\sim v} X_{t^-}(w))\\ - 1 & \text{w.\ prob.\ }\propto \exp( - \beta \sum_{w\sim v} X_{t^-}(w))
    \end{cases}\,,
\end{align}
while keeping all other spins unchanged. In the presence of boundary conditions (fixed assignment of spins to the vertices of $\mathbb Z^2 \setminus \Lambda$),~\eqref{eq:Ising-distribution} and~\eqref{eq:Glauber-update-rule} are evaluated also counting the edges between $\Lambda$ and $\partial \Lambda = \{v\in \mathbb Z^2 \setminus \Lambda: d(v,\Lambda) =1\}$. We use $+$ to denote the all-plus boundary on $\mathbb Z^2 \setminus \Lambda$.  

The mixing time (quantified from a worst-case initial state) is defined as 
\begin{align}
    \tmix  = \inf \big\{t \ge 0: \max_{\sigma \in \{\pm 1\}^{\Lambda}} \|\mathbb P(X_t \in \cdot \mid X_0 =\sigma)  - \pi_{\Lambda,\beta}\|_{\tv} \le 1/4\big\}\,,
\end{align}
and quantifies how long it takes an Ising system to equilibrate. The mixing time of the Ising dynamics on an $n\times n$ box $\Lambda_n$ has seen much attention and undergoes a transition at $\beta_c$ paralleling the static phase transition. At $\beta<\beta_c$, the system mixes rapidly in $c_\beta \log n$ time due to the exponential decay of correlations~\cite{MaOl1,LS-information-percolation}, while at $\beta>\beta_c$ in the absence of boundary conditions, the plus and minus phases of the system are both metastable and it takes $\exp(\Omega(n))$ time to transit between them and equilibrate~\cite{Schonmann-LD,Martinelli-phase-coexistence,CGMS96}. In the last few decades, the study of low-temperature dynamics has focused on in what sense the transit between the phases is the only obstruction to fast mixing, and how that question interplays with the geometry of the underlying graph. 

\begin{figure}[t]
\centering
    \includegraphics[width = 0.27 \textwidth]{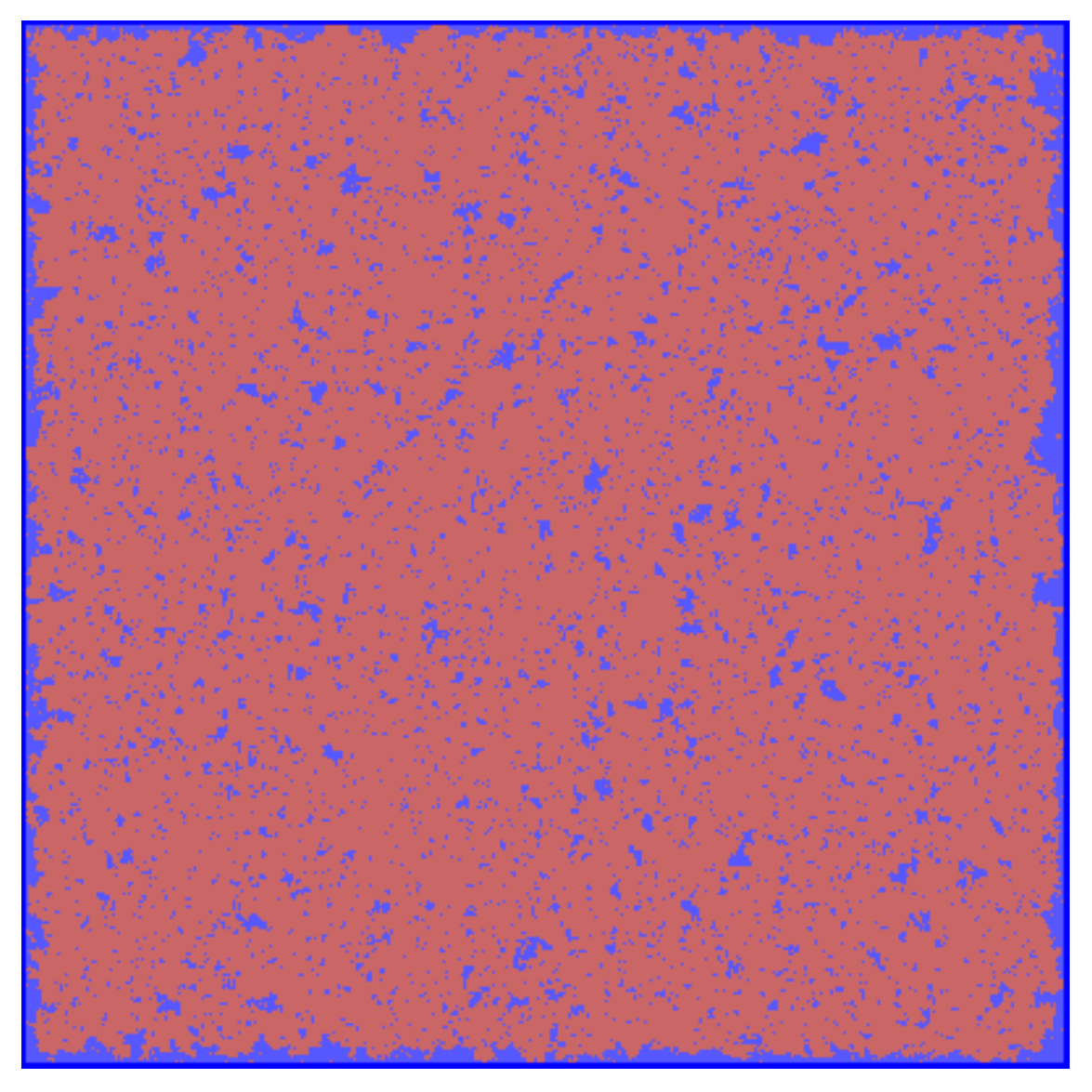} \qquad
    \includegraphics[width = 0.27 \textwidth]{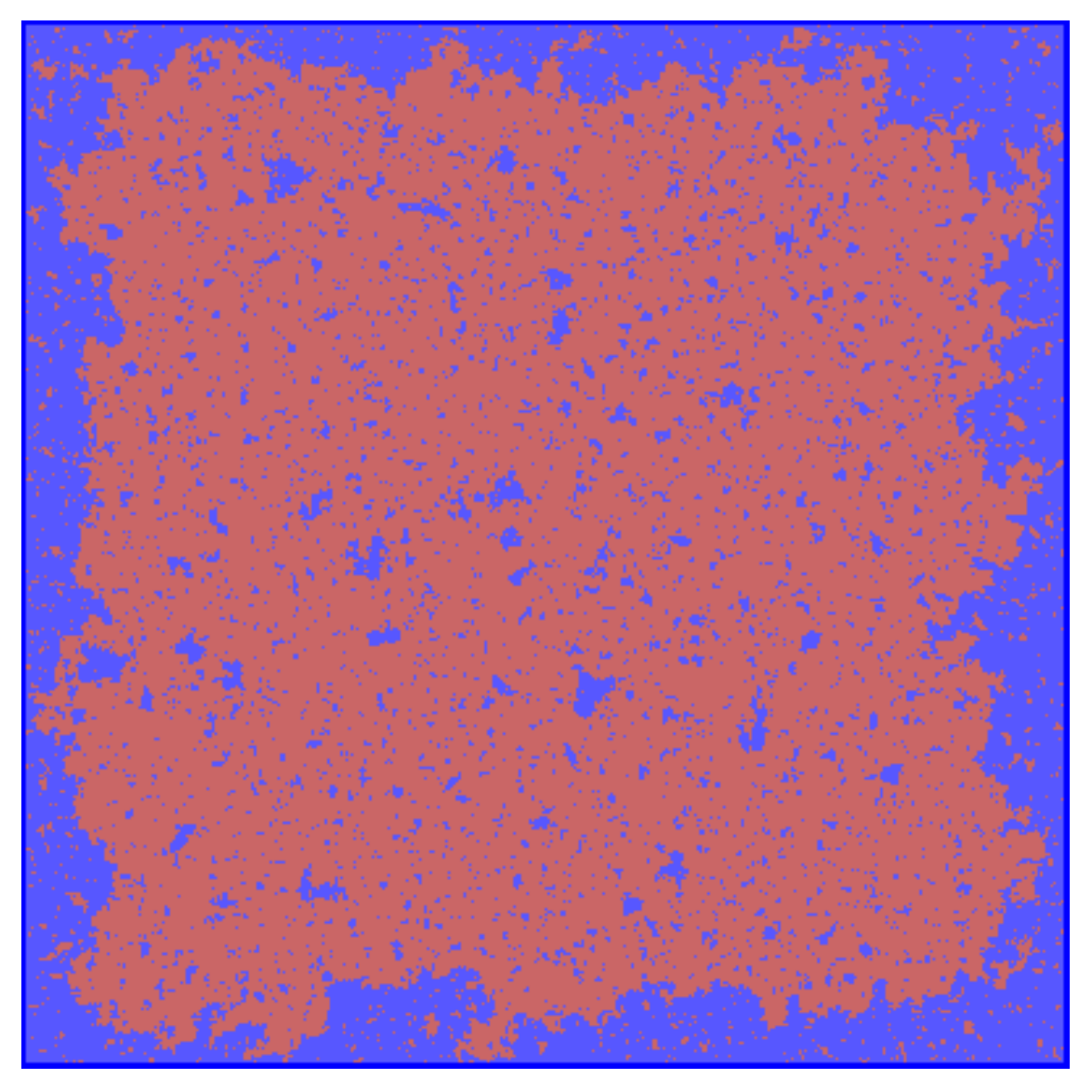} \qquad 
    \includegraphics[width = 0.27 \textwidth]{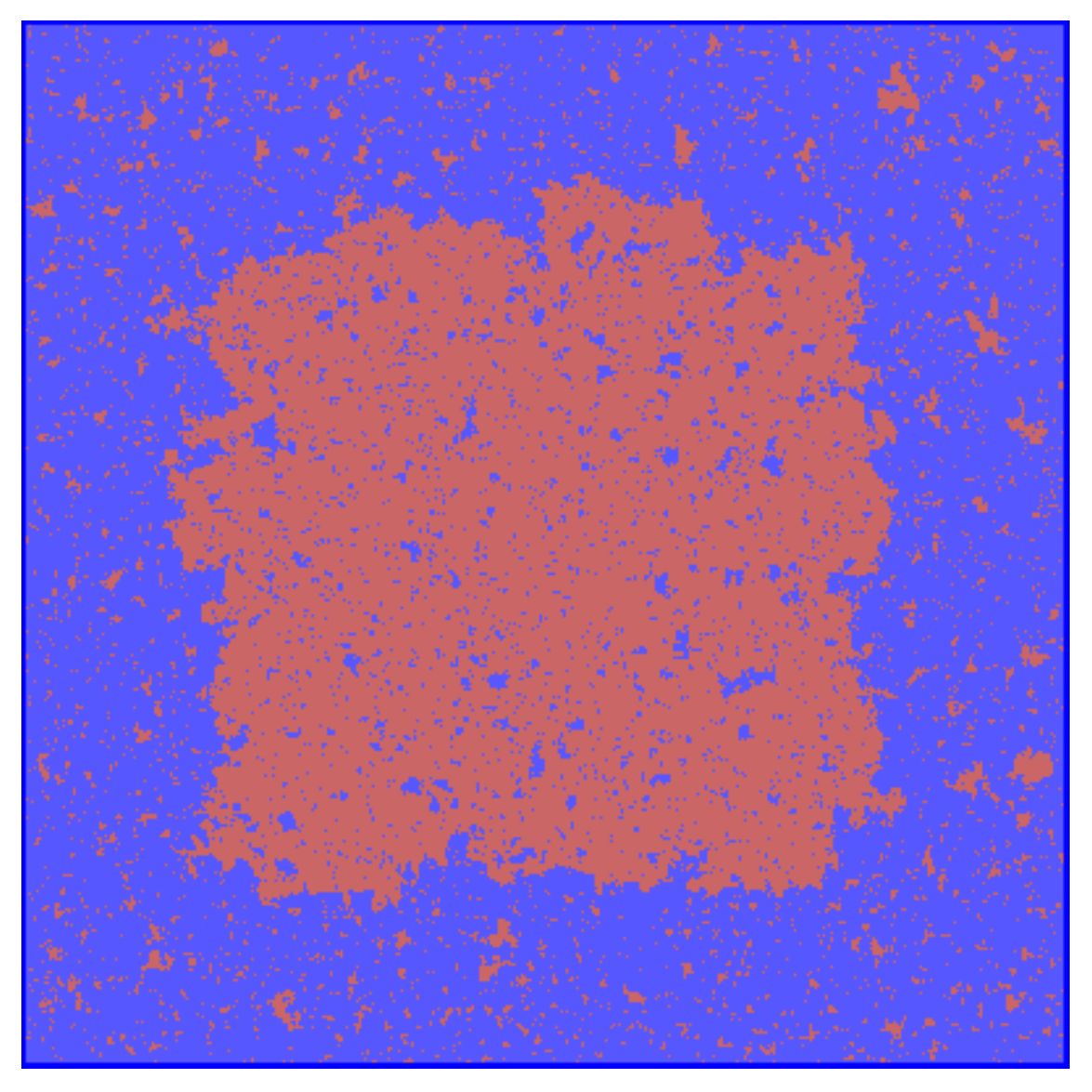}
\caption{Snapshots over time of low-temperature Ising Glauber dynamics on a 400$\times$400 box in $\mathbb Z^2$, initialized from all minus (red) in a box with plus (blue) boundary conditions.}\label{fig:mean-curvature-flow}
\end{figure}

A specific formulation of this is what effect a weak breaking of the symmetry by placing plus boundary conditions on $\Lambda_n$ has on the low-temperature out-of-equilibrium dynamics. 
In the physics literature, Lifshitz~\cite{Lifshitz} predicted that in this setting, the convergence to equilibrium from the all-minus initialization would happen on $O(n^2)$ timescales, as the plus phase slowly encroaches on the locally stable minus droplet: see Figure~\ref{fig:mean-curvature-flow}. 

The aim of providing a rigorous foundation to this \emph{Lifshitz law} has become fundamental both for understanding how rapidly a low-temperature
system equilibrates within a pure phase, and towards understanding dynamics of boundaries between coexisting metastable phases. 
In the first context, this question has become canonical for the task of getting polynomial mixing in a spin system where there is no correlation decay and no positive curvature/contractivity of the dynamics: for instance~\cite{MaTo} went as far as calling it ``a kind of holy grail for the subject".  
In the second context, the setup is connected to a rich predicted theory in statistical physics for interface motion. 
Namely, on timescales~$(t n^2)_{t\in [0,\infty)}$ and rescaled in space by~$n^{-1}$, the boundary of the minus droplet should converge to a deterministic curve evolving by anisotropic mean-curvature flow~\cite{Spohn93}, with its inward velocity expressed in terms of the  surface tension entering the Wulff functional~\cite{KarmaLobkovsky05}. Furthermore, at a finer scale, the fluctuations of the droplet boundary about this deterministic curve are expected to be described by a stochastic PDE~\cite{Spohn93,GrantGunton83}.

Mathematically, however, the understanding has been more limited. As hinted at above, the key mechanism dictating the equilibration of Ising dynamics with plus boundary conditions is the evolution of the interface bounding the minus droplet under Glauber dynamics. 
At \emph{zero-temperature}, it was shown that started from all-minus with plus-boundary conditions, the droplet shrinks via mean-curvature flow and the system absorbs into the all-plus configuration in $n^{2+o(1)}$ time~\cite{lacoin2014zero-2D} (see also earlier bounds in~\cite{FoScSi02,Morris-zero-temp}). Similarly, $n^{2+o(1)}$ bounds are known for the approximating $(1+1)$D solid-on-solid (SOS) model of interfaces~\cite{CPT-monotone-surfaces}. 
On the one hand, at low-but-positive temperature, the Ising interface at large scale looks like a Brownian bridge~\cite{GreenbergIoffe}, and already at mesoscopic scales behaves like its SOS approximation. However, the 
essential challenge is that the microscopic defects of the model at $\beta>0$ (bubbles outside the interface, and overhangs of the interface) can destroy the delicate balance of an exact Laplacian drift on the interface heights that gives the $n^{2+o(1)}$ mixing in the approximating models. 

This challenge has been well-known for the low-temperature Ising dynamics with plus boundary conditions since work of~\cite{Martinelli-phase-coexistence} who proved an $e^{n^{1/2+o(1)}}$ bound on the mixing time at $\beta$ sufficiently large. A lower bound of $n^{2-o(1)}$ on the mixing time was shown in~\cite{Bodineau-Martinelli02}.  Nearly two decades  after~\cite{Martinelli-phase-coexistence},~\cite{MaTo} proved a sub-exponential $e^{n^{o(1)}}$ upper bound through a novel comparison to modified boundary conditions that help drive the interface down. ~\cite{LMST} pushed that idea to its limits, obtaining quasi-polynomial $e^{O((\log n)^2)}$ mixing at all low temperatures $\beta>\beta_c$. (Very recently, towards a different purpose of studying phase ordering, it was shown in~\cite{gheissari2026rapidphaseorderingising} that the quasi-polynomial upper bound could be taken to be uniform over all large $\beta$). 

Our main result in this paper is polynomial mixing of Ising Glauber dynamics on squares in $\mathbb Z^2$ with all-plus boundary conditions at all low temperatures.

\begin{theorem}\label{thm:plus-bdy-mixing}
    Fix $\beta>\beta_c$. There exists a constant $C_\beta>0$ such that Ising Glauber dynamics on the domain $\Lambda_n = [1,n]^2 \cap\mathbb Z^2$ with plus boundary conditions has $\tmix \le n^{C_\beta}$. 
\end{theorem}

(Together with the polynomial bound at criticality of~\cite{LS-critical-2D-Ising} and the high-temperature bound of~\cite{MaOl1},  Theorem~\ref{thm:plus-bdy-mixing} implies Glauber with plus boundary mixes in polynomial time at all fixed~$\beta$ in $d=2$. This resolves Question 1 in Open Problems of the book~\cite{LP}.)

As an immediate  consequence  we also get an effective decay rate for the plus phase autocorrelation in $\mathbb Z^2$ at all low temperatures. Fisher and Huse~\cite{FisherHuse} predicted that the spin-spin autocorrelation in the low-temperature plus phase should decay as $e^{-\sqrt{t}}$. Corollary 4 of~\cite{LMST} had previously shown a super-polynomial relaxation rate of $e^{ -e^{\sqrt{\log t}}}$ in $d=2$. With our Theorem~\ref{thm:plus-bdy-mixing} in hand and following the same proof as that of~\cite[Corollary 1.10]{MaTo}, we deduce stretched exponential autocorrelation decay in $d=2$. 

\begin{corollary}\label{cor:plus-phase-autocorrelation}
    Fix $\beta>\beta_c$. There exist  constants $c,C>0$ such that Ising Glauber dynamics on $\mathbb Z^2$ started from the infinite-volume plus measure $\pi^+_{\mathbb Z^2}$ satisfies for all $t\ge 0$, 
    \begin{align*}
        \text{Cov}_{\pi_{\mathbb Z^2}^+}(X_{0}(0,0), X_{t}(0,0)) \le  C e^{-t^{c}}\,.
    \end{align*}
\end{corollary}

\subsection{Challenges and ideas}

As in previous works on the subject, Theorem~\ref{thm:plus-bdy-mixing} follows fairly straightforwardly from a mixing-time result that holds on domains with Dobrushin boundary conditions, meaning that the boundary conditions are $+$ above a line through $\Lambda_n$, and $-$ below the line. Whereas the preceding works~\cite{Martinelli-phase-coexistence,MaTo,LMST} had used boxes with $-$ boundary on the South side, and $+$ boundary on the other three sides, it ends up being preferable for us to work with the \emph{diagonal Dobrushin boundary conditions}. Those are where the boundary line is $y=n-x$ as depicted in Figure~\ref{fig:diag-Dobrushin-censoring}, left, so the boundary conditions are $-$ on the West and South sides of $\Lambda_n$ and $+$ on the North and East sides.

\begin{theorem}\label{thm:diagonally-Dobrushin-mixing}
    Fix $\beta>\beta_c$. There exists a constant $C$ such that the Ising Glauber dynamics on $\Lambda_n = [1,n]^2 \cap\mathbb Z^2$ with diagonal Dobrushin boundary conditions has $\tmix \le n^{C}$. 
\end{theorem}

At the end of the proof sketch, in Section~\ref{subsec:Dobrushin-to-all-plus} we will give a short deduction of Theorem~\ref{thm:plus-bdy-mixing} from Theorem~\ref{thm:diagonally-Dobrushin-mixing}, so the rest of our discussion will focus on proving polynomial mixing on $n\times n$ boxes with the diagonally Dobrushin boundary conditions.

In the box $\Lambda_n$ with diagonal Dobrushin boundary conditions, the main object of interest is the \emph{interface}, denoted by $\Gamma$, which is the random curve separating the connected $-1$-component of the West and South boundaries of $\Lambda_n$ from the connected $+1$-component of the North and East boundaries (see Figure~\ref{fig:contour-decomposition}, and Definition~\ref{def:contour-decomposition} for a formal definition). Note that this interface can have self-intersections, and it is not a height function as it can have what are called \emph{overhangs}, meaning it backtracks horizontally. Nonetheless, as evident from the left picture in Figure~\ref{fig:diag-Dobrushin-censoring}, such microscopic defects begin to disappear at large scales on which it looks like a random walk bridge of slope $-1$ connecting vertices $(1,n) \to (n,1)$.

The predictions for the mixing time of the Ising model with Dobrushin boundary conditions at low temperatures start from ignoring these microscopic defects. This is done either by taking a zero-temperature limit of this interface, yielding is a uniform down-right Dyck path from $(1,n) \to (n,1)$ (equivalent to the corner flip chain/the symmetric exclusion process) or the softer SOS approximation where it is a random walk bridge with geometric increment distributions. 
The works proving polynomial bounds for these kinds of random interface models stem from the seminal work of~\cite{Wilson}. Wilson's approach was to observe that in a random-walk bridge, when a Glauber update happens at column $i$, the new expected height in that column is exactly the average of the heights of its two neighbors. As a result, the expected height $\mathbb E[X_t \mid X_0]$ solves the discrete heat equation, and its distance to stationarity decays exponentially with rate given by the first eigenvalue of the discrete Laplacian on $\{1,...,n\}$ with Dirichlet boundary conditions, which is $\Theta(1/n^2)$. 

This method has proven to be very powerful and has been used to establish $n^{2+o(1)}$ mixing time for a variety of interface models which are expected to be driven by mean-curvature flow. This includes the zero-temperature Ising interface~\cite{lacoin2014zero-2D}, the $(1+1)D$ SOS model~\cite{MaSi-SOS,CPT-monotone-surfaces}, more general $(1+1)D$ height functions with convex potential~\cite{Lacoin-adjacent-transposition,CLL-nabla-phi,yang2024cutoffsimpleexclusionprocess}. It has also seen application in $(2+1)D$ surface dynamics like the dimer model~\cite{CPT-monotone-surfaces,aggarwal2026mixingtimesglauberdynamics,Laslier-Toninelli,Toninelli-ICM-review}, the Gaussian free field~\cite{GaGhcutoff,bristiel2026mixingprofileglauberdynamics}, and the zero-temperature Ising droplet in higher dimension~\cite{CMST-zero-temperature-3d-ising,Lacoin-Lifshitz-any-dimension}. 

However, given it relies on eigenfunctions of the Laplacian, the method lacks a sense of robustness. Namely, even for models of dynamics on random walk bridges, in the absence of an exact Laplacian representation, there is a notable void of methods for obtaining polynomial bounds on the mixing time. This is evidenced both by the difficulty of getting a polynomial bound for the Ising interface, and by the lack of effective upper bounds for the $(2+1)D$ SOS or integer-valued GFF dynamics. While all of these models roughly have the property that under local resampling, the height in the middle is the harmonic average of the heights on the boundary, this is not \emph{exactly} true due either to perturbations of the interface through bubbles and overhangs (in the Ising case) or rounding effects (in the integer-valued $(2+1)D$ surfaces).

\medskip
At a very high level, our strategy has two main innovations. The first is to leverage a new use of censoring~\cite{PWcensoring} together with what are known as \emph{weak Poincar\'e inequalities}, to restrict our analysis  to Ising configurations that are ``nice" in the sense that all the potentially worrying bubbles and overhangs are truly microscopic (at most $O(\log n)$ in size). This lets us focus on the mixing of just the interface $\Gamma$ as the mixing of the configuration off the interface is easier to show to be fast. 

The second is our way of analyzing the mixing time of the interface $\Gamma$ directly.  
To handle the microscopic defects, we coarsen the dynamics on the interface  to scales at which the path begins to exhibit more Brownian-like properties. More precisely, we resample it on overlapping $C \log n$-width strips for a large enough $C$. At that resolution, the microscopic defects (bubbles and overhangs) start to become negligible and the mesoscopic diffusivity of the interface kicks in. The path exhibits renewal points along with local central limit theorems which allow us to couple interfaces on segments of width $m$ with endpoints that are diffusively separated by $\sqrt{m}$, in the middle $m/3$ of the segment with constant probability. This coupling is used to reduce the width under consideration by a factor of $2/3$, at a constant cost in the mixing time analysis. Recursively applying that $O(\log n)$ times to get from width $n$ to the minimal $O(\log n)$ width gives a $C^{O(\log n)} = \text{poly}(n)$ bound.

\section{Proof sketch}

\subsection{From weak Poincar\'e inequality to mixing time}\label{sec:WPI-to-gap}

A key tool that allows us to formalize the notion that mixing is governed by the convergence time of the interface (the bubbles above and below the interface being negligible because they are sub-critical excitations away from a dominant phase) is going to be to reduce our task to proving a Weak Poincar\'e inequality (WPI) with Dobrushin boundary conditions. To help build intuition, one may interpret a WPI as a Poincar\'e inequality (PI) but only restricted to a high probability part of the state space. 

\begin{definition}\label{def:WPI}
    A Markov chain $P$ on a finite state space $\Omega$ with stationary distribution $\mu$ satisfies an $(\alpha,\delta)$-\emph{Weak Poincar\'e Inequality} if 
    \begin{align*}
        \text{Var}_\mu(f) \le \alpha \mathcal E(f,f) + \delta  \Osc(f)^2\,,
    \end{align*}
    for all test functions $f: \Omega \to\mathbb R$, where $\Osc(f) = \max f - \min f$. 
\end{definition}

Note that an $(\alpha,0)$-WPI is exactly an $\alpha$-PI, or equivalently at least a $1/\alpha$-spectral gap. 
The following is the key polynomial WPI that underpins our proof of Theorem \ref{thm:diagonally-Dobrushin-mixing}.

\begin{lemma}\label{lem:diagonally-Dobrushin-WPI}
    Fix $\beta>\beta_c$. For every $C'$, there exists $C$ such that the Ising Glauber dynamics with diagonal-Dobrushin boundary conditions on $\Lambda_n$ satisfies an $(n^{C},n^{-C'})$-weak Poincar\'e inequality. 
\end{lemma}

\begin{figure}[t]
\centering
    \includegraphics[width = 0.25 \textwidth , align = c]{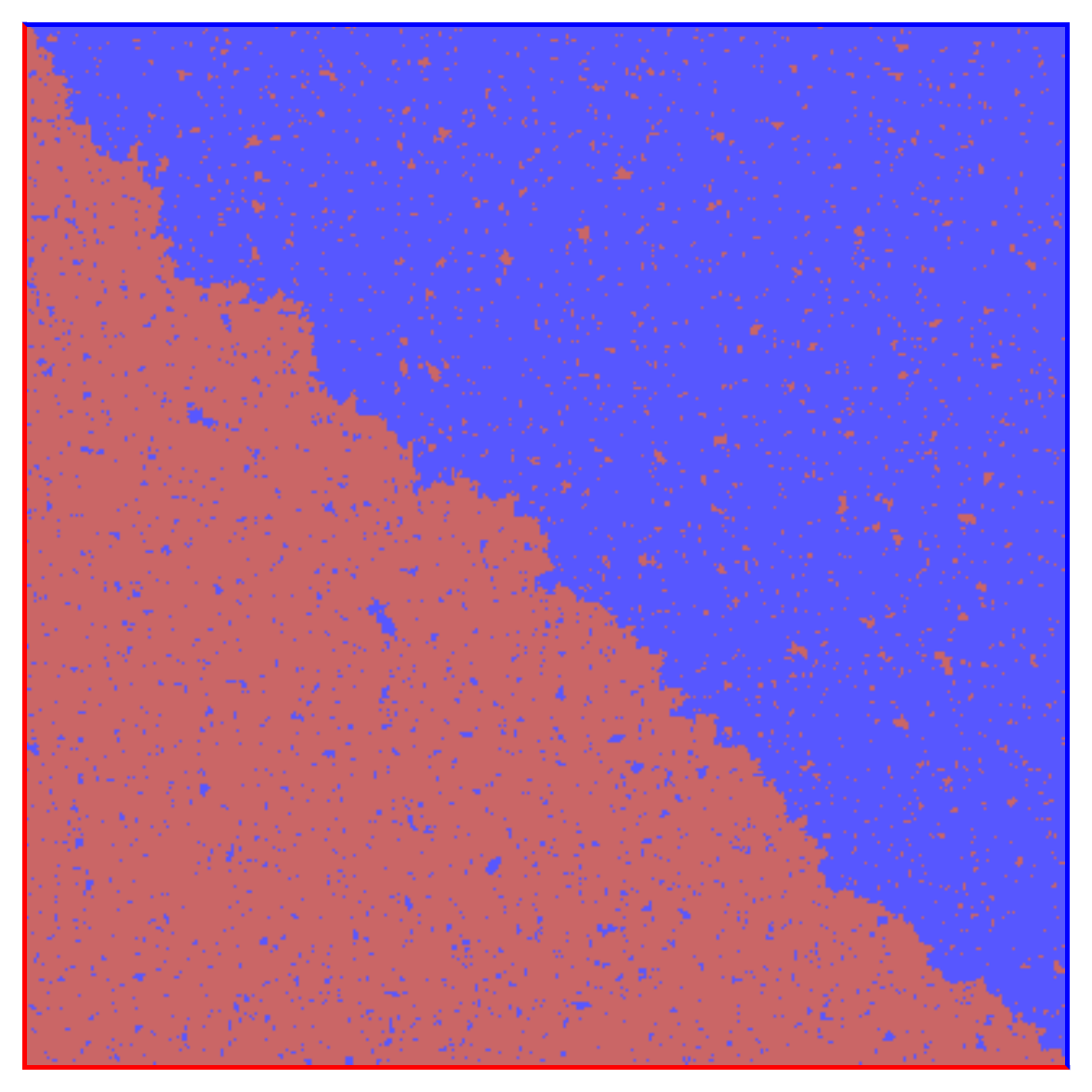} \qquad\qquad 
    \includegraphics[width = 0.55 \textwidth, align= c]{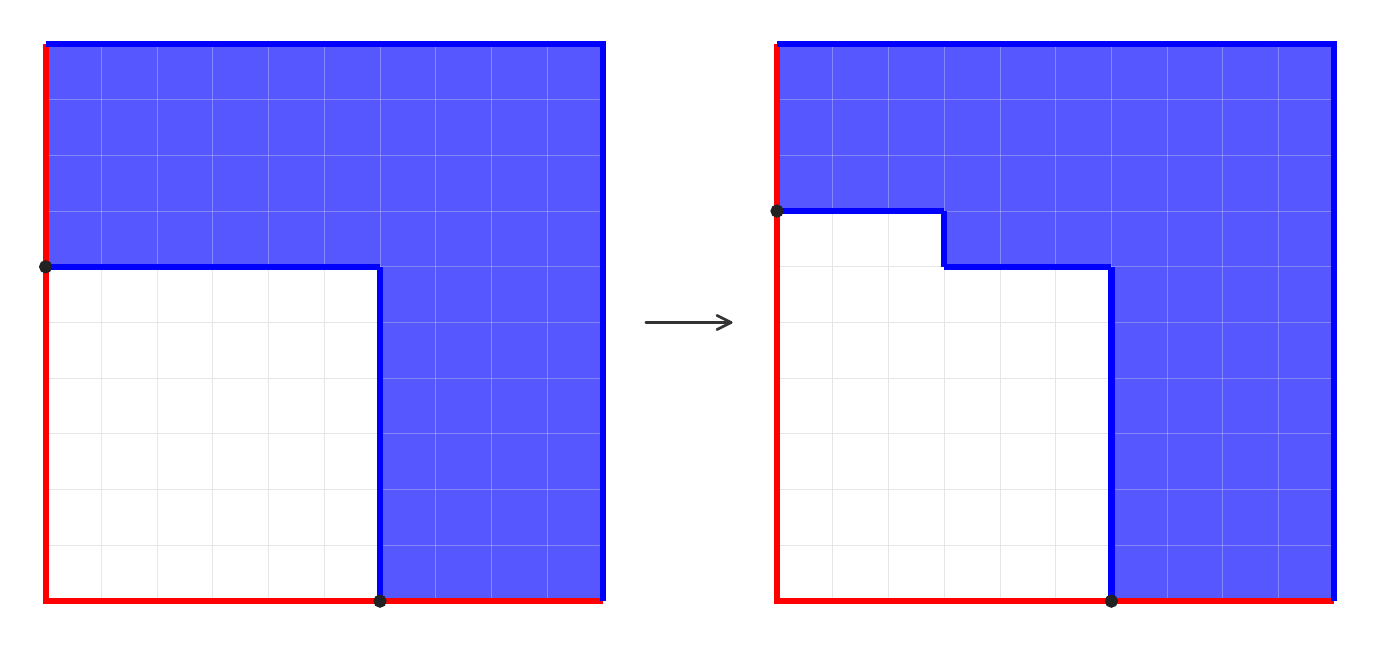} 
\caption{Left: A sample from the Ising model on $\Lambda_n$ for $n=400$ with diagonally Dobrushin boundary conditions. Right: self-reducibility of the diagonally Dobrushin boundary conditions under a censoring scheme that gradually releases vertices from south-west to north-east, one layer at a time.}\label{fig:diag-Dobrushin-censoring}
\end{figure}

Unlike full Poincar\'e inequalities, weak Poincar\'e inequalities do not imply exponential convergence to equilibrium, and in particular, do not control the mixing time. They only imply rapid convergence to equilibrium from initial distributions with Radon--Nikodym to $\pi$ that is bounded by $O(\epsilon^{-1})$. Namely, if a Markov chain $P$ is reversible with respect to $\mu$ and satisfies an $(\alpha,\delta)$-WPI, and $\|\nu/\mu\|_\infty \le B$, then by e.g.,~\cite[Theorem 2.1]{RocknerWang01} and Eq.~(2) of~\cite{WeakPoincareInequalities}, for any $B\sqrt{\delta}< \epsilon<1$, 
\begin{align}\label{eq:WPI-consequence}
    \| P^s\nu  - \mu\|_\tv \le \epsilon \qquad \text{for all }s\ge \alpha \log \frac{B}{\eps^2}
\end{align}
We will informally call this the WPI-time.
The WPI has been a powerful tool in sampling algorithms~\cite{lovasz1993random}, and for simulated annealing schemes in certain high-temperature spin glass settings~\cite{WeakPoincareInequalities,davies2026weakpoincareinequalitiesapproximate,du2026samplingsherringtonkirkpatrickmodelcriticality}. A relevant use-case is that if one has WPI for a sequence of graphs $\Lambda^{(i)}$ that build one-vertex-at-a-time to $\Lambda_n$, then each previous stationary distribution can be used as a $B= O(1)$ warm-start for the next, and if $\delta = o( 1/n^4)$ the error $\varepsilon = o(1/n^2)$ can sustain the union bound over the $n^2$ iterations.  
This strategy which was used in~\cite{lovasz1993random} in the context of sampling from convex bodies, was also recently used in~\cite{RandomFieldIsing} for the random-field Ising model in $(\mathbb Z^d)_{d\ge 2}$. 
For Glauber dynamics in our Ising setting, we can mimic this by using a censoring scheme which allows us to schedule our updates so that the dynamics effectively runs for the WPI time on $\Lambda^{(i)}$ before allowing the next vertex $\Lambda^{(i+1)}\setminus \Lambda^{(i)}$ to also update. 

\begin{definition}\label{def:censoring}
    A \emph{censoring schedule} for Ising Glauber dynamics on $\Lambda$ is a sequence of subsets $(U_i)_{i\ge 1} \subset \Lambda$ and a sequence of increasing times $(t_i)_{i\ge 1}$ such that all clock rings in $\Lambda \setminus U_{i}$ during times $[t_{i-1},t_i)$ are ignored (censored). 
\end{definition}
While more general schedules are permissible, for our purposes, we will simply take the above schedule to be deterministic, fixed in advance, and not depend in any way on any external randomness.

The censoring theorem of~\cite{PWcensoring} then says that if the initialization from which we are performing the censoring scheme is an extremal configuration (all plus or all minus) the mixing time of this process bounds that of the un-censored Glauber dynamics. 
A convenient geometric property of the diagonally Dobrushin boundary conditions is that with a suitably chosen censoring from an extremal configuration, the uncensored region is itself always a  square (or \emph{near-square} of the form  $\{1,...,k\}^2 \cup (\{1,...,j\}\times \{k+1\})$ or $(\{1,...,k\}\times \{1,...,k+1\}) \cup \{k+1\} \times \{j,...,k+1\}$) with diagonally Dobrushin boundary conditions! (See Figure~\ref{fig:diag-Dobrushin-censoring}, right, for how the censoring lets one build out a sequence of near-squares with diagonal Dobrushin boundary conditions.)

In this manner, Lemma~\ref{lem:diagonally-Dobrushin-WPI} together with the censoring scheme  that builds out $\Lambda_n$ using near-squares, one layer at a time, implies the  polynomial spectral gap and mixing time of Theorem~\ref{thm:diagonally-Dobrushin-mixing}. 

\subsection{Benefits of the WPI}
The ability to restrict to proving weak Poincar\'e inequalities essentially means we can restrict our mixing time analysis to a high-probability part of the space. We will use this to restrict to the configurations $\sigma$ that (roughly speaking) have  
\begin{enumerate}
    \item No big bubbles:  Connected components of minus spins above the interface (and plus spins below the interface) have diameter $O(\log n)$. 
    \item ``Regular" interface: All overhangs of the interface have size at most $O(\log n)$. Also, in any vertical strip $J = \{j,...,j+w\}\times \{1,...,n\}$ of width $w$ at least $\Omega(\log n)$, the entry and exit data of the interface through $J$ induce a slope that is in the interval $[-7/4,-1/4]$ around $-1$.
\end{enumerate}
Because all these conditions (and some more technical ones related to the Ornstein--Zernike representation of the interface) have their thresholds set at  $C\log n$ for a large constant $C$, at equilibrium they can be shown to have  probability $1-n^{-10}$. Therefore, to establish Lemma~\ref{lem:diagonally-Dobrushin-WPI}, it will suffice to restrict Glauber dynamics to this nice part of the state space and prove an $n^C$ Poincar\'e inequality: see Lemma~\ref{lem:conditioning-wpi}.  
Still, even with this restriction to nice interfaces and configurations, proving polynomial mixing of the interface remains a challenge for the reasons spelled out earlier that the microscopic defects of the interface preclude the use of the method of~\cite{Wilson}.   

Let us first reduce the task to proving polynomial WPI-time for a chain on interfaces as that will be the crux of the argument. 
If we let $\bar \pi$ be the marginal of $\Gamma$, the law of total variance gives 
\begin{align}\label{eq:intro-total-variance}
    \text{Var}_\pi (f)= \bar \pi \big[\text{Var}_\pi(f \mid \Gamma)] + \text{Var}_{\bar \pi} ( \pi[f\mid \Gamma])\,.
\end{align}
Reminiscent of the projection-restriction approach~\cite{madras2002markov,JSTV04}, the decomposition allows us to bound these individually by proving (weak) Poincar\'e inequalities for a projection chain which makes Glauber-type resampling updates to the interface $\Gamma$ using its marginal distribution, and a Glauber chain that updates the configuration conditional on the interface $\Gamma$. 

\subsection{Mixing off of the interface}
The (significantly) simpler of these two to control is the weak Poincar\'e inequality for Glauber dynamics conditional on an interface $\Gamma$. The law of the Ising model given $\Gamma$ tensorizes into a product of a $+1$-boundary conditions model on the subset of $\Lambda$ above $\Gamma$ and $-1$-boundary model on the subset of $\Lambda$ below $\Gamma$. Since we restricted to bubbles being of size $C\log n$ at the expense of only proving a WPI, the exponential decay of correlations of the Ising model at $\beta>\beta_c$ with monochromatic boundary conditions allows us to bound the mixing time by the mixing time on $C' \log n$-width strips for $C' \gg C$. This will be at most a polynomial by a standard bound of the mixing time of a spin system  by exponential in the cutwidth of the graph: see Lemma~\ref{lem:cutwidth-mixing-bound}.   

\begin{figure}[t]
\centering
    \includegraphics[width = 0.29 \textwidth]{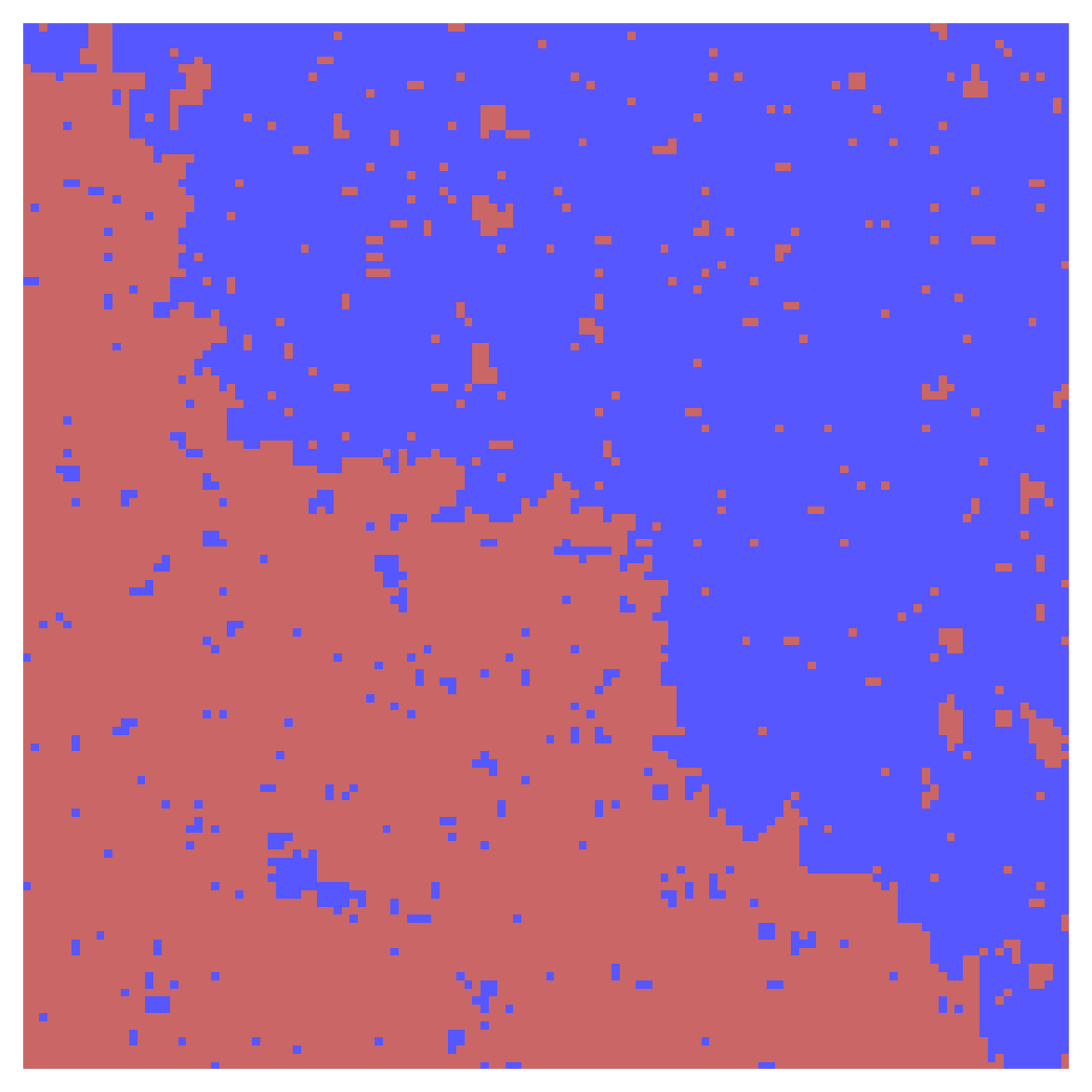} \qquad\qquad\qquad 
    \includegraphics[width = 0.29 \textwidth]{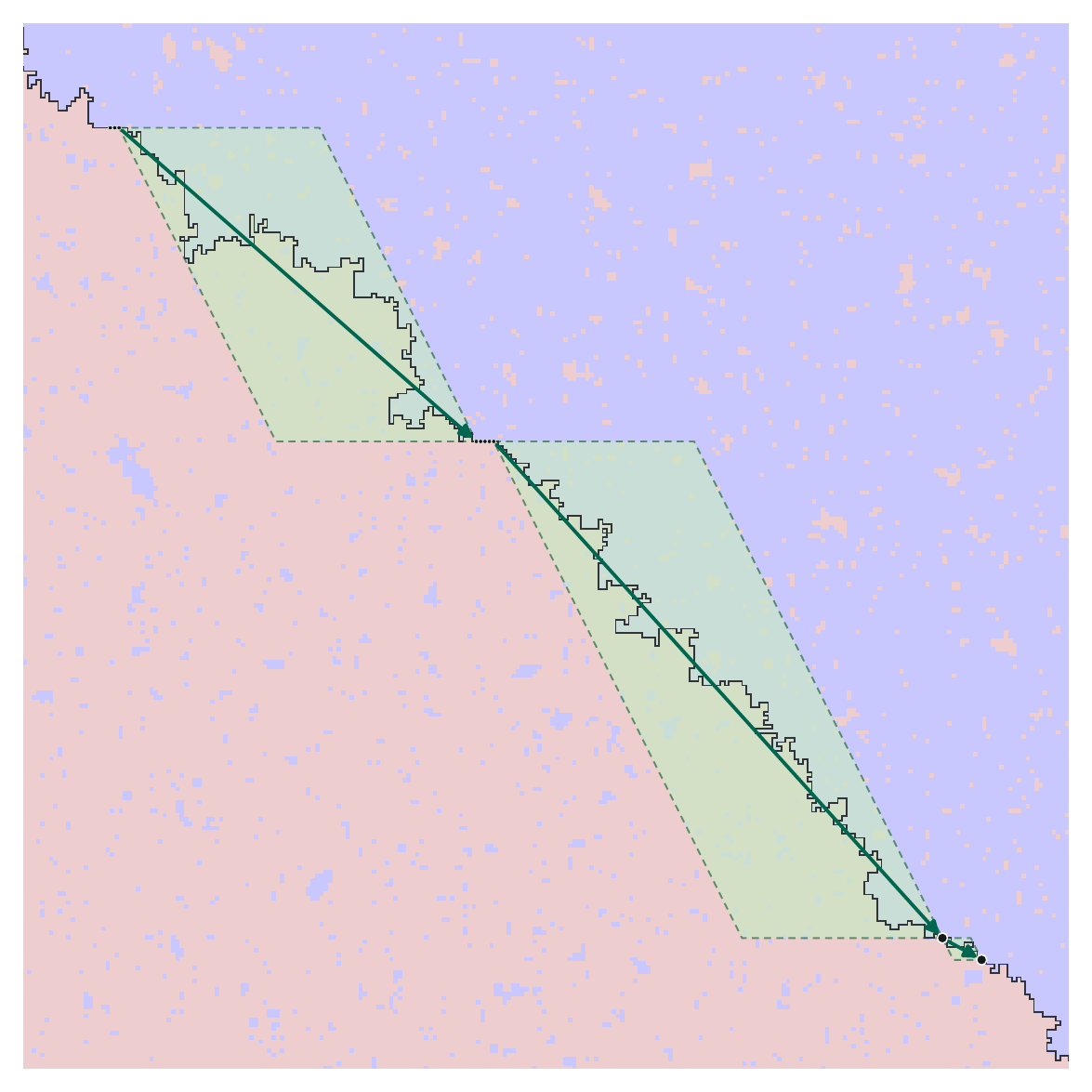} 
\caption{Left: A zoom-in of the sample from the Ising model with diagonal Dobrushin boundary conditions. Right: The black line traces
the interface between the two plus and minus phases. The pale green cones show its Ornstein--Zernike
decomposition.}\label{fig:diag-interface-OZ-decomp}
\end{figure}

\subsection{Mixing of regular Ornstein--Zernike interfaces}\label{sec:proof-sketch-interface-mixing}

We now turn to proving a WPI for the marginal on the interface, corresponding to showing that the projection of the chain to $\Gamma (X_t)$ equilibrates quickly. We begin by reiterating the point that the approaches based on~\cite{Wilson} that were fruitful for approximating models, do not apply because the expected height of an Ising interface at the midpoint of a segment of length $m$ does not exactly satisfy the harmonicity property needed to employ the method---rather we can only guarantee it will be within $O(\sqrt{m})$ of the average. In fact, even an $O(1)$ error here in the expected height at the midpoint would destroy the applicability of the Laplacian eigenfunction.

In an idealized formulation, the prior bounds on the plus boundary Ising mixing time of~\cite{MaTo,LMST} can be described as using a local limit estimate of $1/\sqrt{m}$ that the interface passes exactly through the midpoint of its boundary data, to split the analysis of the interface dynamics into an independent left-half and right-half, both with boundary data of the same slope. (Technically, the authors modified the boundary conditions to make it typical for the interface to go through the midpoint, at a polynomial cost in each step of the recursion.)  Repeating this analysis $O(\log n)$ times to go from segments of width $n$ to segments of width $O(1)$ costs a polynomial at each of the $O(\log n)$ scales, giving the quasipolynomial bound of $e^{ O(\log n)^2}$ that they arrive at. 

To get polynomial bounds, one cannot afford to lose more than a constant per scale. To ensure this we take advantage of the Brownian-like property of the interface at any scale large enough to homogenize the effect of the microscopic defects. Namely, we use the diffusive fluctuations of the interface to couple a pair of interfaces on a segment of length $m$ with boundary data that are separated by $O(\sqrt{m})$ at their midpoint column with constant probability. That allows us to split the analysis into domains of a constant factor smaller width with only a constant cost, but at the expense of having to allow slopes that deviate from the initial $-1$ slope. 
This constant cost, when raised to a $\log n$ power for recursing down in scale from width $n$ to width $O(\log n)$, yields the polynomial WPI bound we end up with. 

To be a bit more precise, consider an idealization of the Ising interface whereby it is a random-walk bridge on a segment of width $m$ and of slope $-1$, and suppose we only have access to crude Gaussian upper bounds on its height fluctuations, and anti-concentration bounds at resolution scales of $C \log n$ and larger. Consider the block dynamics that performs heat-bath resampling updates on this interface with two blocks $B_{\rm West}, B_{\rm East}$ which are the left two-thirds and right two-thirds of the segment, respectively. Starting from two initial  interfaces $\Gamma(0), \Gamma'(0)$, we couple them with $\Omega(1)$-spectral gap by the following two stage scheme: let $\Delta_{\rm West}(t), \Delta_{\rm East}(t)$ be the height deviations of $\Gamma(t)$ from the line connecting its boundary points, at indices $m/3, 2m/3$ respectively. 
\begin{enumerate}
    \item \emph{Reducing the height to diffusive scale.} In an update on $B_{\rm West}$, the deviation $\Delta_{\rm West}(t)$ will in expectation become $\frac{\Delta_{\rm East}(t-1)}{2} + O(\sqrt{m})$ (the error being because we only have access to Gaussian tail bounds around the midpoint), and likewise under an East update, $\Delta_{\rm East}(t)$ will become $\frac{\Delta_{\rm West}(t-1)}{2} + O(\sqrt{m})$. Therefore, $\Delta_{\rm West}, \Delta_{\rm East}$ contract to $O(\sqrt{m})$ scale at $\Omega(1)$ rate. 
    \item  \emph{Coupling diffusively separated interfaces.} Once $\Delta_{\rm East},\Delta_{\rm East}'$ are both $O(\sqrt{m})$, in an update on $B_{\rm West}$, the anticoncentration input ensures that with constant probability they can be coupled at index $m/3$. If that coupling succeeds, by the Markov property, they can be coupled to the left of $m/3$, and in the next update on $B_{\rm East}$, on the entire segment. 
\end{enumerate}

\begin{figure}[t]
\centering
    \includegraphics[width =  \textwidth]{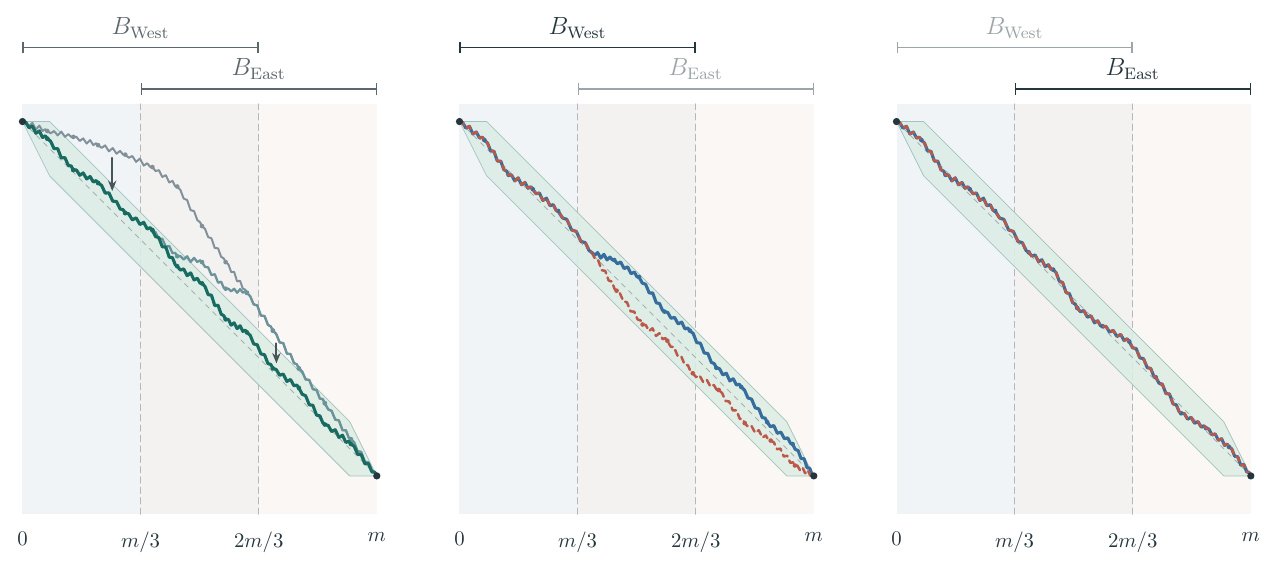} 
\caption{Left: First stage of the interface dynamics brings the interface towards the line connecting boundary points by a constant factor in expectation in each step. Middle and Right: Once two interface heights are $O(\sqrt{m})$ apart on a width $m$ domain, with constant probability when $B_{\rm West}$ updates, they are coupled in the left third, and the next update on $B_{\rm East}$ makes them agree everywhere.}\label{fig:interface-block-dynamics}
\end{figure}

\noindent See Figure~\ref{fig:interface-block-dynamics} for a visualization of these two stages to bound the coupling time of the block dynamics.  

We conclude by describing a few of the technical challenges to implementing this strategy for the  actual interface $\Gamma$ of the Ising model. The Ising interface $\Gamma$ has an \emph{effective random-walk representation} known as the \emph{Ornstein--Zernike (OZ)} representation, which decomposes the interface at \emph{cone-points} (depicted in Figure~\ref{fig:diag-interface-OZ-decomp}, right) into complicated increments between these cone-points. Using this representation---in fact a decorated version of it from~\cite{IOSV-prewetting} which gives the interface a true renewal structure---the jumps between the renewal points of $\Gamma$ obey the law of a \emph{directed two-dimensional random-walk bridge}. In particular, the number of increments of this bridge is random, and heat-bath resampling moves for the interface are only well-behaved if they only resample between renewal points. Still, at only polynomial cost, we can work with a block dynamics that resamples from the nearest renewal point to the edge of the block it is resampling. 

The inputs from Ornstein--Zernike theory for the low-temperature Ising model are sufficiently well-developed that though technical, the idealized argument of 1--2 above essentially works so long as we work with segments of width at
least $C\log n$ and restrict to interfaces that are sufficiently nice. Sufficiently nice here means the increments between renewal points are all $O(\log n)$, and all slopes encountered in block resamplings are uniformly bounded (these are both events of probability $1-n^{-10}$). Finally, paying a polynomial for the mixing time of the $O(\log n)$-width strip, we can obtain the desired polynomial weak Poincar\'e inequality for the Ising interface $\Gamma$ with respect to the Glauber dynamics' Dirichlet form.\\

While implementing the above strategy is the heart of the paper, before embarking on that process, we conclude the proof of Theorem~\ref{thm:plus-bdy-mixing} assuming Theorem~\ref{thm:diagonally-Dobrushin-mixing}. The remainder of the paper will then focus entirely on the proof of the latter.

\subsection{Reducing all-plus boundary to the diagonal Dobrushin boundary}\label{subsec:Dobrushin-to-all-plus}

\begin{proof}[{Proof of Theorem~\ref{thm:plus-bdy-mixing}}]
Let $\Lambda_n$ have all-plus boundary conditions. 
A convenient consequence of the monotonicity of the Ising Glauber dynamics is that the mixing time from worst-case initialization is controlled by the mixing time from the two-extremal initializations, all-plus and all-minus.  Namely, for any domain $\Lambda$ and boundary conditions $\tau$, if $X_{+,t},X_{-,t}$ denote  the Ising Glauber dynamics on $\Lambda$ with $\tau$ boundary, initialized from all-plus and all-minus respectively, one has (using the grand coupling and monotonicity) 
    \begin{align}\label{eq:tv-controlled-by-extremal-inits}
    \max_{\sigma_0}\|\mathbb P(X_t \in \cdot \mid X_0 = \sigma_0) - \pi_{\Lambda}^{\tau} \|_\tv & \le \mathbb P(X_{+,t} \ne X_{-,t}) \nonumber \\ & \le \sum_{v\in \Lambda} (\mathbb P(X_{+,t}(v) =1) - \mathbb P(X_{-,t}(v)=1)) \nonumber\\ 
    & \le 2 \max_{\iota\in \{-,+\}} \sum_{v \in \Lambda}  |\mathbb P(X_{\iota,t} (v) =1) - \pi_{\Lambda}^\tau(\sigma_v =1)| \,.\end{align}
We apply the above with $\Lambda = \Lambda_n$ and with $\tau\equiv +1$. 
It therefore suffices to show that for $t \ge n^{C}$ for suitably chosen $C$, each of $\mathbb P(X_{+,t}(v)=1)$ and $\mathbb P(X_{-,t}(v) =1)$ are within $o(1/n^2)$ of $\pi_{\Lambda_n}^+(\sigma_v =1)$. 

    The fast mixing from the all-plus initialization is straightforward (e.g., the mixing time bounds of~\cite{LMST} can be applied at the logarithmic scale using monotonicity as done in the proof of Proposition~3.3 in~\cite{GhSi22}). For self-containedness, one could also apply our Lemma~\ref{lem:plus-init-mixing-with-plus-bdy} to see that for $t\ge n^{2C}$, one has $|\mathbb P(X_{+,t}(v)=1) - \pi_{\Lambda_n}^+(\sigma_v =1)| \le n^{ - 10}$. 

    We now focus on the all-minus initialization. By rotation symmetry of $\Lambda_n$, it suffices to bound  $|\mathbb P(X_{-,t}(v)  = 1) - \pi_{\Lambda_n}^+(\sigma_v=1)|$ for $v$  in the top-right quadrant of $\Lambda_n$. Fix any such vertex $v$. Shift $\Lambda_n$ to be the subset $\{3n, ...,4n\}^2$ of the domain $D_n = \{1,...,4n\}^2$ and endow $\Lambda_n$ with all-plus boundary, while $D_n$ has the diagonally Dobrushin boundary. Let $Y_{-,t}$ be the Ising Glauber dynamics on $D_n$ initialized from all minus. By monotonicity, 
    \begin{align*}
        \mathbb P(Y_{-,t}(v) =1)  \le \mathbb P(X_{-,t}(v) =1) \le \pi_{\Lambda_n}^+(\sigma_v = 1)\,.
    \end{align*}
    By Theorem~\ref{thm:diagonally-Dobrushin-mixing} (and the standard sub-multiplicativity of total-variation distance after the mixing time), there exists a constant $C= C(\beta)>0$ such that for all $t\ge n^{C}$, 
    \begin{align*}
       0 \le \pi_{D_n}^{\rm Dob}(\sigma_v =1) -  \mathbb P (Y_{-,t}(v) =1) \le n^{-10}\,.
    \end{align*}
    Together the above two inequalities give us the bound that for all $t\ge n^C$ 
    \begin{align*}
        |\mathbb P(X_{-,t}(v) = 1)   - \pi_{\Lambda_n}^+(\sigma_v =1)|\le n^{-10} + |\pi_{\Lambda_n}^+(\sigma_v =1)- \pi_{D_n}^{\rm Dob}(\sigma_v =1)|\,.
    \end{align*}
    The proof is done after the following Lemma~\ref{lem:embed-plus-bdy-in-Dobrushin-stationary} which bounds the difference in stationary measures between $\Lambda_n$ with plus boundary and the containing $D_n$ with its diagonally Dobrushin boundary. 
\end{proof}

\begin{lemma}\label{lem:embed-plus-bdy-in-Dobrushin-stationary}
    Fix $\beta>\beta_c$. Consider $\Lambda_n = \{3n, ...,4n\}^2$ with all-plus boundary conditions and $D_n = \{1,...,4n\}^2$ with diagonally Dobrushin boundary conditions. For all $v\in \{(7/2)n,...,4n\}^2$, 
    \begin{align*}
        \pi_{\Lambda_n}^+ (\sigma_v =1) - \pi_{D_n}^{\rm Dob} (\sigma_v=1) \le e^{ - \Omega(n)}\,.
    \end{align*}
\end{lemma}

Lemma~\ref{lem:embed-plus-bdy-in-Dobrushin-stationary} essentially follows from the fact that the interface in $D_n$ has Gaussian concentration around the diagonal, and is therefore typically order-$n$ distance from $v$. Since the south-west boundary of $\Lambda_n$ is also order-$n$ distance away from $v$ and both induce plus boundary conditions, decay of truncated correlations of the low-temperature Ising model couple them at $v$.

\begin{proof}[Proof of Lemma~\ref{lem:embed-plus-bdy-in-Dobrushin-stationary}]
 Let $\Gamma$ be the interface in $D_n$ with its diagonally Dobrushin boundary conditions. 
Consider the event $A(\Gamma)$ that the interface $\Gamma$ is below the line  $y = 5n-x$, i.e., $D_n$'s boundary line shifted up by $n$. If we condition on the interface $\Gamma$, then 
\begin{align}\label{eq:single-site-split-+-dob}
    \pi_{\Lambda_n}^+(\sigma_v =1 ) - \mathbb E_{\Gamma \sim \pi_{D_n}^{\rm Dob}} [\pi_{D_n}^{\rm Dob}(\sigma_v =1 \mid \Gamma )] \le \big(\pi_{\Lambda_n}^+(\sigma_v =1) - \min_{\Gamma \in A} \pi_{D_n^\uparrow (\Gamma)}^+ (\sigma_v =1)\big) + \pi_{D_n}^{\rm Dob}(A^c)\,,
\end{align}
where $D_n^\uparrow(\Gamma)$ is the part of $D_n$ ``above" the interface, that is all connected components of $[1,4n]^2 \setminus \Gamma$ on which the knowledge of $\Gamma$ induces the plus boundary conditions (see Section~\ref{sec:domains}). 

For the quantity in the parenthetical of~\eqref{eq:single-site-split-+-dob}, fix $\Gamma\in A$. 
Since $D_n^\uparrow (\Gamma)$ and $\Lambda_n$ share their north and east plus boundaries, the set $S = \partial \Lambda_n \cap D_n^\uparrow (\Gamma)$ is of size $2n$ and $d(v,S) \ge n/2$. We will use the truncated decay of correlations under $\pi_{\mathbb Z^2}^+$.  
Enumerate $S=\{x_1,\ldots,x_{|S|}\}$ and set
\[
    \mu_0=\pi_{D_n^\uparrow(\Gamma)}^+,\qquad
    \mu_j=\mu_0\bigl(\,\cdot\mid
        \sigma_{x_1}=\cdots=\sigma_{x_j}=1\bigr).
\]
By the domain Markov property, the marginal of $\mu_{|S|}$
on $\Lambda_n$ is $\pi_{\Lambda_n}^+$.
Each $\mu_j$ is an Ising measure with plus boundary and
additional plus pinnings. The GHS inequality and exponential
decay of truncated correlations in the infinite-volume
plus phase at $\beta>\beta_c$ give, uniformly in $\Gamma\in A$ and $j$,
\begin{equation}\label{eq:ghs-pinned-covariance}
    0\le \operatorname{Cov}_{\mu_j}(\sigma_x,\sigma_y)
    \le \operatorname{Cov}_{\pi_{\mathbb Z^2}^+}
                         (\sigma_x,\sigma_y)
    \le  e^{- \Omega( d(x,y))};
\end{equation}
see \cite[Section~2.3.1, Eqs.~(2.32)--(2.33)]{BCK03}.
Since $\mu_{j-1}(\sigma_{x_j}=1)\ge 1/2$, we get 
\begin{align}
    \pi_{\Lambda_n}^+(\sigma_v=1)-\pi_{D_n^\uparrow(\Gamma)}^+(\sigma_v=1)
    &=
    \sum_{j=1}^{|S|}
    \frac{\operatorname{Cov}_{\mu_{j-1}}
                    (\sigma_v,\sigma_{x_j})}
         {4\mu_{j-1}(\sigma_{x_j}=1)}
     \le  e^{-\Omega(n)}\,,
    \label{eq:plus-domain-comparison}
\end{align}
uniformly over $\Gamma \in A.$ 

The term $\pi_{D_n}^{\rm Dob}(A^c)$ in~\eqref{eq:single-site-split-+-dob} is bounded by $e^{-\Omega(n)}$ by the Gaussian one-point tails of Ising interface heights and a union bound: see e.g.~\cite[Eqs.~(5.4)--(5.9)]{LMST} for a derivation in the case of horizontal Dobrushin boundary conditions, that carries over mutatis mutandis to the $\pi/4$-diagonal Dobrushin boundary conditions due to the sharp triangle inequality for the surface tension holding at all angles. Plugging this, together with~\eqref{eq:plus-domain-comparison}, in to~\eqref{eq:single-site-split-+-dob} completes the proof of the lemma. 
\end{proof}

 \subsection{AI Disclosure}
This project originated in 2023 when the authors, with an eye towards the problem of the Ising dynamics with all-plus boundary conditions, began working on a proof of polynomial mixing for Ornstein--Zernike interfaces (without dealing with the bubbles above and below the interface). The main aim of the project was to develop methods to prove polynomial mixing of random-walk like bridges, without assuming the microscopic resampling properties or relations to Laplacian eigenfunctions that facilitate the Wilson method. By December 2025 the authors had a 43 page draft of a paper proving polynomial mixing time for disordered and finitely dependent SOS models, and more generally random walk bridges only assuming tail bounds on increments, and mesoscopic diffusivity and anticoncentration properties (nothing at the microscopic level). This contained the argument described in Section~\ref{sec:proof-sketch-interface-mixing}, and is one of the main ingredients in this paper for the polynomial mixing for Ising interfaces themselves. However, given that dealing with the latter still would involve many technical challenges, the draft manuscript was temporarily tabled. 

In October 2026, recognizing the capabilities of ChatGPT 6.0 Astra, the authors uploaded their draft document to Astra, and guided it through implementing their already-developed proof strategy for Ising type interfaces via intermediate models of growing complexity, starting from a fixed number of OZ increments, to OZ bridges with random number of increments (but still no bubbles), to Ising itself with the further suggestion of ensuring regularity using monotonicity.
While the idea of restricting to high probability parts of the space to ensure regularity was already present in the authors' initial draft (and also in earlier works like \cite{MaSi-SOS}), Astra identified that this restriction is particularly effectively handled through proving weak Poincar\'e inequalities (WPI). Namely, with diagonal Dobrushin boundary conditions, starting from the extremal configurations allows one to employ the censoring scheme outlined in Section~\ref{sec:WPI-to-gap}, and take advantage of the  self-similar nature of the domains along with their boundary conditions. This reduces the problem to proving a polynomial WPI for the Dobrushin domain, whence regularity on the interface can be assumed. While ideas of gradually growing domains after mixing on a smaller domain through WPI have been around in the sampling literature~\cite{lovasz1993random}, combining it with censoring to get a full Poincar\'e inequality seems to be novel, and we expect it will find future applications.

The entire body of the paper is human written. The technical Appendix~\ref{app:OZ-inputs} proving minor modifications of known Ornstein--Zernike inputs to our specific setting was written with the assistance of ChatGPT 6.0 Astra. ChatGPT 6.2 Astra was also used for assistance with generating figures and in the final proofreading and copy-editing. 

\subsection{Outline of paper}
Section~\ref{sec:contours} describes the low-temperature (contour) representation of the Ising model, the geometry of its interface, and its Ornstein--Zernike representation. In Section~\ref{sec:WPI-off-interface}, we prove the polynomial WPI conditional on the interface. The heart of the paper is Section~\ref{sec:interface-wpi} which proves the polynomial WPI for the projection of the dynamics to the interface. Section~\ref{sec:mixing-Dobrushin-bdy} deduces Lemma~\ref{lem:diagonally-Dobrushin-WPI} and proves Theorem~\ref{thm:diagonally-Dobrushin-mixing}. Finally, Appendix~\ref{appendix:OZ-estimates} proves deferred technical Ornstein--Zernike  inputs. 

\subsection{Acknowledgments}
R.G.\ is grateful to A.\ Sinclair for extensive discussions on this problem in 2019, in which some of the ideas to smoothen the interface dynamics using $O(\log n)$-width strips were already present. 
S.G.'s research was partly supported by NSF grants 1945172 and 2553225 and a grant from the Miller Institute. 
The research of R.G.\ is supported in part by NSF CAREER grant DMS-2440509, NSF DMS grant 2246780, and a Simons Fellowship in Mathematics. Part of this work was done while both S.G.\ and R.G.\ were long-term participants in the program `Spectral Theory beyond Graphs' at the Simons Institute for the Theory of Computing. 
This work was supported in part by a grant of access to OpenAI models through the ChatGPT for Academic Researchers program.

\section{Low-temperature representation and interfaces}\label{sec:contours}
Throughout the paper, we imagine $\beta>\beta_c$ fixed, allow constants $c,C,c_\beta,C_\beta$ to depend on $\beta$, and suppose $n$ is sufficiently large. In particular, unless otherwise indicated, big-$O$, $\Theta,\Omega$ notation is as $n\to\infty$ and can hide constants that depend on $\beta$.

\subsection{Contour decomposition}\label{sec:domains}
{\color{green!50!black}}
In the remainder of the paper, we work with the domain~$\Lambda_n=\{1,\ldots,n\}^2$ with diagonally Dobrushin boundary conditions (those that are minus on the south and west sides of $\Lambda_n$ and plus on the north and east sides) and abbreviate $\pi= \pi_{\Lambda_n}^{\rm Dob}$. 
For an Ising configuration $\sigma\in \{-1,+1\}^{\Lambda_n}$, associate a set of edges in the dual graph $F(\sigma) \subset E({\Lambda_n^*})$ where a dual edge $e^*$ separating vertices $u\sim v$ is in the configuration if and only if $\sigma_u\ne \sigma_v$. Since the boundary conditions change sign exactly twice as one goes around $\Lambda_n$, once at $a_\Lambda = (0,n)$ and again at $b_\Lambda = (n,0)$, the dual graph $F(\sigma) \subset E(\Lambda^*)$ we get will have even degree at all internal vertices of $\Lambda_n^*$ and a source and sink dual vertex of odd degree at $a^*_\Lambda = (\frac{1}{2},n-\frac{1}{2})$ and $b^*_{\Lambda} = (n+\frac{1}{2}, -\frac{1}{2})$. 

The Ising configurations on $\Lambda_n$ are easily seen to be in weight-preserving bijections with the dual subgraphs $F$ if we assign the dual subgraphs weight $$ \pi_{\Lambda}^{\rm Dob}(F) \propto e^{-2\beta |F|}\,.$$

\begin{figure}
\centering
    \includegraphics[width = 0.3\textwidth, align= c]{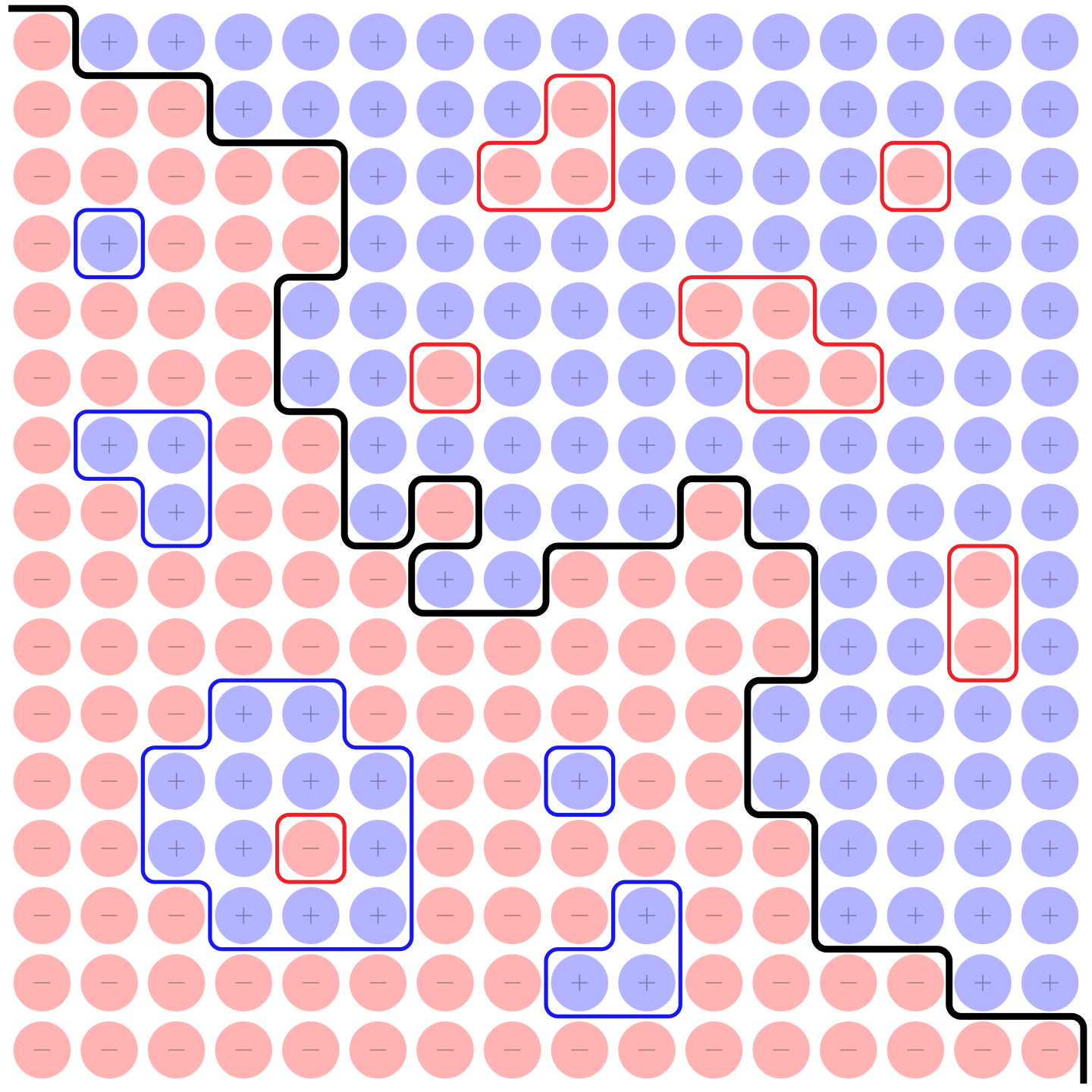} \qquad  \qquad\qquad  
    \includegraphics[width =0.2\textwidth, align= c]{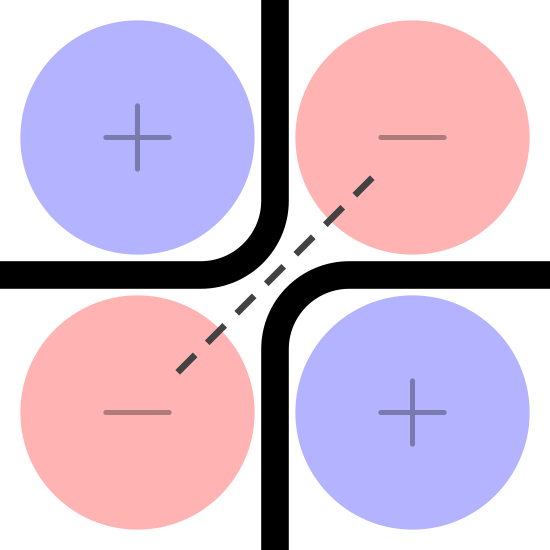}
    \caption{Left: An Ising configuration $\sigma$ on $\Lambda_n$ with diagonally Dobrushin boundary conditions. The dual edges shown are form $F(\sigma)$ and the interface $\Gamma$ is shown in black. Right: Zoom of a vertex near the middle shows the NW/SE splitting when four edges meet.}\label{fig:contour-decomposition}
\end{figure}

\begin{definition}\label{def:contour-decomposition}
For a subset $F\subset E({(\mathbb Z^2)^*})$, we wish to decompose it into a sequence of paths between the vertices incident to an odd number of edges of $F$, and a collection of loops. A priori, such a decomposition is ill-defined due to vertices incident to three or four edges of $F$. To remedy this ambiguity, it is standard to use a \emph{splitting convention} according to which we partition the edges of $F$. For concreteness we use the south-east ($SE$) splitting rule and partition $F$ into loops and paths as follows: for every vertex $v\in (\mathbb Z^2)^*$ incident to three or four edges, if two of those edges are the edges directly south and east of $v$, \emph{split} that pair apart from the remaining edge(s), and analogously if two of the edges are those directly north and west of $v$, split that pair apart from the remaining edge(s). Partition $F$ into $\{\gamma_1 ,\gamma_2 ,...,\gamma_k\}$ by taking the (maximal) connected components of $F$, where edges are connected to one another if they share a vertex \emph{and have not been split apart at that vertex}. The resulting partition $\underline \gamma := \{\gamma_1,\gamma_2,...,\gamma_k\}$ consists only of edge-disjoint paths and loops, and is called the \emph{contour decomposition}. See Figure~\ref{fig:contour-decomposition} for a depiction. 
\end{definition}

\begin{definition}\label{def:interface}
With the splitting of Definition~\ref{def:contour-decomposition}, for any configuration $\sigma$ on $\Lambda_n$, we can define the interface $\Gamma$ as the unique contour in the decomposition of $F(\sigma)$ that connects $a^*_{\Lambda} \to b^*_{\Lambda}$. 
\end{definition}

For a subgraph $F \subset E(\Lambda^*)$, we use $\partial F$ to denote the dual vertices that have odd degree in $F$. In what follows, an admissible path from dual vertices $a\to b$ in a domain $\Lambda^*$ will mean a contour $\gamma \subset \Lambda^*$ having $\partial \gamma = \{a,b\}$. 

Knowledge that some admissible $\Gamma$ is in the contour decomposition of $\sigma$ implies that the edges of $\Gamma$ are in $F$, but due to our splitting rule, it also forces certain edges incident to it to not be in $F$. For instance, a path passing through a vertex \emph{reserves} both edges of each grouping $\{N,W\}$ and $\{S,E\}$ it uses, because our splitting rule means if the interface goes from $\{(-1,0), (0,0)\}\to \{(0,0),(0,-1)\}$ then necessarily the $N$ and $E$ edges can't be in $F(\sigma)$ because if they were they would have been part of $\Gamma$. We let $Q_\Gamma$ be the union of the edges of $\Gamma$ with the edges it ``reserves" and note that crudely it has size at most $3|\Gamma|$ where   $|\Gamma|$ is the number of edges in $\Gamma$.

Then by the bijection of Definition~\ref{def:contour-decomposition}, in the domain $\Lambda$ with its diagonal Dobrushin boundary conditions, given interface $\Gamma$, the law of the remainder of $F(\sigma)$ has weights proportional to $e^{-2\beta |F|}$ and partition function  
\begin{equation}\label{eq:contour-length-deletion}
 Z_0(H)=\sum_{F\subset E(H):\,\partial F=\varnothing}e^{-2\beta|F|}\,, \qquad H = E(\Lambda^*) \setminus Q_\Gamma\,.
\end{equation}
I.e., it is a weighted distribution over even subgraphs on the graph $\Lambda^* \setminus Q_\Gamma$, which is the same as an Ising model on $H^*$ with monochromatic boundary conditions.

\begin{definition}\label{def:collar}
Let~$\partial_v\Gamma\subset\Lambda$ be the endpoints of the primal edges dual to~$Q_\Gamma$. The fact that $\Gamma$ is an interface prescribes the spin values on exactly these bonds to be $+$ ``above" $\Gamma$ and $-$ ``below" $\Gamma$.  Namely, starting from the fixed boundary spins and filling in the configuration $\sigma = F^{-1}(\Gamma)$ determines an assignment of spins~$s_\Gamma$ on~$\partial_v\Gamma$. It is easy to see that 
\begin{equation}\label{eq:collar-size}
 |\partial_v\Gamma|\le 8|\Gamma|+4,
 \qquad
 \partial_v\Gamma\subset\{x\in\Lambda:\operatorname{dist}_\infty(x,\Gamma)\le2\}\,.
\end{equation}
Let $\Lambda^\uparrow _n(\Gamma)$ be the union of  connected components of $\Lambda\setminus \partial_v \Gamma$ that get $+1$-boundary under $s_\Gamma$, let $\Lambda^{\downarrow}_n(\Gamma)$ be the union of those that get $-1$ boundary under $s_\Gamma$. 
\end{definition}

\begin{observation}\label{obs:conditional-spin-law}
For every admissible interface~$\Gamma$, one has 
    $\pi(\cdot \mid \Gamma) = \pi_{\Lambda^\uparrow (\Gamma)}^+\otimes \pi_{\Lambda^\downarrow (\Gamma)}^-$. 
\end{observation}

Since each of the connected components of $\Lambda \setminus \partial_v \Gamma$ get the monochromatic boundary conditions, because of truncated decay of correlations at $\beta>\beta_c$, given $\Gamma$, closed contours are exponentially unlikely in number of edges. Namely, we will use the following precise equilibrium input in arbitrary domains with monochromatic boundary conditions. For every $\beta>\beta_c$, there exist constants $C_\beta,m_\beta$ such that for every finite~$A\subset\mathbb Z^2$,  and~$x,y\in A$,
\begin{equation}\label{eq:loops-covariance-input}
 0\le\cov_{\pi_{A}^+}(\sigma_x,\sigma_y)
 \le\cov_{\pi_{\mathbb Z^2}^+}(\sigma_x,\sigma_y)
 \le C_\beta e^{-m_\beta\|x-y\|_\infty}\,.
\end{equation}
The first inequality is FKG, the second is the GHS inequality, and the last estimate is exponential decay of truncated correlations in the two-dimensional plus phase. The bounds in~\eqref{eq:loops-covariance-input} are recorded e.g.,~in~\cite[Section 2.3.1, (2.32)--(2.33)]{BCK03}.

\subsection{Random line representation}
We refer the reader e.g., to the survey~\cite[Sections 3.1--3.2]{OttVelenikSurvey} for more details on what follows. 
Recalling~\eqref{eq:contour-length-deletion}, by  a spin flip on $\Lambda_n^\downarrow(\Gamma)$, the probability of the interface being $\Gamma$ can be expressed using the following \emph{random-line weights}. For any $\Lambda$, and admissible path $\gamma\subset \Lambda^*$, 
\begin{equation}\label{eq:q-Lambda}
 q_\Lambda(\gamma)=e^{-2\beta|\gamma|}\frac{Z_0(\Lambda^*\setminus Q_\gamma)}{Z_0(\Lambda^*)}\,.
\end{equation}

Abusing notation slightly, one may directly use a dual domain $H\subset E((\mathbb Z^2)^*)$ as the subscript on $q$ in which case one does not take the dual in $Z_0$. For a set of dual edges $A$, we will say $\Lambda\setminus A$ denotes deletion of the dual edges in~$A$ from $\Lambda^*$. Then under Dobrushin boundary conditions on $\Lambda_n$, the interface marginal is exactly proportional to the  random-line weight: for $\Gamma \subset \Lambda_n^*$,
\begin{align*}
    \bar\pi_{\Lambda}^{\rm Dob}(\Gamma) = \frac{q_\Lambda(\Gamma)}{\sum_{\Gamma':a_\Lambda^*\to b_\Lambda^*:\Gamma' \subset \Lambda}q_\Lambda(\Gamma')}\,.
\end{align*}

We next wish to take an infinite-volume limit of the weights~\eqref{eq:q-Lambda} when $\beta>\beta_c$. 
For any finite dual graph~$K$ containing vertices $u,v$ but not containing the edge~$e=uv$, the dual \emph{high-temperature expansion} of the Ising model gives 
\begin{equation}\label{eq:edge-addition}
 \frac{Z_0(K\cup e)}{Z_0(K)}=1+\tanh(\beta^*)\langle \sigma_u,\sigma_v\rangle_{K,\beta^*}\,,
\end{equation}
where $\langle \sigma_u,\sigma_v\rangle_{K,\beta^*}$ is the expectation of $\sigma_u\sigma_v$ under $\pi_{K,\beta^*}$. This duality means high-temperature decay of correlations can be used to study these ratios. In particular, this factor lies in~$[1,1+e^{-2\beta}]$ and converges under growing sequences of $K$ by the GKS inequality. Restoring any fixed finite set of edges one at a time therefore gives a positive limit for its partition-function ratio. In particular, for any fixed $\Gamma$, the limit of~\eqref{eq:q-Lambda} as $\Lambda\uparrow \mathbb Z^2$ exists and we denote it by
\begin{align}\label{eq:q-infty}
q_\infty(\Gamma)=\lim_{\Lambda\uparrow(\mathbb Z^2)^*}q_\Lambda(\Gamma)\,.\end{align}

\subsection{Geometry of the interface}\label{sec:interface-geometry}

Suppose that $\Gamma$ has source and sink points $a^*_{\Lambda},b^*_{\Lambda}$. The \emph{Ornstein--Zernike representation} of the interface is its decomposition into what are called \emph{cone-points} so that the increments between these cone-points behave like a weakly correlated directed 2-dimensional random walk. 

\begin{figure}
\centering
\includegraphics[width  = 0.9\textwidth]{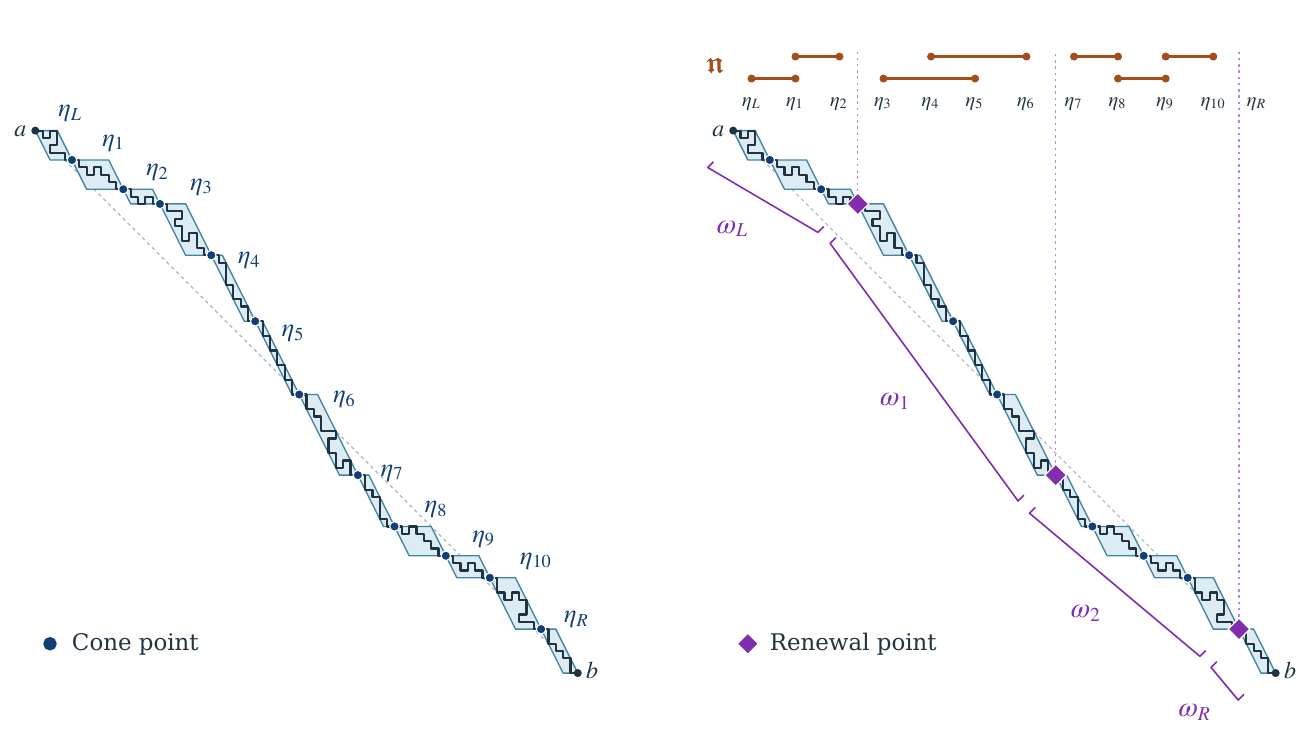}
    \caption{Left: The irreducible decomposition of a path $\Gamma$ from $a\to b$, with the cones shaded in blue, cone points marked by circles, and increments $\eta_L,\eta_1,...,\eta_{10},\eta_R$. Right: the interval percolation $\mathfrak{n}$ on top of the increments of the irreducible decomposition gives renewal points that are cone points not spanned by open intervals of $\mathfrak{n}$, separating the renewal components $\omega_i$, e.g., $\omega_1 = \eta_{3}\circ \cdots \circ \eta_6$.}\label{fig:OZ-decomp-decorations}
\end{figure}

 Since the diagonal Dobrushin boundary conditions have slope $-1$, throughout the paper we use the following \emph{cones}: 
\begin{equation}\label{eq:physical-cones}
 \mathcal C=\{(x,y):x\ge0,\ -2x\le y\le0\}\,.
\end{equation}
The diagonal ray of slope $-1$ from $(0,0)$ is evidently in the interior of $\mathcal C$. 
Following~\cite[Definition 5.1]{IOSV-prewetting}, a path from~$a$ to~$b$ is said to be forward-confined if it lies in~$a+\mathcal C$, backward-confined if it lies in~$b-\mathcal C$, and diamond-confined if it lies in~$(a+\mathcal C)\cap(b-\mathcal C)$. See  Figure~\ref{fig:diag-interface-OZ-decomp}, right, and Figure~\ref{fig:OZ-decomp-decorations}, left. 
\begin{definition}\label{def:cone-point}
    A cone-point~$v$ of~$\Gamma$ has its entire past (in the natural edge-by-edge parametrization of $\Gamma$ starting at $a$ and ending at $b$) in~$v-\mathcal C$ and its entire future in~$v+\mathcal C$.
\end{definition} 

Since~$\mathcal C$ does not contain a North, South, or West step, the path must enter a cone-point from the West and leave to the East. Cutting a path $\Gamma$ at each of its cone-points decomposes it into a sequence of \emph{increments} $(\eta_i)$ as 
\begin{align}\label{eq:Gamma-irreducible-rep}
 \Gamma=\eta_L\circ\eta_1\circ\cdots\circ\eta_m\circ\eta_R\,,\qquad m\ge0\,,
\end{align}
where $\circ$ denotes concatenation, 
 each $\eta_i$ is diamond-confined and delineated by a pair of cone-points, $\eta_L$ is backward-confined by the first cone-point, and $\eta_R$ is forward-confined by the last cone-point of $\Gamma$. The decomposition of Eq.~\eqref{eq:Gamma-irreducible-rep} is unique, and known as the \emph{irreducible representation} of $\Gamma$.

\paragraph{Surface tension.}\label{sec:cone-times}
When $\beta>\beta_c$, the dual inverse temperature~$\beta^* =\operatorname{arctanh}(e^{-2\beta})$ is subcritical. In particular there is a unique infinite-volume Gibbs measure $\pi_{\mathbb Z^2,\beta^*}$ at parameter $\beta^*$ and we can define the \emph{Green's function} $$G_\beta(x_1,x_2) = \langle \sigma_{(0,0)} \sigma_{(x_1,x_2)}\rangle_{\mathbb Z^2,\beta^*}$$ for the infinite-volume two-point function of this dual model. (Note that by translation invariance, $\langle \sigma_u \sigma_v \rangle_{\mathbb Z^2,\beta^*}$ only depends on the difference $v-u$.) 
The duality of the planar Ising model gives the following connection between interface weights at temperature $\beta$ and high-temperature two-point functions:  for $a,b\in (\mathbb Z^2)^*$ and $\beta>\beta_c$, 
\begin{align}\label{eq:random-line-rep-duality}
    \sum_{\Gamma:a\to b} q_\infty(\Gamma) = \langle \sigma_{a}\sigma_b\rangle_{(\mathbb Z^2)^*,\beta^*}= G_\beta(b-a)\,,
\end{align}
where the sum runs over all admissible paths $\Gamma\subset \mathbb Z^2$ from $a$ to $b$.
We define the inverse correlation length~$$\xi(x) = - \lim_{n\to\infty} \frac{1}{n} \log G_\beta(nx_1,nx_2)\,.$$ By duality~\eqref{eq:random-line-rep-duality}, $\xi$ is exactly the surface tension $\tau_\beta$ of the low-temperature model. In the notation of~\cite{CIV03}, the closed convex body associated with~$\xi$ is
\begin{equation*}
 \mathcal K_{\beta}=\{h\in\mathbb R^2:h\cdot x\le\xi(x)
                   \text{ for every }x\in\mathbb R^2\}\,.
\end{equation*}
By~\cite[Theorems 1.1 and B]{CIV03}, $\xi$ is a strictly convex norm whose unit sphere has positive curvature, and
\begin{equation}\label{eq:dual-convex-set}
 G_\beta(x)\le e^{-\xi(x)}\,,\qquad
 \operatorname{int}\mathcal K_{\beta}
 =\operatorname{int}\Big\{h:\sum_xe^{h\cdot x}G_\beta(x)<\infty\Big\}\,.
\end{equation}
For a slope~$s$, the unique (\cite[Section 3.4]{CIV03}) vector~$h(s)$ dual to~the associated vector $(1,s)$ is 
\begin{equation}\label{eq:dual-vector}
 h(s)\in\partial\mathcal K_{\beta},\qquad
 h(s)\cdot(1,s)=\xi(1,s)\,.
\end{equation}

In this paper, we fix the admissible slope range $\mathfrak{S}=[-7/4,-1/4]$. 
All slope-dependent Ornstein--Zernike estimates below will hold uniformly on any fixed compact subinterval of~$(-2,0)$, with constants depending on that subinterval and~$\beta$.

\subsection{Renewal representation of the interface}\label{sec:boundary-geometry}
The irreducible representation of~\eqref{eq:Gamma-irreducible-rep} has increments that are not independent of one another, but rather have exponentially decaying correlations. A robust tool that gives it a true renewal structure is the interval percolation  of~\cite[Theorem 5.9 and Appendix A.3]{IOSV-prewetting} which represents the correlations using an extra layer of independent randomness, which we call \emph{decorations}. 

\begin{definition}\label{def:decoration}
Given $\Gamma$, in its irreducible representation~\eqref{eq:Gamma-irreducible-rep}, call~$\eta_0=\eta_L$ and~$\eta_{m+1}=\eta_R$. The decoration process $\mathfrak{n}$ sets each interval~$I=[i,j]$, $0\le i<j\le m+1$, to \emph{open} independently with probability
\begin{equation}\label{eq:conditional-mark-law}
 p_I(\Gamma)=\frac{f_I}{1+f_I}\,,\qquad \text{where}\qquad 
 1+f_I=\frac{q_\infty(\eta_i\circ\cdots\circ\eta_j)\,q_\infty(\eta_{i+1}\circ\cdots\circ\eta_{j-1})}{q_\infty(\eta_{i+1}\circ\cdots\circ\eta_j)\,q_\infty(\eta_i\circ\cdots\circ\eta_{j-1})}\,.
\end{equation}
Let~$\mathcal P_\Gamma$ represent the conditional distribution of $\mathfrak{n}$. A cone-point separating $\eta_i$ and $\eta_{i+1}$ is a \emph{renewal point} of $(\Gamma,\mathfrak{n})$ if no open interval contains both indices. 
If the cone points separating $\eta_i,\eta_{i+1}$ and $\eta_{j},\eta_{j+1}$ are consecutive renewal points, then $\omega = \eta_{i+1} \circ \cdots \circ \eta_{j}$ is called a \emph{renewal component}. We will see that the renewal components are independent conditionally on the  renewal points. Each $\omega$ retains the information of its constituent increments and its open intervals in $\mathfrak{n}$. See Figure~\ref{fig:OZ-decomp-decorations}. 
\end{definition}

Just as in the horizontal case treated in~\cite{IOSV-prewetting}, there are constants $C,c$ such that~$0\le f_{[i,j]}\le Ce^{-c(j-i)}$. (Indeed, $\eta_i$ and~$\eta_j$ lie in opposite cones based at their internal endpoints, separated by at least~$j-i-1$ columns, so Lemma~\ref{lem:deletion} gives this bound.)

For a renewal component~$\omega=\eta_i\circ\cdots\circ\eta_j$, and $h\in \mathbb R^2$,  we define its $h$-tilted weight as
\begin{equation}\label{eq:decorated-component-weight}
 e^{h\cdot D(\omega)}\prod_{i'=i}^j q_\infty(\eta_{i'})
 \prod_{I\in\mathfrak n: I \subset [i,j]}f_I\,.
\end{equation}
where~$D(\omega)=(T(\omega),Y(\omega))$ is the \emph{displacement} between the left and right endpoints of $\omega$. For fixed $h$, we write~$\mathsf p_h$ for the probability distribution given by~\eqref{eq:decorated-component-weight} on possible renewal components. We write~$\rho^L_{\Lambda,h},\rho^R_{\Lambda,h}$ for the measures on renewal components that force them to stay in $\Lambda$ if started from $a$ or ending at $b$, respectively, so that these give the weights on the boundary renewal components. Every internal renewal component is diamond-confined, so its heights lie between its endpoint heights, and it needs no further restriction to $\Lambda$ since its diamond necessarily lies in $\Lambda$. 

Recall that~$\bar\pi_\Lambda$ denotes the interface marginal of~$\pi_\Lambda^{\mathrm{Dob}}$ and let~$\bar \pi^\infty = \bar\pi^\infty_\Lambda$ be the probability distribution proportional to~$q_\infty$ from~\eqref{eq:q-infty} on the same set of admissible paths in $\Lambda$ connecting $a_\Lambda^*\to b_\Lambda^*$. We use hats to denote their respective decorated laws, 
\begin{equation}\label{eq:marked-reference-definition}
 \widehat\pi_\Lambda(\Gamma,\mathfrak n)
 =\bar\pi_\Lambda(\Gamma)\mathcal P_\Gamma(\mathfrak n)\,,\qquad
 \widehat\pi^\infty(\Gamma,\mathfrak n)
 =\bar\pi^\infty(\Gamma)\mathcal P_\Gamma(\mathfrak n)\,.
\end{equation}

\begin{definition}\label{def:regular-interface}
A decorated interface is \emph{exceptional} if it only has $2$ or fewer renewal components. As long as it is not exceptional, it decomposes uniquely into renewal components as~$\omega_L\circ\omega_1\circ\cdots\circ\omega_k\circ\omega_R$, for some $k\ge1$, and we call $\omega_1,...,\omega_k$ its \emph{internal or bulk renewal components} while $\omega_L,\omega_R$ are its boundary renewal components. 
\end{definition}

Let~$R_{\Lambda,h}$ be the measure on exceptional interfaces with weights given by the product of the component weights in~\eqref{eq:decorated-component-weight}, and~$e^{h\cdot D(\Gamma)}q_\infty(\Gamma)$ for a path without any interior cone-point.

For~$s\in\mathfrak{S}$, put~$h=h(s)$ as in~\eqref{eq:dual-vector} and define the tilted partition function over paths in $\Lambda$ but using the homogeneous infinite-volume weights:  
\begin{equation*}
 \widehat Z_{\Lambda,h}(x)=e^{h\cdot x}\sum_{\Gamma:a_{\Lambda}^*\to b_{\Lambda}^*: \Gamma \subset \Lambda^*}q_\infty(\Gamma),
 \qquad x=b_\Lambda^*-a_\Lambda^*\,.
\end{equation*}
The next two propositions rewrite these weights on paths as products of weights over the renewal components, representing $\Gamma$ as a 2D directed random walk. 
The constants in what follows are  uniform in~$n$, the endpoint source and sink edges, and slopes in compact intervals of $(-2,0)$. These are essentially the same as~\cite[Appendix A]{IOSV-prewetting} except at slope $-1$ instead of slope $0$.  

\begin{proposition}\label{prop:domain-renewal}
Fix $s\in\mathfrak S$ and let $h=h(s)$ as in~\eqref{eq:dual-vector}. For every non-exceptional decorated interface
$(\Gamma,\mathfrak n)$ in $\Lambda$ with renewal decomposition
$\omega_L\circ\omega_1\circ\cdots\circ\omega_k\circ\omega_R$, we have
\begin{equation}\label{eq:decorated-interface-factorization}
 e^{h\cdot D(\Gamma)}q_\infty(\Gamma)\mathcal P_\Gamma(\mathfrak n)
 =
 \rho^L_{\Lambda,h}(\omega_L)
 \Bigl(\prod_{i=1}^k \mathsf p_h(\omega_i)\Bigr)
 \rho^R_{\Lambda,h}(\omega_R)\,.
\end{equation}
 In particular,  under~$\widehat\pi^\infty$, conditional on the ordered list of all renewal
points, the renewal components of the 
interface are independent. 
\end{proposition}

As a corollary, we get the following directed random-walk description of the renewal points of the decorated interface (with any tilt). 

\begin{corollary}\label{cor:OZ-random-walk}
Fix $s\in\mathfrak S$ and let $h=h(s)$ as in~\eqref{eq:dual-vector}. For a measure $\mu$ on renewal components, let $\mu^D$ denote its marginal on the component's displacement, and let $x=b_\Lambda^*-a_\Lambda^*$.
Summing~\eqref{eq:decorated-interface-factorization} over non-exceptional
decorated interfaces from $a_\Lambda^*$ to $b_\Lambda^*$ and adding
exceptional interfaces gives
\begin{equation}\label{eq:domain-renewal}
 \widehat Z_{\Lambda,h}(x)
 =
 \Bigl[
   (\rho^L_{\Lambda,h})^D
   *\Bigl(\sum_{k\ge1}(\mathsf p_h^D)^{*k}\Bigr)
   *(\rho^R_{\Lambda,h})^D
 \Bigr](x)
 +R_{\Lambda,h}^D(x)\,,
\end{equation}
where $*$ denotes convolution on $\mathbb Z^2$.
\end{corollary}

\begin{proposition}\label{prop:domain-component-bounds}
For some~$\sigma>0$, the renewal component distributions all satisfy
\begin{equation}\label{eq:domain-component-moments}
 \sum_{\omega} e^{\sigma|\omega|}\mu(\omega)=O(1)
 \qquad\text{where} \qquad \mu\in\{\mathsf p_h,\rho^L_{\Lambda,h},\rho^R_{\Lambda,h},R_{\Lambda,h}\}\,.
\end{equation}
The boundary measures, while not probability measures, each have total mass that is $\Omega(1)$. 
Finally, the part of~$\partial\mathcal K_{\beta}$ at slopes in $(-2,0)$ has an analytic parametrization~$h=(\lambda(\theta),\theta)$, whose mean slope~$\rho(\theta)=\E_{\mathsf p_h}[Y]/\E_{\mathsf p_h}[T]$ is strictly increasing and covers~$\mathfrak{S}$, with~$\rho(\theta)=s$ when~$h=h(s)$.
\end{proposition}

\begin{remark}\label{rem:near-squares-modification}
    Analogues of Propositions~\ref{prop:domain-renewal}--\ref{prop:domain-component-bounds} also hold for slight modifications of the domain $\Lambda$ (e.g., the near-squares of Section~\ref{sec:near-squares}) and under restrictions on their two boundary components with the natural modifications, so long as (a) the boundary distribution has an element of uniformly bounded length that gets uniformly positive tilted weight, and (b) once the two boundary components are fixed, every bridge of internal renewal components joining their renewal points stays in $\Lambda$.
\end{remark}

For completeness, we give proofs of Propositions~\ref{prop:domain-renewal}--\ref{prop:domain-component-bounds} in Appendix~\ref{proof:renewal-expansion}. We finally note that the exceptional paths where there are no internal renewal components are very unlikely: namely, because $R_{\Lambda,h}$ must span $n$ many columns, by Proposition~\ref{prop:domain-component-bounds} and Markov's inequality,
\begin{equation}\label{eq:exceptional-paths-unlikely}
 R_{\Lambda,h}(b_\Lambda^*-a_\Lambda^* \,,\,|\Gamma |\ge r)
 \le Ce^{-\sigma\max\{n,r\}}\,.
\end{equation}

\section{Mixing off of the interface}\label{sec:WPI-off-interface}

Our main result in this section is that conditioned on any interface $\Gamma$, the Ising Glauber dynamics satisfies a polynomial WPI on the complement of the spins fixed by the interface.

\begin{lemma}\label{lem:monochromatic-boundary-WPI}
Fix~$\beta>\beta_c$. For every~$C'>0$ there exists~$C$ such that, for every~$N$ and every finite domain~$U\subset\mathbb Z^2$ with~$|U|\le N$, Ising Glauber dynamics on~$U$ with all-plus boundary conditions satisfies an~$(N^C,N^{-C'})$-weak Poincar\'e inequality.
\end{lemma}

We prove a polynomial Poincar\'e inequality for a block dynamics that resamples blocks of side-length~$O(\log N)$ restricted to spin configurations in which every minus cluster has diameter at most~$\ell$, with~$\ell$ a sufficiently large multiple of~$\log N$. This restricted part of the state space has probability $1-n^{-100}$ under~$\pi_U^+$. 
To relate this block dynamics to the single-site Glauber dynamics, we  will use a crude mixing time bound on boxes at the logarithmic scale, coming from a general bound on mixing times of spin systems on graphs in terms of the cut-width of the underlying graph, that we reproduce here. Recall that the cut-width of a graph $W= (V(W),E(W))$ is defined by 
\begin{align}\label{eq:cutwidth-def}
 \operatorname{cw}(W)=\min_{(v_1,\ldots,v_{|V(W)|})}\max_{0\le i\le |V(W)|}
 \#\bigl\{\{u,v\}\in E(W):u\in\{v_1,\ldots,v_i\},\
 v\notin\{v_1,\ldots,v_i\}\bigr\}\,,
\end{align}
where the min is over labelings of the $|V(W)|$ vertices. 
The following is deduced from the canonical path method of~\cite{JS}, see e.g., Section~2.1 of~\cite{BKMP-trees}. 

\begin{lemma}\label{lem:cutwidth-mixing-bound}
    Fix any $\beta>0$, any finite graph $G= (V,E)$, and any boundary conditions $\tau \in \{\pm 1\}^{V\setminus W}$ for a subset $W \subset V$. Then Ising Glauber dynamics on $G^\tau$ has inverse spectral gap at most $(1+e^{8\beta}) |V| e^{ 4\beta {\rm cw}(G)}$. 
\end{lemma}

It will also be useful to have the following lemma which shows that an $\alpha$-Poincar\'e inequality for a chain restricted to a part of the state space of probability $1-\delta$ transfers to an $(\alpha,\delta)$-WPI to the unrestricted chain. 

\begin{lemma}\label{lem:conditioning-wpi}
Consider a state space $\Omega$ and a stationary distribution $\mu$. Let~$A\subset\Omega$ have probability~$p_A=\mu(A)>0$, suppose~$\mathcal E_\mu$ and~$\mathcal E_A$ are Dirichlet forms for Markov chains $P$ and $P_A$ reversible with respect to $\mu$ and $\mu(\cdot \mid A)$ respectively, with   
\begin{equation}\label{eq:conditioning-energy}
 \mathcal E_A(f|_A,f|_A)\le p_A^{-1} \mathcal E_\mu(f,f)
 \qquad\text{for all }f:\Omega\to\mathbb R\,.
\end{equation}
If $P_A$ satisfies an $(\alpha,\epsilon)$-WPI, then $P$ satisfies an $(\alpha,p_A \epsilon + (1-p_A))$-WPI. 
\end{lemma}
\begin{proof}
    Centering $f$ at~$\mu(f \mid A)$ bounds the contribution from~$A^c$ by~$\mu(A^c)\Osc(f)^2$, so  
    \begin{align*}
         \var_\mu(f)&\le p_A\var_{\mu_A}(f)+(1-p_A)\Osc(f)^2\,. 
    \end{align*}
Plugging this in to the $(\alpha,\eps)$-WPI for $P_A$ and using~\eqref{eq:conditioning-energy} gives the claim.   
\end{proof}

In particular, an~$\alpha$-Poincar\'e inequality for a chain $P_A$ that either rejects updates that leave~$A$ or does heat-bath resampling conditioned on $A$, together with~$\mu(A^c)\le\delta$ implies an~$(\alpha,\delta)$-weak Poincar\'e inequality on~$\Omega$. We are now ready to prove Lemma~\ref{lem:monochromatic-boundary-WPI}.

\begin{proof}[Proof of Lemma~\ref{lem:monochromatic-boundary-WPI}]
Fix $U\subset \mathbb Z^2$ and let~$\pi=\pi_U^+$. We say vertices $u$ and $v$ are $\star$-connected if $d_\infty(u,v)\le 1$. For an integer~$\ell\ge2$ and a finite domain~$A\subset\mathbb Z^2$, define
\begin{equation*}
 G_\ell(A)=\{\text{every minus $\star$-cluster in $A$ has $\ell^\infty$-diameter at most $\ell$}\},
 \qquad \mu_A=\pi_A^+(\,\cdot\mid G_\ell(A))\,.
\end{equation*}
By~\eqref{eq:loops-covariance-input} (specifically the exponential tails it implies on the dual FK-clusters which must surround any minus $\star$-cluster) and a union bound, we have for constants $c,C$ depending on $\beta$: 
\begin{align}
 \pi_A^+(G_\ell(A)^c)&\le C|A|e^{-c\ell}\,,\label{eq:off-interface-tail}\\
 0\le \pi_A^+(\sigma_x)-\pi_{A\cup\{y\}}^+(\sigma_x)&\le Ce^{-c\|x-y\|_\infty}\,,
 \qquad x\in A,\ y\notin A\,.\label{eq:off-interface-influence}
\end{align}

By monotonicity, namely adding one site at a time and telescoping~\eqref{eq:off-interface-influence}, for~$S\subset A\cap B$
\begin{equation*}
 \|\pi_A^+(\sigma_S\in \cdot)-\pi_B^+(\sigma_S \in \cdot)\|_{\tv}
 \le C\sum_{x\in S}\sum_{y\in A\oplus B}e^{-c\|x-y\|_\infty}\,.
\end{equation*}
where $A \oplus B$ denotes their symmetric difference. 
Conditioning a measure on an event changes it in total variation by the probability of the complement, so if~$|D_i|\le N$, $S\subset D_1\cap D_2$ and~$\dist_\infty(S,D_1\oplus D_2)\ge\ell$, then~\eqref{eq:off-interface-tail} and the triangle inequality give
\begin{equation}\label{eq:off-interface-local-TV}
 \|\mu_{D_1}(\sigma_S\in \cdot)-\mu_{D_2}(\sigma_S\in \cdot)\|_{\tv}\le\eta_\ell\,,
 \qquad \eta_\ell=C_1N^2e^{-c_1\ell}\,.
\end{equation}
Choose~$\ell= K\log N$ with~$K=K(\beta,C')$ large enough that~$\eta_\ell\le 1/10$ and~$CNe^{-c\ell}\le N^{-C'}$.

We will follow the block dynamics approach~\cite{Martinelli-notes}. Suppose~$A\subset R:=([u,v]\times[w,z])\cap\mathbb Z^2$, $|A|\le N$, and consider the distribution~$\mu_A= \pi_A^+(\cdot \mid G_\ell(A))$. Without loss of generality, its longer side is horizontal and has more than~$16\ell$ sites. Set~$m=(u+v)/2$ and define the overlapping halves
\begin{equation}\label{eq:off-interface-halves}
 B_L=A\cap\{(x_1,x_2): x_1\le m+2\ell\}\,,\qquad B_R=A\cap\{(x_1,x_2): x_1\ge m-2\ell\}\,.
\end{equation}
We are going to approximately tensorize variance across the blocks $B_L, B_R$, by showing that because of the plus boundary conditions and restriction to $G_{\ell}(A)$, the specific configuration on $B_L^c$ will not have much influence on $B_R^c$. 

Let $H_L$ be the set of all spins in $B_L^c$ together with all $\star$-connected minus components of $B_L^c$. Revealing $H_L$ also reveals its $+1$ $\star$-boundary, and we can let $D_L$ be all remaining sites it does not reveal information about. Define~$H_R,D_R$ symmetrically, as in Figure~\ref{fig:off-interface-halves}.

\begin{figure}[t]
\centering
\includegraphics[width  = 0.6 \textwidth]{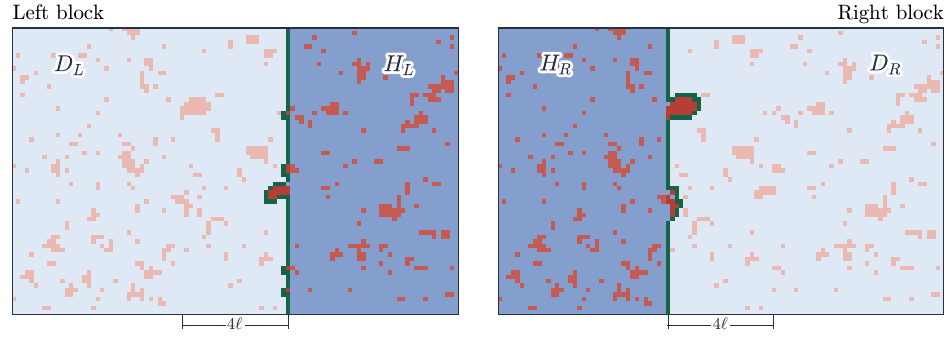}\caption{The minus $\star$-components of $B_L^c$ and $B_R^c$ are recorded by~$H_L,H_R$. The green boundaries are exposed to be plus spins, and~$D_L,D_R$ are the unrevealed sites.}
\label{fig:off-interface-halves}
\end{figure}

 Given a possible realization of~$H_L=h$, the domain Markov property gives plus boundary conditions on~$D_L(h)$, and the restriction of $G_{\ell}(A)$ becomes simply~$G_\ell(D_L(h))$. Hence
\begin{equation}\label{eq:off-interface-exact-law}
 \mu_A(\sigma_{D_i}\in\cdot\mid H_i=h)=\mu_{D_i(h)}\,,\qquad i\in\{L,R\}\,.
\end{equation}
On~$G_\ell(A)$ every explored minus cluster has span at most~$\ell$, and its plus neighbors add at most one layer. Therefore~$H_L,H_R$ are functions of the spins in the respective deterministic sets
\begin{equation*}
 S_L=A\cap\{(x_1,x_2): x_1\ge m+\ell-1\}\,,\qquad S_R=A\cap\{(x_1,x_2): x_1\le m-\ell+1\}\,.
\end{equation*}
Take two possible realizations of $H_L$,~$h,h'$, and put~$D=D_L(h)$, $D'=D_L(h')$. Their symmetric difference lies in~$S_L$, while~$S_R\subset D\cap D'$ and~$\dist_\infty(S_R,S_L)\ge2\ell-2\ge\ell$. By~\eqref{eq:off-interface-exact-law} and~\eqref{eq:off-interface-local-TV}, 
\begin{equation*}
 \sup_{h,h'}\|\mu_A(H_R\in\cdot\mid H_L=h)-\mu_A(H_R\in\cdot\mid H_L=h')\|_{\tv}
 \le\eta_\ell\,.
\end{equation*}
Averaging this bound over~$h'$ gives
$\|\mu_A[g\mid H_L]-\mu_A[g]\|_\infty
 \le 2\eta_\ell\|g\|_\infty$
for every bounded function~$g$ of~$H_R$. The same argument applies
with~$L$ and~$R$ interchanged. Since~$\eta_\ell\le 1/10$, applying
\cite[Lemma~3.1]{BertiniCancriniCesi02} to~$H_L,H_R$,
with~$p=\infty$ and~$\varepsilon=2\eta_\ell$, gives 
\begin{equation}\label{eq:off-interface-split}
 \var_{\mu_A}(f)\le 2\bigl\{\mu_A[\var(f\mid H_L)]+\mu_A[\var(f\mid H_R)]\bigr\}\,.
\end{equation}

Apply~\eqref{eq:off-interface-split} recursively, conditioning on each revealed history and averaging. In each step of this recursion, the longer side-length of our domain is cut by a factor of at least $3/4$, and therefore every two steps of the recursion reduce the maximal side-length by at least a factor of $3/4$. We then get 
\begin{align}\label{eq:off-interface-block-PI}
    \text{Var}_{\mu_U}(f) \le 2^{2 \log_{4/3} N} \sum_{R_i} \mu_U [\text{Var}(f\mid \sigma_{U\setminus R_i})]
\end{align}
where the sum runs over the possible rectangles we get from splitting the initial rectangle $B$ via $B_{L}$ or $B_R$, rotating, and iterating until the rectangle has side-length at most $16 \ell$. This is a Poincar\'e inequality for the block dynamics which does heat-bath updates on the minimal-size blocks $R_i\cap U$, with constant that is at most polynomial in $N$.

For a deterministic block~$Q = R_i \cap U$ and an outer configuration~$\tau \in \{-1,+1\}^{U\setminus Q}$, let~$p_\tau=\pi(G_\ell(U)\mid\sigma_{U\setminus Q}=\tau)$. By variance decomposition on each $\tau$ with~$p_\tau>0$, we have 
\begin{equation*}
 p_\tau\var_{\pi}(f\mid\tau,G_\ell(U))\le \var_\pi(f\mid\tau)\,.
\end{equation*}
Averaging over $\tau$ verifies~\eqref{eq:conditioning-energy} for these block updates. Lemma~\ref{lem:conditioning-wpi} and~\eqref{eq:off-interface-block-PI} now give the WPI,
\begin{equation*}
 \var_\pi(f)\le 2^{2\log_{4/3} N}\sum_{R_i}\pi[\var(f\mid\sigma_{U\setminus  R_i})]+\pi(G_\ell(U)^c)\Osc(f)^2\,.
\end{equation*}
Each~$U\cap R_i$ lies in a square of side~$16\ell$ and has cut-width at most~$16\ell$. Lemma~\ref{lem:cutwidth-mixing-bound} bounds its relaxation time, uniformly over $\tau$, by~$t_\ell=(1+e^{8\beta})(16\ell)^2e^{4\beta(16\ell+1)}$. Applying this bound on each block and summing the single-site Dirichlet forms, each site is counted only a bounded number of times. We conclude that
\begin{equation*}
 \var_\pi(f)\le C2^{2\log_{4/3} N} t_\ell\,\mathcal E_{\rm Gl}(f,f)+\pi(G_\ell(U)^c)\Osc(f)^2\,.
\end{equation*}
Since~$\ell=O(\log N)$ and~$\pi(G_\ell(U)^c)\le N^{-C'}$ by~\eqref{eq:off-interface-tail}, this proves the lemma.
\end{proof}

\begin{corollary}\label{cor:given-Gamma-WPI}
    Fix $\beta>\beta_c$. For any $C'$ there exists $C$ such that the following holds. For any interface $\Gamma$ on $\Lambda_n^{\rm dob}$, the Ising Glauber dynamics on $\Lambda_n \setminus \partial_v\Gamma$ with its induced boundary conditions (plus boundary on $\Lambda_n^\uparrow(\Gamma)$ and minus on $\Lambda_n^\downarrow(\Gamma)$) satisfies an $(n^{C},n^{-C'})$-weak Poincar\'e inequality. 
\end{corollary}
\begin{proof}
By Observation~\ref{obs:conditional-spin-law} the regions $\Lambda^\uparrow_n(\Gamma),\Lambda^\downarrow_n(\Gamma)$ carry independent Ising measures with monochromatic boundary conditions. Tensorize the variance across these two factors and apply Lemma~\ref{lem:monochromatic-boundary-WPI} to each with~$N=n^2$ (using spin flip symmetry for the minus region). The two single-site Dirichlet forms add, and choosing the error in each factor at most~$\tfrac12 n^{-C'}$ gives the claim. 
\end{proof}

\subsection{Pure phase relaxation from all plus}

Using the above weak Poincar\'e inequality together with the use of censoring to go from weak Poincar\'e to mixing, we already get polynomial mixing \emph{from all-plus initialization} on arbitrary domains in $\mathbb Z^2$ with all-plus boundary conditions. This result is alternatively easy to see by monotonicity and a cut-width bound on $\log n$-width domains (e.g., the $d=2$ results of~\cite{GhSi22}) but we give a version of the proof using the above WPI and censoring as another example of that combination. 

\begin{lemma}\label{lem:plus-init-mixing-with-plus-bdy}
    Fix $\beta>\beta_c$. Consider an arbitrary domain $\Lambda \subset \mathbb Z^2$ with all-plus boundary conditions. The Ising Glauber dynamics on $\Lambda^+$, initialized from all-plus, i.e., $X_{+,0}(v) =+1$ for all $v\in \Lambda$, has 
    \begin{align*}
        \|\mathbb P(X_{+,t}\in \cdot ) - \pi^+_{\Lambda}\|_{\tv} \le |\Lambda|^{-5} \qquad \text{for $t\ge |\Lambda|^C$}\,.
    \end{align*}
\end{lemma}
\begin{proof}
    Fix an arbitrary ordering of the vertices of $\Lambda$ as $v_1,...,v_{|\Lambda|}$, and let $\Lambda^{(i)}$ be the graph induced on $v_1,...,v_{i}$. Let $\overline{X}_{+,t}$ be the censored version of $X_{+,t}$ with censoring scheme (recall Definition~\ref{def:censoring}) given by the growing domains $(\Lambda^{(i)})_{i\le |\Lambda|}$ and times $t_0 =0$ and $t_i - t_{i-1} = C_0 |\Lambda|^{C_1}$ where $C_1$ is bigger than the constant $C$ of Lemma~\ref{lem:monochromatic-boundary-WPI} for $C' =100$, and $C_0$ is a sufficiently large constant. Under the monotone coupling, 
    \begin{align*}
        \|\mathbb P(X_{+,t}\in \cdot ) - \pi^+_{\Lambda}\|_{\tv}\le \sum_{v\in \Lambda} \Big(\mathbb P(X_{+,t}(v)=1) - \pi_{\Lambda}^+(\sigma_v =1)\Big)\,.
    \end{align*}
    By the censoring theorem of~\cite{PWcensoring}, each summand on the right-hand side is at most $\mathbb P(\overline{X}_{+,t}(v) =1) - \pi_{\Lambda}^+(\sigma_v=1)$. We claim inductively that, 
    \begin{align}\label{eq:mixing-from-all-plus-inductive}
        \|\mathbb P(\overline X_{+,t_{i}} (\Lambda^{(i)})\in \cdot ) - \pi_{\Lambda^{(i)}}^+\|_{\tv } \le i |\Lambda|^{-10}\,.
    \end{align}
    The base case is trivially true since it is a single vertex with all-plus boundary conditions, and the censoring effectively induces the same all-plus boundary conditions, so the system mixes as soon as a clock rings. Suppose~\eqref{eq:mixing-from-all-plus-inductive} holds for some $i$. Then except with probability $i |\Lambda|^{-10}$, we can couple  $\overline{X}_{+,t_i}(\Lambda^{(i+1)})$ to the initial state of a chain $Y_{0}^{(i+1)}\sim \pi_{\Lambda^{(i)}}^+ \sqcup +$ where $\pi_{\Lambda^{(i)}}^+\sqcup +$ is the concatenation of a sample from $\pi_{\Lambda^{(i)}}^+$ with a plus site at $\Lambda^{(i+1)}\setminus \Lambda^{(i)}$. Using the identity coupling of the dynamics, this gives a coupling of $\overline{X}_{+,t_i+s}(\Lambda^{(i+1)})$ to $Y_s^{(i+1)}$ for all $0 \le s \le C_0 |\Lambda|^{C_1}$. The initialization $Y_{0}^{(i+1)}$ has Radon--Nikodym derivative to $\pi_{\Lambda^{(i+1)}}^+$ bounded by some constant $C_\beta$. Therefore by~\eqref{eq:WPI-consequence} and Lemma~\ref{lem:monochromatic-boundary-WPI}, $Y_{C_0 |\Lambda|^{C_1}}^{(i+1)}$ has total-variation distance at most $|\Lambda|^{-10}$ to $\pi_{\Lambda^{(i+1)}}^+$. Combining these two errors proves the inductive claim~\eqref{eq:mixing-from-all-plus-inductive} at $i+1$. Taking the final $i=|\Lambda|$ in~\eqref{eq:mixing-from-all-plus-inductive}, and multiplying by $|\Lambda|$ for the sum over vertices in $\Lambda$, implies the desired. 
\end{proof}

\section{Mixing of Ornstein--Zernike Interfaces}\label{sec:interface-wpi}

In this section, we show fast mixing for the OZ interfaces of the Ising model under $O(\log n)$-width block resampling moves. While we eventually prove what amounts to a WPI for functions that only depend on the interface $\Gamma$ with respect to the Glauber Dirichlet form, most of the work of the section focuses on analyzing the following Markov chain on interfaces $\Gamma$ of $\Lambda_n$ with diagonally Dobrushin boundary conditions.

\begin{definition}\label{def:interface-dynamics-finite-volume}
    Fix a block-width $w\ge 1$. For each $j=0,...,n-w$ define the block $J_j = \{j,...,j+w\}\times \{1,...,n\}$. For a block $J$, let $\mathcal F_{J}$ be the $\sigma$-algebra generated by all increments of the OZ irreducible decomposition~\eqref{eq:Gamma-irreducible-rep} of $\Gamma$ that intersect $\Lambda_n \setminus J$. The \emph{middle portion} $\eta_J$ of $\Gamma$ is the set of increments connecting the two outermost cone-points in $J$. Recalling that $\bar\pi$ is the marginal of $\pi$ on $\Gamma$, the \emph{finite-volume interface dynamics} is the Markov chain that
    \begin{enumerate}
    \item Assigns each $j \in \{0,...,n-w\}$ a rate-one Poisson clock
        \item When the clock at $j$ rings, it resamples $\eta_{J_j}$ from $\bar \pi(\cdot \mid \mathcal F_{J_j})$. 
    \end{enumerate}
    (Note that the number of increments in the middle portion $\eta_{J_j}$ can change under resampling.) We use $P_{\rm Int}$ for the semigroup of this chain, and note it has the following Dirichlet form: 
\begin{equation}\label{eq:fine-form}
 \mathcal E_{\mathrm{Int}}(g,g)=\sum_{0\le j \le n-w}\bar\pi\bigl[\var(g\mid\mathcal F_{J_j})\bigr]\,.
\end{equation}
\end{definition}
One may be worried that this chain is not ergodic on the set of all interfaces. Because we only need to prove a WPI, we may take $w = O(\log n)$ and restrict attention to the set of interfaces having no increment of size larger than $w/3$, say, in which case it is easily seen to be ergodic on that part of its state space. 
The main result of this section is the following WPI for the interface dynamics. 

\begin{lemma}\label{lem:interface-WPI}
For every~$\beta>\beta_c$ and $C'$, there are $C, C''$ such that the finite-volume interface dynamics on $\Lambda_n$ with block-width $w = C'' \log n$ satisfies an $(n^{C},n^{-C'})$-WPI. That is, for every function $g$ of the interface $\Gamma$, 
\begin{align*}
    \text{Var}_{\bar \pi}(g) \le n^C \mathcal E_{\rm Int}(g,g) + n^{ - C'} \Osc (g)^2\,.
\end{align*}
\end{lemma}

Lemma~\ref{lem:interface-WPI} will follow the idealized strategy described in Section~\ref{sec:proof-sketch-interface-mixing} (although technically through a chain that updates the renewal representation instead of the irreducible representation). In order to apply this in the variance decomposition of~\eqref{eq:intro-total-variance} to get the weak Poincar\'e for the full Glauber dynamics, we then compare the Dirichlet form $\mathcal E_{\rm Int}(g,g)$ to $\mathcal E_{\rm Gl}(g,g)$ at a further comparison polynomial cost by Lemma~\ref{lem:cutwidth-mixing-bound} because the blocks $J$ have $O(\log n)$ width.

\subsection{Ornstein--Zernike inputs}\label{sec:renewal-construction}
Before getting to the dynamical analysis, we record the equilibrium inputs on the OZ interface distribution. These follow from a mixture of lemmas in~\cite{VelenikThesis,PfisterVelenik,CIV03,IOSV-prewetting,Ioffe14}. In Appendix~\ref{app:OZ-inputs}, we establish the lemmas in the specific formulations we present here.

The main ingredients will be exponential tails for the lengths of the renewal components, a comparison of the finite-volume and infinite-volume weights on interfaces, Gaussian tail bounds on height fluctuations away from the chord connecting the boundary points, and anti-concentration for the height fluctuations in the form of local limit theorems for the renewal process. 

\subsubsection{Exponential tails and domain modifications}

The first lemmas we recall are standard tail bounds on the (un-normalized) random-line weights $q_\Lambda(\Gamma)$ of~\eqref{eq:q-Lambda}. The first is an exponential tail bound on the weight of random-line paths of length at least $m$, and follows directly from summing~\cite[Lemma~4.2.6]{VelenikThesis}.

\begin{lemma}\label{lem:raw-contour-tail}
Fix $\beta>\beta_c$. There exist $c_\beta, C_\beta$ such that for any~$m\ge1$,  any $a,b$, and any $\Lambda$, 
\begin{equation}\label{eq:raw-path-length-tail}
 \sum_{\Gamma:a\to b,\ |\Gamma|\ge m}q_\Lambda(\Gamma)
 \le C_\beta e^{-c_\beta m}\,,
\end{equation}
where the sum is over admissible paths with prescribed endpoints at~$a,b$. 
\end{lemma}

The next lemma controls the total weight of paths connecting two vertices by their dual two-point function. It is a direct consequence of the duality relation~\eqref{eq:random-line-rep-duality} and the GKS inequality. 

\begin{lemma}\label{lem:raw-source-sum}
For distinct endpoints~$a,b$, for every $\Lambda$, we have 
\begin{equation}\label{eq:raw-source-sum}
 \sum_{\Gamma:a\to b}q_\Lambda(\Gamma)\le G_\beta(b-a) \le e^{ - \xi (b-a)}\,.
\end{equation}
\end{lemma}

Finally, we recall the following random-line BK inequality~\cite[Lemma 6.5]{PfisterVelenik}. 
\begin{lemma}\label{lem:physical-through-point}
For every vertex~$z\ne a,b$, and every $\Lambda$, we have 
\begin{equation}\label{eq:through-point}
 \sum_{\Gamma:a\to b,\ z\in\Gamma}q_\Lambda(\Gamma)
 \le  G_\beta(z-a)G_\beta(b-z)
 \le e^{-\xi(z-a)-\xi(b-z)}\,.
\end{equation}
\end{lemma}

All three bounds above extend to $\Lambda = \mathbb Z^2$ by taking a limit $\Lambda \uparrow \mathbb Z^2$.

 We also need the following domain decoupling of~\cite[Lemma~6.8]{PfisterVelenik} and~\cite[Lemma~3.1]{CIV03}, even when arbitrary edges have been deleted from the domain. Recall $Q_\eta$ from the paragraphs following Definition~\ref{def:interface}  ($\eta$ together with edges it reserves due to the splitting rule).

\begin{lemma}\label{lem:domain-weight-comparison}
Fix~$\beta>\beta_c$. There are~$c_\beta,C_\beta>0$ such that the following holds. Consider finite~$\Lambda_1,\Lambda_2 \subset (\mathbb Z^2)^*$  and let $B=E(\Lambda_1)\oplus E(\Lambda_2)$ be their symmetric difference. Suppose that a path~$\eta$ is admissible in both graphs.  Then
\begin{equation}\label{eq:domain-weight-comparison}
 \left|\log\frac{q_{\Lambda_1}(\eta)}{q_{\Lambda_2}(\eta)}\right|
 \le C_\beta \sum_{e\in Q_\eta, f\in B} e^{ - c_\beta d(e,f)}\,.
\end{equation}
If~$\Lambda_1\subset\Lambda_2$, the logarithm in~\eqref{eq:domain-weight-comparison} is nonnegative. These bounds apply both to entire interfaces and to (concatenations of) increments, and by a limiting procedure also to $q_\infty$. 
\end{lemma}

\begin{lemma}\label{lem:domain-law-comparison}
With the notation of Lemma~\ref{lem:domain-weight-comparison}, and arbitrary reweighting factors $a(\eta)$, let~$\nu_i$ be the probability distributions with weights~$a(\eta)q_{\Lambda_i}(\eta)$ on a family of~$\eta \in \mathcal A$. If~$\sum_{e\in Q_\eta,f\in B} e^{- c_\beta d(e,f)}\le M$ for all~$\eta \in \mathcal A$, then we have the Radon--Nikodym bound, 
\begin{equation}\label{eq:domain-radon-nikodym}
 e^{-2C_\beta M}
 \le\frac{d\nu_1}{d\nu_2}
 \le e^{2C_\beta M}\,.
\end{equation}
In particular, changing at most~$k$ edges of the background graph costs at most~$e^{C_\beta k}$ in the weights on $\eta\in \mathcal A$. The flexibility of having the additional weights $a$ is to allow for common $h$-tilts $e^{ h\cdot D(\eta)}$.
\end{lemma}

If every path in~$\mathcal A$ has length at most~$R$ and its reserved edges are at distance at least~$r$ from~$B$, the exponents in~\eqref{eq:domain-radon-nikodym} are at most~$C_\beta(R+1)e^{-c_\beta r}$. The proof of these lemmas is in Appendix~\ref{proof:domain-weight-comparison}.

\subsubsection{Diffusivity and local limit theorems in the renewal representation}

Recall the representation~\eqref{eq:domain-renewal} of the law over decorated interfaces from Section~\ref{sec:boundary-geometry}, and recall that the internal renewal comonents are independent before conditioning on their displacement. 
Write $\mathsf p_s=\mathsf p_{h(s)}$ for the law~\eqref{eq:decorated-component-weight} of an internal renewal component~$\omega$, with displacement $D(\omega)=(T(\omega),Y(\omega))$ and edge length $|\omega|$ with tilt $h(s)$ from~\eqref{eq:dual-vector}. Recall $\widehat \pi^\infty$ is the law on decorated interfaces in $\Gamma$ from~\eqref{eq:marked-reference-definition}. We define a cutoff threshold $Q$ and restrict attention to the decorated interfaces that are nice relative to that threshold. 

\begin{definition}\label{def:B_Q-set}
    For a cutoff $Q$, let $\mathcal B_Q$ be the event that the decorated interface $(\Gamma,\mathfrak{n})$ on $\Lambda_n$ is not exceptional (it has at least three renewal components) and every renewal component, including the two boundary components, has total number of edges $|\omega_i|\le Q$. 
\end{definition}

As with Propositions~\ref{prop:domain-renewal}--\ref{prop:domain-component-bounds}, the following extend to near-squares and other geometries as long as (a)--(b) in Remark~\ref{rem:near-squares-modification} hold. The following shows that if $Q \ge C \log n$ for a large constant $C$, then $\mathcal B_Q^c$ is polynomially rare. 

\begin{lemma}\label{lem:OZ-exponential-tails}
There are constants $c,C$ such that on a square of width $n$,
\begin{equation}\label{eq:OZ-cutoff-error}
    \widehat\pi^\infty(\mathcal B_Q^c)
       \le Cn^{3/2}e^{-cQ},\qquad \text{for all}\qquad  1\le Q\le n\,.
\end{equation}
\end{lemma}
Lemma~\ref{lem:OZ-exponential-tails} is essentially a consequence of Proposition~\ref{prop:domain-component-bounds}, and is proved in Appendix~\ref{proof:renewal-estimates}.

We will consider the distribution over decorated interfaces conditioned on the event $\mathcal B_Q$. For~$s\in\mathfrak S$ and~$Q$ sufficiently large, there is a vector~$h_Q(s)$ such that
\begin{equation}\label{eq:tilted-p-distribution}
 \mathsf p_{s,Q}(\omega)
 =\mathsf p_s(\omega)\mathbf1_{\{|\omega|\le Q\}}
    e^{(h_Q(s)-h(s))\cdot D(\omega)}
\end{equation}
is a probability distribution and has mean slope~$\E [Y(\omega)]/\E [T(\omega)]=s$. (We also apply the same $Q$-cutoff and tilt to the boundary measures $\rho^L, \rho^R$ and normalize those.) On conditioning the concatenation of the renewal components to join points~$a\to b$, the exponents of $(h_Q(s)- h(s)) \cdot D(\omega)$ add up to a constant depending only on $b-a$ since that is the total displacement. Thus the distribution on decorated interfaces from $a\to b$ does not depend on the choice of tilt. We write~$\widehat\pi^\infty_{a,b,Q}$ for this distribution on decorated interfaces from $a\to b$ conditioned on $\mathcal B_Q$.

The probability that a directed 2D random walk started at $(0,0)$ with displacement distribution given by~$\mathsf p_{s,Q}^D$ has a renewal point at~$(x,y)$ is exactly
\begin{equation}\label{eq:OZ-renewal-function}
 U_{s,Q}(x,y)=\sum_{j\ge0}(\mathsf p_{s,Q}^D)^{*j}(x,y)\,.
\end{equation}
(where we are decomposing by the number of renewal components, and $*j$ denotes the $j$-fold convolution of the marginals on their displacements). 
We use the following local limit estimate of~\cite[Section~3.1, Proposition~3.2 and (89)--(96)]{Ioffe14}.

\begin{lemma}\label{lem:OZ-local-limit}\label{prop:renewal-estimates}
For every fixed~$A$, there exist constants $c_A, C$ such that uniformly over~$s\in\mathfrak S$, $Q\ge Q_0$ and~$x\ge L_AQ$ for sufficiently large~$Q_0,L_A$,
\begin{equation}\label{eq:OZ-local-limit}
 \sup_y U_{s,Q}(x,y)\le Cx^{-1/2},\qquad
 U_{s,Q}(x,y)\ge c_Ax^{-1/2}\quad\text{if }|y-sx|\le A\sqrt x\,.
\end{equation}
The bounds also hold after convolution with either or both normalized boundary measures $\rho^L_{s,Q}, \rho^R_{s,Q}$. Without the cutoff at $Q$, the $c/\sqrt{x}$ lower bound holds for sufficiently large $x$. 
\end{lemma}

If we force the two boundary renewal components at $\Omega(1)$ cost, then conditionally apply~\eqref{eq:OZ-local-limit} we also conclude the following useful lower bound on the total weight of $a\to b$ paths~\eqref{eq:domain-renewal}: 
\begin{equation}\label{eq:physical-reference-normalizer}
 \widehat Z_{\Lambda,h(s)}(b-a)\ge c n^{-1/2}\,.
\end{equation}
The same lower bound holds for the tilted finite-volume sum $\sum_{\Gamma \subset\Lambda: \Gamma: a\to b} e^{h(s) \cdot (b-a)} q_{\Lambda}(\Gamma)$ because~$q_\Lambda(\Gamma)\ge q_\infty(\Gamma)$ for every $\Gamma\subset \Lambda$.

The following complements the anticoncentration of Lemma~\ref{lem:OZ-local-limit} with a Gaussian concentration bound for the renewal points. Such a Gaussian tail bound on the maximum deviation of the interface from the line connecting its endpoints was shown in~\cite[Theorem 5.3 and Eq. (5.3)]{LMST} when the slope is $0$ and without any truncation or cutoff by $Q$. 
For completeness we give a proof with the truncation and at slope $s\in \mathfrak{S}$; since the renewal components of the OZ bridge are exchangeable and have good moment bounds, this will follow from standard martingale arguments.   

\begin{lemma}\label{lem:OZ-bridge-concentration}\label{prop:bridge-fluctuations}
For $a = (a_1,a_2)$ and $b= (b_1,b_2)$, let~$\ell=b_1-a_1\ge LQ$, $Q\ge Q_0$, and~$s=(b_2-a_2)/\ell$. Write~$\Delta_{a,b}$ for the largest vertical distance of a renewal point from the chord joining~$a$ to~$b$. For~$s\in\mathfrak S$ and~$u\ge0$,
\begin{equation}\label{eq:OZ-bridge-concentration}
 \widehat\pi^\infty_{a,b,Q}(\Delta_{a,b}\ge u)
 \le C\exp(-c\min(u^2/\ell,u))\,.
\end{equation}
In particular, $\E[\Delta_{a,b}]\le C\sqrt\ell$ for a uniform constant $C$. Furthermore, there is a different constant $C$ such that for every $Q \ge Q_0$ and every slope $s$ in $[-2,0]$, one has  
\begin{equation}\label{eq:OZ-bridge-rough-moment}
 \E_{\widehat\pi^\infty_{a,b,Q}}[\Delta_{a,b}^2]\le CQ\ell\,.
\end{equation}
\end{lemma}

These anticoncentration and concentration estimates are proved in Appendix~\ref{proof:renewal-estimates}. Using the latter, we can show that the diagonal interface in $\Lambda_n$ stays away from the boundary, except near the endpoints, and therefore the law is not sensitive to small modifications of the shape of the domain.

\begin{lemma}\label{lem:OZ-finite-volume-comparison}
On the admissible paths in~$\Lambda$, set $\psi_\Lambda(\Gamma)=q_\Lambda(\Gamma)/q_\infty(\Gamma)$. Then $\psi_\Lambda\ge1$ and
\begin{equation}\label{eq:OZ-finite-volume-density}
    \frac{d\bar\pi}{d\bar\pi^\infty}
       =\frac{\psi_\Lambda}{\bar\pi^\infty[\psi_\Lambda]}\,.
\end{equation}
There are constants $b,c,C>0$ such that, on every square of width~$n$ and for every $R\ge1$,
\begin{equation}\label{eq:OZ-likelihood-tail}
    \max\{\bar\pi(\psi_\Lambda>e^{bR}), \bar\pi^\infty(\psi_\Lambda>e^{bR})\}
       \le Cn^{5/2}e^{-cR}\,.
\end{equation}
The same statements hold under the decorated laws, since both use the same conditional law $\mathcal P_\Gamma$ for the decoration. The bounds also hold mutatis mutandis on the near-squares of  Definition~\ref{def:near-squares}. 
\end{lemma}

The proof of Lemma~\ref{lem:OZ-finite-volume-comparison} is in Appendix~\ref{proof:boundary-influence}.

\subsection{The infinite-volume interface dynamics} 

The following ``infinite-volume chain" is more homogeneous than Definition~\ref{def:interface-dynamics-finite-volume} as it uses $\mathbb Z^2$-weights on the increments rather than $\Lambda_n$-weights. Note that interfaces drawn from the distributions $ \widehat \pi^\infty$ were still confined to $\Lambda_n$; the difference was they used the translation invariant weights $q_\infty$ from~\eqref{eq:q-infty} on the internal increments. This chain will also break up the correlations between the increments of $\Gamma$ by allowing the Markov chain to also resample the decorations $\mathfrak{n}$ which were introduced to factorize the correlations. Whereas the finite-volume chain resampled between outermost cone-points of blocks $J$ from Definition~\ref{def:interface-dynamics-finite-volume}, the infinite-volume chain will resample between outermost renewal points of $J$. 

\begin{definition}\label{def:infinite-volume-interface-dynamics}
    The \emph{infinite-volume interface  dynamics} is the following Markov chain on decorated interfaces. It follows the update rule of Definition~\ref{def:interface-dynamics-finite-volume} except when the clock at block $J$ rings, it fixes the information $\mathcal H_J$ which is all renewal components intersecting $\Lambda_n \setminus J$ (all the increments constituent to these components, and their decorations). So long as $\mathfrak{n}$ has at least three distinct components, it resamples the decorated interface $(\Gamma,\mathfrak{n})$ from the conditional distribution $\widehat\pi^\infty(\cdot\mid\mathcal H_J)$. 

    We use $\widehat P_{\rm Int}$ to denote the semigroup of this chain, and note that it has Dirichlet form 
    \begin{align*}
        \widehat{\mathcal E}_{\rm Int}(g,g) =    \sum_{j=0}^{n-w}\widehat\pi^\infty[\var(g\mid\mathcal H_{J_j})]\,.
    \end{align*}
\end{definition}

\begin{lemma}\label{lem:infinite-volume-WPI}
    Fix $\beta>\beta_c$. 
    There exist constants $p,c,L,Q_0>0$ such that, for all $Q_0\le Q \le n/L$, and block-width $w\ge LQ$, for any function $g$ of decorated interfaces $(\Gamma,\mathfrak{n})$, 
\begin{equation}\label{eq:reference-recursion}
    \var_{\widehat\pi^\infty}(g)
       \le n^p\widehat{\mathcal E}_{\rm Int}(g,g)
            +n^pe^{-cQ}\Osc(g)^2\,.
\end{equation}
In particular, for every $C'$, there are constants $C, C''$ such that if  $w= C'' \log n $, then the infinite-volume interface dynamics satisfies an $(n^C,n^{-C'})$-WPI. 
\end{lemma}

We prove Lemma~\ref{lem:infinite-volume-WPI} using a block dynamics~\cite{Martinelli-notes} style reduction in scale. 
Fix any two endpoints $a,b$  at slope in $\mathfrak{S}$ and recall the law $\widehat\pi^\infty_{a,b,Q}$ using the cutoff threshold $Q$ of Definition~\ref{def:B_Q-set}.
By translation, take $a=(0,0)$ and $b=(x,\rho x)$, where
$\rho\in\mathfrak S=[-7/4,-1/4]$. Define
\begin{align}\label{eq:West-East-blocks}
 B_{\rm West}=[0,2x/3]\times\mathbb R\,,
 \qquad B_{\rm East}=[x/3,x]\times\mathbb R\,.
\end{align}
For either block $B \in \{B_{\rm West},B_{\rm East}\}$, let $\mathcal H_B$ be the sigma-algebra generated by the renewal components
that are not wholly contained in $B$ (both the constituent increments and
decorations).

\begin{definition}\label{def:block-interface-dynamics}
In the \emph{two-block interface dynamics}, each block $B_{\rm West}, B_{\rm East}$ has an independent
rate-one clock. When the clock of $B$ rings, resample the renewal
components wholly contained in $B$ (increments and their decorations),
from $\widehat\pi^\infty_{a,b,Q}(\cdot\mid\mathcal H_B)$.
\end{definition}

The following shows that the two-block interface dynamics has constant spectral gap. When recursively applied $\log n$ many times starting from width $x=n$ and slope $\rho =-1$, it reduces to the infinite-volume interface dynamics at a polynomial cost in the spectral gap. 
\begin{lemma}\label{lem:two-block-interface-constant-gap}
    Fix $\beta>\beta_c$ and let $\rho$ belong to the set of admissible slopes $\mathfrak{S} = [-7/4,-1/4]$. There are constants $Q_0,L,C_0$ such that, whenever $Q\ge Q_0$ and $x\ge LQ$, the two-block dynamics satisfies
    \begin{align}\label{eq:reference-two-window}
        \text{Var}_{\widehat \pi^\infty_{a,b,Q}}(g) \le C_0 \widehat \pi_{{a},b,Q}^\infty\big[\var(g\mid\mathcal H_{B_{\mathrm{West}}})
                         +\var(g\mid\mathcal H_{B_{\mathrm{East}}})\big]
    \end{align}
\end{lemma}
\begin{proof}
 We will use $z$ to denote a state of the decorated interface from $(0,0)\to (x,\rho x)$. For any renewal point $(w,\phi_z(w))$ in $z$, we think  of $\phi_z$ as the function which takes the $x$-coordinate of a renewal point of $z$ and returns its \emph{height}. We also want to specify the (random) renewal locations nearest the boundaries of the
blocks $B\in \{B_{\rm West}, B_{\rm East}\}$, and use those to define heights at $x$-coordinates $x/3$ and $2x/3$. Let
\begin{equation}\label{eq:interface-height-V}
\begin{aligned}
 &{\mathsf{u}}(z)=\min\{u\ge x/3:u\text{ is a renewal of }z\}\,, \\ 
 &H_{\sf u}(z)=\phi_z({\mathsf{u}(z)})-\rho {\mathsf{u}(z)},\qquad V_{\sf{u}}(z)=1+\frac{|H_{\sf u}(z)|}{\sqrt x}\,.
\end{aligned}
\end{equation}
and similarly define ${\mathsf{v}(z)} = \max\{u\le 2x/3: u \text{ is a renewal of }z\}$ and $H_{\sf v}(z), V_{\sf{v}}(z)$. Here, $\mathsf{u}$ is giving the left-most renewal in $B_{\rm East}$, $H_\mathsf{u}$ is its relative height to the line $(0,0)\to (x,\rho x)$, and $V_{\mathsf{u}}$ measures that relative height on the diffusive scale.
The $Q$-cutoff implies $x/3 \le {\sf u} \le x/3 +Q$. 

The below  first shows that $V_{\sf u}(z)$ contracts to the $O(1)$ scale at constant rate, then couples two interfaces $z,z'$ which have $V_{\sf u}(z) + V_{\sf u}(z') = O(1)$ in a further $O(1)$ block updates. 
We will first prove the bound for a chain that performs sweeps sequentially updating $B_{\rm East}$ then $B_{\rm West}$, then at the end transfer that spectral gap bound to the continuous-time two-block chain. 

\begin{figure}[t]
\centering
\includegraphics[width  = 0.95 \textwidth]{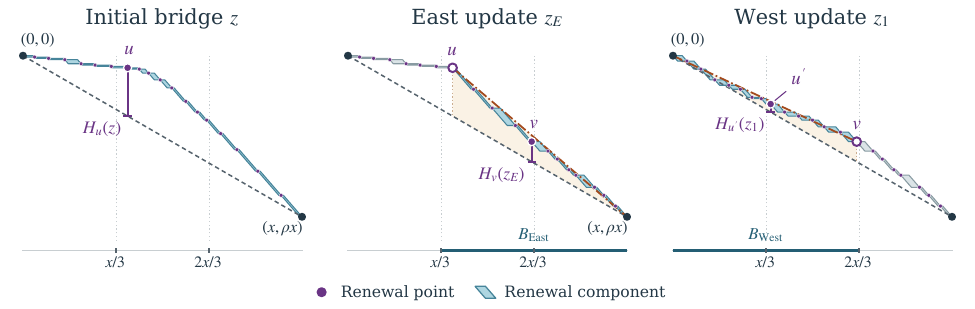}\caption{One sweep of the two-block dynamics, starting from configuration (decorated interface) $z$, after an update on $B_{\rm East}$, one gets $z_E$, and after a subsequent update on $B_{\rm West}$, one gets $z_1$. The update on $B_{\rm East}$ resamples all renewal components to the right of renewal point $u$ (the left-most one in $B_{\rm East}$ in $z$). By the Gaussian concentration, after one sweep, if the height in $z$ at the renewal point near $x/3$ was super-diffusive, in expectation, the height in $z_1$ at the renewal point near $x/3$ will be a factor smaller than in $z$.}
\label{fig:oz-sweep}
\end{figure}

\smallskip
\noindent \emph{Step 1: Contraction of the relative height.}  
Starting from a state~$z$, consider an East update immediately followed by a West update. Let~$z_E$ be the state after the East update and~$z_1$ the state after the subsequent West update. Let~$u=\mathsf u(z)$, $v=\mathsf v(z_E)$, and~$u'=\mathsf u(z_1)$. Then we can define deviations $Z_E, Z_W$ from the expected chords such that 
\begin{align}
 H_{\mathsf v}(z_E)
 &=\frac{x-v}{x-u}H_{\mathsf u}(z)+Z_E,
 &\frac{x-v}{x-u}
 &\le\frac12+\frac {3Q}x,
 \label{eq:actual-renewal-drift}\\
 H_{\mathsf u}(z_1)
 &=\frac{u'}vH_{\mathsf v}(z_E)+Z_W,
 &\frac{u'}v
 &\le\frac12+\frac {3Q}x\,.
 \label{eq:right-renewal-drift}
\end{align}
(The bounds on the width ratios are due to $\sf{u},\sf{v}$ always being at most $Q$ away from $x/3, 2x/3$ respectively.) We refer to Figure~\ref{fig:oz-sweep} for a depiction of one East-West sweep and how it contracts $H_{\mathsf{u}}$ in expectation. 
Similarly, we can calculate the respective slopes under the East update and the subsequent West update as 
\begin{equation*}
 \rho-\frac{H_{\mathsf u}(z)}{x-u},
 \qquad
 \rho+\frac{H_{\mathsf v}(z_E)}v\,.
\end{equation*}
Since~$\mathfrak{S} = [-7/4, -1/4]$, height discrepancies of magnitude at most~$\delta x$ for $\delta$ small keep all slopes in a fixed compact subset of~$(-2,0)$.
If~$|H_{\mathsf u}(z)|\le\delta x$, the Gaussian tail bound of Lemma~\ref{lem:OZ-bridge-concentration}, namely~\eqref{eq:OZ-bridge-concentration}, applied on this neighborhood, gives~$\E[|Z_E|\mid z]\le C\sqrt x$ for a uniform constant $C$. Otherwise, the rough estimate~\eqref{eq:OZ-bridge-rough-moment} gives
\begin{equation*}
 \E[|Z_E|\mid z]
 \le C\sqrt{Qx}
 \le\frac{C}{\delta\sqrt L}|H_{\mathsf u}(z)|\,.
\end{equation*}
The same two estimates hold for the West update, conditionally on~$z_E$, with~$Z_E,H_{\mathsf u}(z)$ replaced by~$Z_W,H_{\mathsf v}(z_E)$. For large~$L$,~\eqref{eq:actual-renewal-drift}--\eqref{eq:right-renewal-drift} therefore imply, for a fixed~$\alpha<1$ and constant $C$,
\begin{equation}\label{eq:half-refresh-drift}
\begin{aligned}
 \E[|H_{\mathsf v}(z_E)|\mid z]
 &\le\alpha|H_{\mathsf u}(z)|+C\sqrt x\qquad \text{and} \qquad \E[|H_{\mathsf u}(z_1)|\mid z_E]
 &\le\alpha|H_{\mathsf v}(z_E)|+C\sqrt x\,.
\end{aligned}
\end{equation}
Combining the two updates yields constants~$\lambda<1$ and~$c_0<\infty$ such that after one sweep, 
\begin{align}\label{eq:sweep-drift}
 \E[V_{\mathsf u}(z_1)]
 &\le\lambda \E[V_{\mathsf u}(z)]+c_0,
\end{align}

\smallskip
\noindent \emph{Step 2: Coupling at diffusive height.}
We now wish to couple the process started from two initial states $z, z'$. 
Fix~$D<\infty$ and suppose~$V_{\mathsf u}(z)+V_{\mathsf u}(z')\le D$. Put~$u=\mathsf u(z), u' = {\mathsf{u}}(z')$. We claim that in an East update, these two can be coupled with constant probability on the suffix (that is, all renewal components to the right of $B_{\rm West}$). 

Let~$\mathsf q$ sample independent renewal components from the tilted and truncated distribution $\mathsf p_{\rho,Q}$ of~\eqref{eq:tilted-p-distribution}, backwards from~$(x,\rho x)$ until the first renewal~$(r,y)$ with~$r\le 2x/3$. (If the initial right endpoint of $b_{\Lambda_n}^*$ is the right of our domain, then $\mathsf{q}$ first samples from the normalized boundary renewal component measure $\bar\rho^R_{\rho,Q}$, then uses the independent increments of $\mathsf p_{\rho,Q}$.) 
Necessarily~$2x/3-Q\le r\le 2x/3$, and the sampled suffix includes the component crossing the right edge of $B_{\rm West}$. By the fact that the expected displacement of the internal renewal component is at slope $\rho$, and the uniform second moment bound of~\eqref{eq:OZ-bridge-rough-moment}, by Chebyshev, there is a fixed~$M$ such that
\begin{equation}\label{eq:prefix-central-event}
 E=\{|y-\rho r|\le M\sqrt x\},
 \qquad \mathsf q(E)\ge\tfrac12\,.
\end{equation}

Under the distribution when $B_{\rm East}$ is resampled in the two-block dynamics from state $z$, with left endpoint $(u,\phi_z(u))$, the Radon--Nikodym derivative of this suffix relative to~$\mathsf q$ is
\begin{equation}\label{eq:prefix-density}
 f(z)=
 \frac{U_{\rho,Q}(r-u,y-\phi_{z}(u))}
 {U^R_{\rho,Q}(x-u,\rho x-\phi_{z}(u))},
 \qquad
 U^R_{\rho,Q}=U_{\rho,Q}*\bar \rho^R_{\rho,Q}\,,
\end{equation}
where we recall $U_{\rho,Q}(x,y)$ from~\eqref{eq:OZ-renewal-function}. 
(The~$\bar \rho^R_{\rho,Q}$ is the law for the final boundary increment if the domain under consideration touches the right-edge of $\Lambda_n$, and is~$\delta_0$ otherwise.) This captures that given the pair of delineating renewal points, the laws are the same, so the Radon--Nikodym is only a function of the renewal points, and we can use local limit theorems to couple. On~$E$, the remaining width~$r-u$ is $x/3 - O(x) = \Theta(x)$, while
\begin{equation*}
 |y-\phi_{z}(u)-\rho(r-u)|
 \le(M+D)\sqrt x\,.
\end{equation*}
By the local limit theorem of Lemma~\ref{lem:OZ-local-limit}, specifically~\eqref{eq:OZ-local-limit}, so long as $Q\ge Q_0$, 
\begin{equation}\label{eq:prefix-common-mass}
 f(z)\wedge f(z')\ge\eta_D\quad\text{on }E\,,
\end{equation}
for a uniformly positive constant $\eta_D$. 
We therefore can take an optimal coupling of the suffixes of $z,z'$ under the $B_{\rm East}$ resampling using the common submeasure~$\eta_D\mathbf1_E\mathsf q$, and complete the remainder of the decorated interface in $B_{\rm East}$ independently, say, conditional on the respective suffixes. On success of this coupling, when $B_{\rm West}$ is updated, the boundary data provided by $z_{E}, z_{E}'$ are the same, and therefore the identity coupling on the West update gives 
\begin{align}\label{eq:diffusive-coupling-probability}
    \mathbb P(z_1 = z_1') \ge p_D:=\eta_D/2\,.
\end{align}

\smallskip
It remains to deduce the Poincar\'e inequality from the above and pass it to the continuous-time two-block dynamics. 
For distinct states $z,z'$, let~$S=V_{\mathsf u}(z)+V_{\mathsf u}(z')$. Given the individual height contractions of the interfaces when they are far from the line connecting $(0,0)\to (x,\rho x)$, and the constant coupling probability when they are only $O(\sqrt{x})$ from that line, we wish to show contraction of their difference in a suitable metric. Towards that, define the following distance  between $z$ and $z'$: 
\begin{align*}
    \mathsf{d}(z,z')=  (A_0 + V_{\sf u}(z) + V_{\sf u}(z'))\mathbf 1\{z\ne z'\}
\end{align*}
for a large constant $A_0$ to be determined. We show contraction under one sweep of East resampling followed by West resampling. Let $(Z_1,Z_1')$ be the outputs of that sweep, under a coupling that is arbitrary if $S > D_0$ and is the coupling of Step 2 that achieved~\eqref{eq:diffusive-coupling-probability} if $S\le D_0$.  
The contraction of $V_{\sf u}$~\eqref{eq:sweep-drift}, and the diffusive coupling probability of~\eqref{eq:diffusive-coupling-probability} imply
\begin{equation}\label{eq:coupling-decrease}
 \mathsf{d}(z,z')-\E[\mathsf d(Z_1,Z'_1)]\ge
 \begin{cases}
 (1-\lambda)S-2c_0,&S>D_0,\\
 p_{D_0}A_0+(1-\lambda)S-2c_0,&S\le D_0\,.
 \end{cases}
\end{equation}
Choose~$D_0>4c_0/(1-\lambda)$ and then~$A_0>4c_0/p_{D_0}$. Dividing by~$A_0+S$ gives for constant~$\varepsilon_0>0$, 
\begin{equation}\label{eq:weighted-contraction}
 \E[\mathsf d(Z_1,Z'_1)]
 \le(1-\varepsilon_0)\mathsf{d}(z,z')\,.
\end{equation}
Note that all constants here were independent of~$x,\rho,Q$.

Consider the positive self-adjoint operator~$T=P_{\rm East}P_{\rm West}P_{\rm East}$, where $P_{\rm East}$ and $P_{\rm West}$ are the resampling operators on $B_{\rm East}$ and $B_{\rm West}$. The last East update preserves~$V_{\mathsf u}$, and if $z= z'$ keeps them the same via the identity coupling, so extending the sweep coupling by one East update, preserves~\eqref{eq:weighted-contraction}. Hence, for any interface function $g$, if we define the seminorm  
\begin{equation*}
 \|g\|_{\mathsf{d}}=\max_{z\ne z'}
 \frac{|g(z)-g(z')|}{\mathsf{d}(z,z')},
 \qquad \text{ we have } \qquad 
 \|Tg\|_{\mathsf{d}}\le(1-\varepsilon_0)\|g\|_{\mathsf{d}}\,.
\end{equation*}
This contraction implies a spectral gap for the discrete-time chain $T$, $\text{Var}(g) \le \epsilon_0^{-1}\langle g,(I-T)g\rangle$ (see e.g., Theorem 13.1 of~\cite{LP}). Since $P_{\rm East},P_{\rm West}$ are orthogonal projections, 
\begin{align*}
 \langle g,(I-T)g\rangle
 &=\|g-P_{\rm East}g\|_2^2
   +\|(I-P_{\rm West})P_{\rm East}g\|_2^2 \le3\left(
   \E\var(g\mid\mathcal H_{B_{\rm East}})
   +\E\var(g\mid\mathcal H_{B_{\rm West}})
   \right)\,.
\end{align*}
The inequality follows by writing
\begin{equation*}
 (I-P_{\rm West})P_{\rm East}g
 =(I-P_{\rm West})g
  +(I-P_{\rm West})(P_{\rm East}g-g)\,,
\end{equation*}
 using $(a+b)^2 \le 2 a^2 + 2b^2$ and that projections are contractions. This proves~\eqref{eq:reference-two-window}, with the choice of constant~$C_0=3/\varepsilon_0$. 
\end{proof}

In order to go from this to a WPI for the interface dynamics of Definition~\ref{def:infinite-volume-interface-dynamics}, we recursively apply the previous lemma. 

\begin{proof}[Proof of Lemma~\ref{lem:infinite-volume-WPI}]
    Fix $Q$ large. We first condition on the $Q$-cutoff event $\mathcal B_Q$ (recall its definition from Definition~\ref{def:B_Q-set}) and write $\widehat\pi^\infty_{\Lambda,Q}$ for the conditional law. In what follows, suppose $L$ is sufficiently large, and set $W_Q=  LQ $.
We will continue recursing down in width $x$ so long as the slope~$\rho$ satisfies
\begin{equation}\label{eq:recursive-slope-interval}
 \operatorname{dist}\bigl(\rho,[-3/2,-1/2]\bigr)\le B\sqrt{Q/x}\,.
\end{equation}
For a fixed $B$ as long as~$L$ is sufficiently large, then always~$B/\sqrt{L}<1/4$ and $\rho \in \mathfrak{S}$.  

We construct a tree representing the distributions we encounter along the recursion. The root bridge has distribution~$\widehat\pi^\infty_{\Lambda,Q}$. The tree is constructed recursively as follows. For a parent bridge of width~$x>LQ$ whose slope satisfies~\eqref{eq:recursive-slope-interval}, take its endpoints as~$a,b$ in Definition~\ref{def:block-interface-dynamics} and for each of~$B\in\{B_{\rm West},B_{\rm East}\}$ and every possible realization~$h$ of~$\mathcal H_B$, introduce a child carrying the parent's law conditioned on $\mathcal H_B =h$. For each vertex $w$ on the tree, we use $\nu_{w}$ to denote this conditional distribution. 
This is again a distribution of the form $\widehat{\pi}_{a_h,b_h,Q}^\infty$ with some width and some endpoints $a_h,b_h$ determined by $h$. If~$\nu_{p(w)}$ denotes its parent's law, we also give this child \emph{weight}~$\nu_{p(w)}(h)$. By applying~\eqref{eq:reference-two-window} of Lemma~\ref{lem:two-block-interface-constant-gap}, for any parent, 
\begin{equation}\label{eq:variance-tensorization-across-block-dynamics}
 \var_\nu(g)
 \le (2C_0) \cdot \frac{1}{2}\sum_{B\in\{B_{\rm West},B_{\rm East}\}}
       \sum_h \nu(h)\var_{\nu(\cdot\mid h)}(g)\,.
\end{equation}
We continue building out the tree, except when we encounter either a vertex whose width is at most $L Q$, or a vertex whose domain has width $x$ but whose slope $\rho$ violates~\eqref{eq:recursive-slope-interval} (in which case we will crudely bound its conditional variance by $\Osc(g)^2$). 

We can iterate~\eqref{eq:variance-tensorization-across-block-dynamics} down the tree, representing the final result as an average of products of tokens accumulated by a random walk directed down the tree where at each step, one of East or West is chosen with probability $1/2$ each, and then one of the $h$ realizations available at the parent are chosen with probability $\nu(h)$. The random walk $(X_{k})_{k\ge 0}\sim \mathbf{P}$ starts at the root and picks up $C_0$ at every step it descends down the tree. Letting $\mathcal L_{\text{good}}$ denote the leaves whose domain widths are at most $LQ$, and assigning the conditional law at $\ell\in \mathcal L_{\text{good}}$ to $\nu_\ell$, the random walk picks up $\Var_{\nu_\ell}(g)$ if it ends at $\ell \in \mathcal L_{\text{good}}$ and picks up $\Osc(g)^2$ if it ends at a bad vertex. Let $X_\infty$ be its ending vertex.  Note that at any stage, if the parent had domain width $x$, each child has domain width~$x'=2x/3+O(Q)$, which is necessarily in $[x/2,3x/4]$ as long as $L$ is large. Therefore, the tree has depth at most $\log_{4/3}n$. As a result, if the random walk expectation is denoted by $\mathbf{E}$, 
\begin{align}
 \var_{\widehat\pi^\infty_{\Lambda,Q}}(g)
 &\le (2C_0)^{\log_{4/3} n}
 \mathbf E[
   \mathbf 1_{\{X_\infty \in \mathcal L_{\rm good}\}}\var_{\nu_{X_\infty}}(g)
   +\mathbf 1_{\{X_\infty \notin \mathcal L_{\rm good}\}}\Osc(g)^2
 ]\notag
\end{align}
In order for $X_\infty \notin \mathcal L_{\rm good}$, at some point the walker went from a vertex which had domain-width $x$ and slope satisfying~\eqref{eq:recursive-slope-interval} to child with domain-width $x'\in [x/2,3x/4]$ but slope not satisfying~\eqref{eq:recursive-slope-interval}, and in particular that means the child's slope differs from its parents' slope by at least $B(\sqrt{Q/x'} - \sqrt{Q/x})$. That is to say that for that child, its renewal points nearest $x/3$ or $2x/3$ deviated from the chord connecting the parents' boundary points by at least $c B \sqrt{Qx}$, and the probability of this is at most $C e^{ - cQ}$ by the Gaussian tail bound of~\eqref{eq:OZ-bridge-concentration}. Altogether, we obtain for some constants $c,C$, for all $Q\ge Q_0$, 
\begin{align}\label{eq:WPI-recursion-interface}
 \var_{\widehat\pi^\infty_{\Lambda,Q}}(g)
 &\le n^C
 \mathbf{E}[\mathbf 1_{\{X_\infty \in \mathcal L_{\rm good}\}}\var_{\nu_{X_\infty}}(g)]
   + n^C e^{ -  c Q}\Osc(g)^2
 \,.
\end{align}
We finally want to compare the variances on the random domains of size at most $LQ$ given by the leaves, to the fixed windows $J_j = \{j,...,j+w\}\times \{1,...,n\}$ of Definition~\ref{def:infinite-volume-interface-dynamics}. 
If a leaf $\ell$ has boundary columns at $u_\ell,v_\ell$, choose a  window~$J_j$ such that $\{u_\ell,...,v_\ell \}\subset \{j,...,j+w\}$ (possible because $w\ge LQ$). 
For each window $J_j$, let $G_j$ be the event that the leaf $X_\infty$ stops at gives window $J_j$ and note that $\mathcal H_{J_j}$ is already conditioned on in the conditional law $\nu_{X_\infty}$. We thus write 
\begin{align*}
    \mathbf{E}[\mathbf 1_{\{X_\infty \in \mathcal L_{\rm good}\}}\var_{\nu_{X_\infty}}(g)] & = \sum_{j=0}^{n-w} \mathbf E[\mathbf 1_{\{G_j\}} \mathbf 1_{\{X_\infty \in \mathcal L_{\rm good}\}} \text{Var}_{\nu_{X_\infty}}(g)] \\
    & \le \sum_{j=0}^{n-w} \widehat \pi_{\Lambda,Q}^\infty[\text{Var}(g \mid \mathcal H_J)]\,,
\end{align*}
where the inequality was using that on each good leaf $\ell$ assigned to $J_j$, the exterior $\mathcal H_{J_j}$
is fixed, so $m_j:=\widehat\pi^\infty_{\Lambda,Q}[g\mid\mathcal H_{J_j}]$
is constant there and $\var_{\nu_\ell}(g)\le\nu_\ell[(g-m_j)^2]$.
Averaging over these leaves and then dropping their indicator bounds the
$j$th contribution by $\widehat\pi^\infty_{\Lambda,Q}[(g-m_j)^2]$. Combined with~\eqref{eq:WPI-recursion-interface}, we get an $(n^C, n^Ce^{ - cQ})$-WPI for the infinite-volume interface dynamics conditioned on the $Q$-cutoff. This transfers to an  $(n^C, n^Ce^{ - cQ} + C n^{3/2} e^{ - cQ})$-WPI for the infinite-volume interface dynamics without the conditioning via Lemma~\ref{lem:conditioning-wpi} and Lemma~\ref{lem:OZ-exponential-tails}. 
\end{proof}

\subsection{Transferring to finite-volume interface dynamics}

Recall the finite-volume interface dynamics from Definition~\ref{def:interface-dynamics-finite-volume}. Our aim is to use comparison tools, and the fact that the laws $\bar \pi^\infty$ and $\bar \pi$ on interfaces $\Gamma$ in $\Lambda_n$ have comparable Radon--Nikdoym derivative (at least for nice interfaces) to transfer the WPI from Lemma~\ref{lem:infinite-volume-WPI}.

We begin with a comparison lemma for Poincar\'e inequalities and weak Poincar\'e inequalities between the infinite-volume interface dynamics and the finite-volume dynamics. If the resampling laws had a bounded Radon--Nikodym derivative, this would be completely standard, but in this situation, they resample conditional on different sigma-algebras ($\mathcal F_J$ vs.\ $\mathcal H_J$) and the Radon--Nikodym derivative is only bounded on a high probability set of interfaces.

\begin{lemma}\label{lem:good-set-comparison}
Fix width~$w$, and consider both the finite and infinite-volume interface dynamics from Definitions~\ref{def:interface-dynamics-finite-volume} and \ref{def:infinite-volume-interface-dynamics}. 
Suppose the infinite-volume interface dynamics satisfies an~$(A,\delta)$-WPI. Let~$\psi_\Lambda(\Gamma)$ be the ratio of weights $q_\Lambda(\Gamma)/q_\infty(\Gamma)$ and let~$Z=\bar\pi^\infty[\psi_\Lambda]$. Then for every~$K\ge1$, all functions $g$ of the interface $\Gamma$, the finite-volume interface dynamics has 
\begin{equation}\label{eq:good-set-comparison}
 \var_{\bar\pi}(g)
 \le KA\mathcal E_{\rm Int}(g,g)
 +\left(\frac{K\delta}{Z}+\bar\pi(\psi_\Lambda(\Gamma) >K)\right)\Osc(g)^2\,.
\end{equation}
\end{lemma}

\begin{proof}
Fix a function $g$ depending only on the interface $\Gamma$. 
For each block~$J$, the finite-volume conditional expectation 
$m_J=\bar\pi[g\mid\mathcal F_J]$ is determined by~$\mathcal H_J$.
Using the fact that $\text{Var}(f) = \mathbb E[(f-\E[f])^2]\le \E[(f-c)^2]$ for all $c\in \mathbb R$, we obtain
\[
 \widehat\pi^\infty[\var(g\mid\mathcal H_J)] 
 \le \bar\pi^\infty[(g-m_J)^2]
 \le Z\bar\pi[\var(g\mid\mathcal F_J)]\,,
\]
where the last inequality used 
$d\bar\pi/d\bar\pi^\infty=\psi_\Lambda/Z$ and~$\psi_\Lambda\ge1$ from Lemma~\ref{lem:OZ-finite-volume-comparison}.
Since both chains update each $(J_j)_{j=0}^{n-w}$ at rate-one, summing over~$j$ gives the Dirichlet form comparison 
\begin{align}\label{eq:interface-Dirichlet-form-comparison}
\widehat{\mathcal E}_{\rm Int}(g,g)
 \le Z\mathcal E_{\rm Int}(g,g)\,.\end{align}

On the event~$\{\psi_\Lambda(\Gamma) \le K\}$
the density $d\bar \pi/d\bar \pi^\infty$ is at most~$K/Z$. On its complement
we can bound the contribution to the variance by $\Osc(g)^2$. Thus, 
\[
 \var_{\bar\pi}(g)\le \frac KZ\var_{\bar\pi^\infty}(g)
     +\bar\pi(\psi_\Lambda>K)\Osc(g)^2\,.
\]
Since $\var_{\bar\pi^\infty}(g)=\var_{\widehat\pi^\infty}(g)$ (as $g$ does not depend on the decorations),
applying the infinite-volume interface dynamics' $(A,\delta)$-WPI and then~\eqref{eq:interface-Dirichlet-form-comparison} proves the claim.
\end{proof}

We can now apply Lemma~\ref{lem:good-set-comparison} to deduce Lemma~\ref{lem:interface-WPI} from Lemma~\ref{lem:infinite-volume-WPI}.

\begin{proof}[Proof of Lemma~\ref{lem:interface-WPI}]\label{sec:marginal}
Let $b,c$ be the constants in Lemma~\ref{lem:OZ-finite-volume-comparison}, and let $c_0$ be the decay rate $c$ multiplying $Q$ in Lemma~\ref{lem:infinite-volume-WPI}. Choose a fixed $k$ with $c_0k>b+c$. Take $R\ge R_0$, where $R_0$ is large enough that $Q= kR\ge Q_0$, and take $w= L_0R$ with $L_0$ sufficiently large. The infinite-volume interface dynamics' WPI from Lemma~\ref{lem:infinite-volume-WPI} has coefficients
\[
    A\le n^p,\qquad \delta\le n^pe^{-c_0Q}\,.
\]
Here and below $p$ is a constant that may change from line to line. 

We now apply Lemma~\ref{lem:good-set-comparison} to the $(A,\delta)$-WPI of Lemma~\ref{lem:infinite-volume-WPI}. Let $Z=\widehat\pi^\infty[\psi_\Lambda]\ge1$ and  set $K=e^{bR}$. Lemma~\ref{lem:OZ-finite-volume-comparison} bounds the probability $\bar \pi^{\infty}(\psi_{\Lambda}(\Gamma)>e^{bR})$ by $Cn^{5/2}e^{-cR}$, so we deduce an $(A',\delta')$-WPI for the finite-volume interface dynamics for $A' \le n^p e^{bR}$ and 
\[
    \delta' \le K\delta/Z+\bar\pi(\psi_\Lambda(\Gamma)>K)\le n^{p'}e^{-cR}\,.
\]
Taking $R = K'\log n$ for $K'$ sufficiently large proves the lemma. 
\end{proof}

\subsection{WPI with respect to the Glauber Dirichlet form}

The aim of this section is to transfer the weak Poincar\'e inequality we derived for the interface dynamics to one where the right-hand side entails the Glauber Dirichlet form rather than the interface resampling Dirichlet form. This will come only at a polynomial cost since the resampling strips are of width at most $O(\log n)$. 

\begin{definition}\label{def:block-Glauber}
    For width-parameter $1\le w\le n$, the spin \emph{block dynamics} is the Markov chain that performs continuous-time heat-bath resampling updates of the full spin-configuration with blocks $S\in \mathcal S$ where $\mathcal S = \{\{r,...,s\}\times \{1,...,n\}: 1\le r<s\le n\,; s-r\le w\}$ is the set of all possible strips of width at most $w$. We let $P_{\rm Gl-Block}$ and $\mathcal E_{\rm Gl-Block}$ be its semigroup and Dirichlet form. 
\end{definition}

 We first go from the interface dynamics to the Dirichlet form for the above spin \emph{block dynamics}.

\begin{lemma}\label{lem:pinned-vertices-Dirichlet-form-comparison}
Fix a width~$w$ for the interface dynamics and use width parameter~$w+4$
for the spin block dynamics. For truncation level~$L\ge0$, set
$a=(1+e^{8\beta})^{-C(L+1)}$ for a sufficiently large constant~$C$.
Then, for every spin function~$f:\{-1,+1\}^{\Lambda_n}\to \mathbb R$ with conditional expectation 
$g(\Gamma)=\pi[f\mid\Gamma]$, the corresponding Dirichlet forms satisfy
\begin{equation}\label{eq:cylinder-lift}
 \mathcal E_{\mathrm{Int}}(g,g)
 \le3 n \Big(2a^{-1}\bar \pi[\text{Var}(f\mid \Gamma)]+a^{-2}\mathcal E_{\mathrm{Gl\text{-}Block}}(f,f)
       +n^3 e^{c_1 w - c_2 L}\Osc(f)^2\Big)\,.
\end{equation}
\end{lemma}

Lemma~\ref{lem:pinned-vertices-Dirichlet-form-comparison} is a comparison between the Dirichlet form of (finite-volume) interface dynamics and that of the spin block dynamics. For any fixed block $J$ for the interface dynamics and exterior increments $\mathcal F_J$, there is a strip $S \in \mathcal S$ which aligns with the two outermost cone-points in $J$ delineating the middle portion $\eta_J$ which gets resampled. The remaining difference is captured by the Radon--Nikodym derivative on spin configurations given all of $\Gamma$ or all of $\Gamma \setminus \eta_J$, which is easily seen to be exponential in $|\eta_J|$. The following Lemma~\ref{lemma:short-middle} bounds the probability that the middle portion consists of many more edges $L$ than the strip-width $w$; this bad event can be absorbed by the oscillation error in~\eqref{eq:cylinder-lift}. 

\begin{lemma}\label{lemma:short-middle}
Fix any $w$ large, fix a block $J = \{j,...,j+w\}\times \{1,...,n\}$  and let~$\eta_J$ be the \emph{middle portion} of the interface $\Gamma$ being resampled under Definition~\ref{def:interface-dynamics-finite-volume}. There are constants~$C_\beta,c_1,c_2>0$, such that 
\begin{equation}\label{eq:middle-tail}
 \bar\pi(|\eta_J|>L)\le\varepsilon_{n,w,L},\qquad\text{where}\qquad 
 \varepsilon_{n,w,L}:=C_\beta n^{5/2} e^{c_1w-c_2L}\,.
\end{equation}
\end{lemma}

\begin{proof}
Write the interface~$\Gamma=\alpha\circ\eta\circ \omega$ where $\alpha$ is all increments to the left of and including the one through $\{j\}\times \{1,...n\}$ and $\omega$ is all increments to the right of and including the one through $\{j+w\}\times \{1,...,n\}$ and $\eta$ is the \emph{middle portion}. Let~$a,b$ be the endpoints of~$\eta$. The cone constraints at the cone-points give
$0\le b_1-a_1\le w$ and $-2(b_1-a_1)\le b_2-a_2\le0$
so~$\|b-a\|_1\le 3w$. 

For fixed~$a,b$, put~$H=\Lambda\setminus Q_\alpha$ and~$K=H\setminus Q_\eta$ where we recall that $Q_\alpha, Q_\eta$ are the sets of dual edges reserved by the increments of $\alpha,\eta$ respectively. The Markov property gives~$q_\Lambda(\Gamma)=q_\Lambda(\alpha)q_H(\eta)q_K(\omega)$.  Take~$h=h(-1)$. Since~$h\cdot x\le\xi(x)$, the bound of Lemma~\ref{lem:raw-source-sum} gives
\begin{equation*}
  \sum_{\alpha:a_\Lambda^*\to a}e^{h\cdot D(\alpha)}q_\Lambda(\alpha)\le C_\beta\,,
 \qquad \sum_{\omega:b\to b_\Lambda^*}e^{h\cdot D(\omega)}q_K(\omega)\le C_\beta\,.
\end{equation*}
uniformly in $\eta$.  The bound of Lemma~\ref{lem:raw-contour-tail} and~$\|b-a\|_1\le3w$ give uniformly over $\alpha$, that 
\begin{equation*}
 \sum_{\substack{\eta:a\to b\\|\eta|>L}}
 e^{h\cdot D(\eta)}q_H(\eta)
 \le C_\beta e^{C_\beta w-c_\beta L}\,.
\end{equation*}
Summing first over~$\omega$, then~$\eta$, then~$\alpha$, there are at most~$(n+1)^2$ choices for~$a$ and~$Cw^2$ for~$b-a$. Dividing by the normalization constant of~\eqref{eq:physical-reference-normalizer} and absorbing~$w^2$ into~$e^{c_1w}$ we obtain~\eqref{eq:middle-tail}.
\end{proof}

\begin{proof}[Proof of Lemma~\ref{lem:pinned-vertices-Dirichlet-form-comparison}]
Recall the blocks $J_j$ of width $w$ for the interface dynamics of
Definition~\ref{def:interface-dynamics-finite-volume}, and its Dirichlet form
\begin{align*}
    \mathcal E_{\rm Int}(g,g)
    = \sum_j \bar\pi[\operatorname{Var}_{\bar\pi}(g\mid\mathcal F_{J_j})]\,.
\end{align*}
If $A_j=\{|\eta_{J_j}|\le L\}$ is the event that the middle portion
$\eta_{J_j}$ has length at most $L$, then we can condition on that event
at the expense of an oscillation term. Namely, by the variance bound in
the proof of Lemma~\ref{lem:conditioning-wpi},
\begin{align*}
    \mathcal E_{\rm Int}(g,g)
    &\le \sum_j \bar\pi\big[
        \bar\pi(A_j\mid\mathcal F_{J_j})
        \operatorname{Var}_{\bar\pi}(g\mid\mathcal F_{J_j},A_j)
        +\bar\pi(A_j^c\mid\mathcal F_{J_j})\operatorname{osc}(g)^2
    \big] \\
    &\le \sum_j \bar\pi\big[
        \bar\pi(A_j\mid\mathcal F_{J_j})
        \operatorname{Var}_{\bar\pi}(g\mid\mathcal F_{J_j},A_j)
    \big]
    +n\max_j\bar\pi(A_j^c)\operatorname{osc}(g)^2\,.
\end{align*}
Now fix a block $J=J_j$. For $\Gamma,\Gamma'$ with the same realization $h$ of the information revealed by $\mathcal F_J$, define
\begin{align*}
    M_j(\Gamma,\Gamma')
    = \bar\pi(h)\bar\pi(A_j\mid h)
      \bar\pi(\Gamma\mid h,A_j)\bar\pi(\Gamma'\mid h,A_j)\,.
\end{align*}
Set $M_j(\Gamma,\Gamma')=0$ for pairs of $\Gamma,\Gamma'$ with different
exteriors $\mathcal F_J$. Using this to write out the summands above, and applying
Lemma~\ref{lemma:short-middle} to the last term, gives
\begin{align}\label{eq:int-Dirichlet-form-to-Glauber}
    \mathcal E_{\rm Int}(g,g)
    \le \frac12\sum_j\sum_{\Gamma,\Gamma'}M_j(\Gamma,\Gamma')
         (g(\Gamma)-g(\Gamma'))^2
         +n\varepsilon_{n,w,L}\operatorname{osc}(g)^2\,.
\end{align}
Fix a pair $\Gamma,\Gamma'$ contributing to the above sum, and let $u,v$ be the common endpoints
of their two middle portions. Take $S$ to be the $(v_1 - u_1+4)$-width strip between
the columns of $u,v$, enlarged by two columns on each side and intersected
with $\Lambda_n$. Every middle portion $\eta$ compatible with exterior $h$ lies between the two
endpoint columns, so by~\eqref{eq:collar-size} all the spins it gives information about (i.e., $\partial_v \eta$) lie in $S$.
Outside $S$, $h$ prescribes the same sites and the
same spins, determined by tracing the common prefix back from $u$ and
the common suffix forward from $v$.

Inside $S$, each of the two middle portions prescribe $O(L)$ spins
by~\eqref{eq:collar-size}. Due to the cones, the prefix before $u$ and suffix after $v$ only prescribe $O(1)$ many sites in $S$. Therefore the full information of $\Gamma, \Gamma'$ only prescribes $O(L)$ many spins in the strip $S$. 

We now define a resampling procedure for spin configurations given interfaces that interpolates between (a) performing finite-volume interface dynamics and then drawing spins conditional on the interfaces, and (b) a spin resampling on the strip $S$ given the spins on $\Lambda_n \setminus S$. 
Let $\mathbf P$ be the joint law on $(\sigma,\sigma')$ sampled as follows. 
First draw the common configuration
$\sigma_{\Lambda_n\setminus S}
=\sigma'_{\Lambda_n\setminus S}$ from $\pi(\cdot\mid h)$ on
$\Lambda_n\setminus S$, then sample 
\[
    (\sigma_S,\sigma'_S)\sim
    \pi(\cdot\mid \sigma_{\Lambda_n\setminus S},\Gamma)
    \otimes\pi(\cdot\mid \sigma_{\Lambda_n\setminus S},\Gamma')\,.
\]
Note that given $h$, the knowledge of $\Gamma$ is equivalent to the prescribing of the $O(L)$ spins described above in $S$ and therefore, forcing those spins iteratively, each with probability at least $(1+e^{8\beta})^{-1}$,
gives $\pi(\Gamma\mid \sigma_{\Lambda_n\setminus S})\ge  (1+e^{8\beta})^{-C(L+1)}= a$, and likewise for $\Gamma'$.

Given $\sigma_{\Lambda_n\setminus S}$, the first marginal of $\mathbf P$ and $\pi(\cdot\mid\Gamma)$
have the same law. Their Radon--Nikodym derivative is therefore the ratio
of the probabilities they assign to $\sigma_{\Lambda_n\setminus S}$. By Bayes' rule and since $\Gamma$
determines $h$, 
\[
    \frac{\pi(\sigma_{\Lambda_n\setminus S}\mid h)}{\pi(\sigma_{\Lambda_n\setminus S}\mid\Gamma)} = \frac{\pi(h\mid \sigma_{\Lambda_n\setminus S})\,\pi(\sigma_{\Lambda_n\setminus S})}{\bar\pi(h)} \frac{\bar\pi(\Gamma)}{\pi(\Gamma\mid \sigma_{\Lambda_n\setminus S})\,\pi(\sigma_{\Lambda_n\setminus S})}
    =\frac{\pi(h\mid \sigma_{\Lambda_n\setminus S})\bar\pi(\Gamma\mid h)}{\pi(\Gamma\mid \sigma_{\Lambda_n\setminus S})}
    \le a^{-1}.
\]
Applying the same bound to the second marginal gives
\begin{align}\label{eq:Radon-Nikodym-interface-spin}
    \frac{\mathbf P(\sigma\in\cdot)}{\pi(\cdot\mid\Gamma)}
    \vee
    \frac{\mathbf P(\sigma'\in\cdot)}{\pi(\cdot\mid\Gamma')}
    \le a^{-1}\,.
\end{align}

Now recalling $g(\Gamma)=\pi[f\mid\Gamma]$, we can write
\begin{align*}
    (g(\Gamma)-g(\Gamma'))^2
    \le 3(f(\sigma)-g(\Gamma))^2
       +3(f(\sigma)-f(\sigma'))^2
       +3(g(\Gamma')-f(\sigma'))^2\,.
\end{align*}
Averaging $(\sigma,\sigma')\sim\mathbf P$, we have
by~\eqref{eq:Radon-Nikodym-interface-spin}
\begin{align*}
    \mathbf E[(f(\sigma)-g(\Gamma))^2]
    \le a^{-1}\pi[(f-g(\Gamma))^2\mid\Gamma]
    =a^{-1}\operatorname{Var}_\pi(f\mid\Gamma)\,,
\end{align*}
and similarly for $\Gamma'$. The marginals of $M_j$
are $\bar\pi(\Gamma)\mathbf1_{A_j}(\Gamma)$, so averaging either bound
against $M_j$ gives at most
$a^{-1}\bar\pi[\operatorname{Var}_\pi(f\mid\Gamma)]$.

It remains to bound the middle term $\mathbf E[(f(\sigma)-f(\sigma'))^2]$
by the spin-block Dirichlet form. For every $\sigma,\sigma'$ compatible
with $\Gamma,\Gamma'$ and with common exterior $\sigma_{\Lambda_n\setminus S}$,
\begin{align*}
    \frac{M_j(\Gamma,\Gamma')\mathbf P(\sigma,\sigma')}
         {\pi(\sigma_{\Lambda_n\setminus S})\pi(\sigma_S\mid \sigma_{\Lambda_n\setminus S})\pi(\sigma'_S\mid \sigma_{\Lambda_n\setminus S})}
    = 
       \frac{M_j(\Gamma,\Gamma')\,\pi(h\mid \sigma_{\Lambda_n\setminus S})}
            {\pi(h)\pi(\Gamma\mid \sigma_{\Lambda_n\setminus S})\pi(\Gamma'\mid \sigma_{\Lambda_n\setminus S})}
    \le a^{-2}\,,
\end{align*}
since $M_j(\Gamma,\Gamma')\le\bar\pi(h)=\pi(h)$.
For fixed $j$ and $S$, each pair $(\sigma,\sigma')$ determines $\Gamma,\Gamma'$, so it is counted at most once in the sum over
$\Gamma,\Gamma'$. Summing the above first over the compatible $(\sigma,\sigma')$ and then
over $S\in\mathcal S$ therefore gives
\[
    \frac12\sum_{\Gamma,\Gamma'}M_j(\Gamma,\Gamma')
    \mathbf E[(f(\sigma)-f(\sigma'))^2]
    \le a^{-2}\mathcal E_{\rm Gl\text{-}Block}(f,f)\,.
\]
Plugging this bound into~\eqref{eq:int-Dirichlet-form-to-Glauber},
summing over at most $n$ choices of $j$, and using
$\operatorname{osc}(g)\le\operatorname{osc}(f)$, 
\begin{align*}
    \mathcal E_{\rm Int}(g,g)
    \le 
        3n\Big(a^{-1}\bar\pi[\operatorname{Var}_\pi(f\mid\Gamma)]
          +a^{-2}\mathcal E_{\rm Gl\text{-}Block}(f,f)\Big)
        +n\varepsilon_{n,w,L}\operatorname{osc}(f)^2\,.
\end{align*}
Equation~\eqref{eq:middle-tail} now gives~\eqref{eq:cylinder-lift}.
\end{proof}

Lastly, we have to go from the spin block dynamics' Dirichlet form to the true single-site Glauber Dirichlet form. This is a standard comparison which incurs only a polynomial cost when the width of the blocks is $w = O(\log n)$, due to the cut-width bound of Lemma~\ref{lem:cutwidth-mixing-bound}.   

\begin{lemma}\label{lem:Glauber-block-comparison}
    Fix any $w\ge 1$. For every $\beta>0$, there exists a constant $A_\beta$ such that 
    $$\mathcal E_{\rm Gl-Block}(f,f) \le A_\beta n e^{A_\beta w} \mathcal E_{\rm Gl} (f,f) $$
\end{lemma}

\begin{proof}
    Every strip $S \in \mathcal S$ contains at most~$n w$ sites and has cutwidth at most~$w$. Apply Lemma~\ref{lem:cutwidth-mixing-bound} to get for each strip $S = \{r,...,s\}\times \{1,...,n\}$, and any boundary conditions $\tau$ which are the concatenation of $\sigma(\Lambda_n\setminus S)$ with the diagonal Dobrushin boundary of $\Lambda_n$, 
    \begin{align*}
        \text{Var}_{\pi_{S}^\tau}(f) \le 4(1+e^{8\beta}) nw e^{4\beta w} \mathcal E_{{\rm Gl},S^\tau}(f,f)
    \end{align*}
    where $\mathcal E_{{\rm Gl},S^\tau}$ is the Dirichlet form for Glauber dynamics on $S$ with boundary conditions $\tau$. 
    Summing over $S\in \mathcal S$, and taking expectation with respect to $\pi$, we get exactly the Dirichlet form $\mathcal E_{\rm Gl-Block}$ on the left-hand side, and so 
    \begin{align*}
        \mathcal E_{\rm Gl-Block}(f,f) \le 4(1+e^{8\beta}) nw e^{4\beta w} \sum_{S\in \mathcal S} \E_{\pi}[ \mathcal E_{{\rm Gl},S^\tau}(f,f)]\,.
    \end{align*}
    The right-hand side is summing over squared gradients of $f$ for single-spin resamplings, and each vertex in $\Lambda_n$ is counted at least once and at most $w^2$ many times. Thus, the right-hand side is at most $4(1+e^{8\beta}) n w^3 e^{4\beta w}\mathcal E_{\rm Gl}(f,f)$ as desired.  
\end{proof}

\section{Mixing with Dobrushin boundary conditions}\label{sec:mixing-Dobrushin-bdy}

In this section, we combine the above weak Poincar\'e inequalities to deduce the polynomial WPI of Lemma~\ref{lem:diagonally-Dobrushin-WPI} for the square with diagonal Dobrushin boundary conditions. We then explain why the same proof extends to minor modifications of the geometry, and use those WPIs together with a censoring scheme to continually ensure warm-starts, to deduce Theorem~\ref{thm:diagonally-Dobrushin-mixing}.

\subsection{WPI with Dobrushin boundary conditions}

\begin{proof}[Proof of Lemma~\ref{lem:diagonally-Dobrushin-WPI}]
    We begin by putting the interface WPI and off-interface WPIs together. Fix any function $f$ and let $g(\Gamma) = \mathbb {\pi}[f \mid \Gamma]$. Recall the decomposition by the law of total variance from~\eqref{eq:intro-total-variance}. 
    We begin with bounding the second term $\text{Var}_{\bar \pi}(g)$. By Lemma~\ref{lem:interface-WPI}, for every $C'$, we can choose constants $C,C''$ large enough such that if $w= C'' \log n$, one has 
    \begin{align*}
        \text{Var}_{\bar \pi}(g) \le n^C \mathcal E_{\rm Int}(g,g) + n^{- C'}\Osc(g)^2\,.
    \end{align*}
    By Lemma~\ref{lem:pinned-vertices-Dirichlet-form-comparison}, if we then take $L = C_L \log n$ for a sufficiently large constant $C_L$ (depending on $C''$), followed by Lemma~\ref{lem:Glauber-block-comparison}, we see that for any $C'$, there exists $C$ such that 
    \begin{align*}
        \text{Var}_{\bar \pi}(g) \le n^C \bar \pi[\text{Var}_\pi(f\mid \Gamma)] + n^C \mathcal E_{\rm Gl}(f,f) + n^{-C'} \Osc(f)^2 + n^{-C'}\Osc(g)^2\,.
    \end{align*}
    Noting that $\Osc(g) \le \Osc(f)$, and plugging this into~\eqref{eq:intro-total-variance}, we obtain 
    \begin{align}\label{eq:diagonal-WPI-almost-there}
        \text{Var}_{\pi}(f) \le (1+n^C) \bar \pi[\text{Var}_\pi(f\mid \Gamma)] + n^C \mathcal E_{\rm Gl}(f,f) +  2 n^{-C'} \Osc(f)^2\,.
    \end{align}
    By Corollary~\ref{cor:given-Gamma-WPI}, for every $C'$ there exists $C$ such that 
    \begin{align*}
        \bar \pi[\text{Var}_\pi(f\mid \Gamma)] \le n^C \mathcal E_{\rm Gl}(f,f) + n^{ - C'}\Osc(f)^2\,,
    \end{align*}
    and plugging that in to~\eqref{eq:diagonal-WPI-almost-there} concludes the proof, up to adjusting the constant $C$. 
\end{proof}

\subsection{Modifications to above arguments for near-squares}\label{sec:near-squares}

Recall from the proof sketch Section~\ref{sec:WPI-to-gap}, we need the polynomial WPI of Lemma~\ref{lem:diagonally-Dobrushin-WPI} for a sequence of near-square domains that grow to $\Lambda_n$ one vertex at a time. 

\begin{definition}\label{def:near-squares}
    A \emph{near-square} in $\Lambda_n$ is a domain of the form $\{1,...,k\}^2 \cup (\{1,...,j\}\times \{k+1\})$ or $(\{1,...,k\}\times \{1,...,k+1\}) \cup \{k+1\}\times \{j,...,k+1\}$ for $1\le j\le k\le n$. 
    
    The near-squares in $\Lambda_n$ can be enumerated in growing manner as  $\Lambda^{(1)} \subset \cdots \subset \Lambda^{(n^2)}$ in the following way: for each $1\le k\le n$, let $\Lambda^{(k^2)} = \{1,...,k\}^2$ be the square subset of $\Lambda_n$ of side-length $k$ in its bottom-left corner. For $j = 1,...,2k+1$, obtain $\Lambda^{(k^2 +j)}$ from  $\Lambda^{(k^2 +j-1)}$ by adding the ``next" available vertex in the L-shaped strip $\Lambda^{((k+1)^2)}\setminus \Lambda^{(k^2)}$, where the enumeration is by starting at $(0,k+1)$, going east to $(k+1,k+1)$ then south to $(k+1,0)$. Recall Figure~\ref{fig:diag-Dobrushin-censoring}. 
\end{definition}

Note that each \emph{near-square} $\Lambda^{(k^2+j)}$ contains $\Lambda_k$ and is contained in $\Lambda_{k+1}$, and furthermore if the configuration in $\Lambda_n\setminus \Lambda^{(k^2+j)}$ is fixed to all-plus, then the diagonally Dobrushin boundary on $\Lambda_n$ induces diagonally Dobrushin boundary on $\Lambda^{(k^2+j)}$. Our aim is to reason that with minimal geometric modifications, our previous arguments also give the polynomial weak Poincar\'e inequality on the near-squares in $\Lambda_n$. 
The reason for that is that at constant cost in the WPI, we can pin corner vertex spins so that the interface in the near-square $\Lambda^{(k^2 + j)}$ is necessarily contained in the square $\Lambda^{(k^2)}$. 
The following lemma gives a comparison when freezing $O(1)$ spins to certain values. 

\begin{lemma}\label{lem:remove-pins-wpi}
Fix any~$\beta>0$ and~$k\ge1$. For an Ising distribution~$\pi$ on  finite~$U\subset\mathbb Z^2$ with boundary condition $\tau$,  let~$\nu=\pi(\cdot\mid\sigma_A=\xi)$ for any $\xi$ and any $A\subset U$ with $|A|\le k$. There is a constant~$C_{\beta,k}<\infty$, such that the rate-one single-site Dirichlet forms satisfy 
\begin{equation}\label{eq:pin-comparison}
 \var_\pi(f)\le C_{\beta,k}\bigl[\var_\nu(f)+\mathcal E_{\mathrm{Gl}}^\pi(f,f)\bigr],
 \qquad \mathcal E_{\mathrm{Gl}}^\nu(f,f)\le C_{\beta,k}\mathcal E_{\mathrm{Gl}}^\pi(f,f)\,,
\end{equation}
where the Dirichlet form under~$\nu$ only updates spins on ~$U\setminus A$. In particular, an~$(\alpha,\delta)$-WPI for~$\nu$  transfers to a~$(C_{\beta,k}(\alpha+1),C_{\beta,k}\delta)$-WPI for $\pi$.
\end{lemma}

\begin{proof}
    Set~$c=(1+e^{8\beta})^{-1}$ and let~$T$ be the map that sets the spins of~$A$ to~$\xi$. Every single-spin conditional probability is at least~$c$, so~$\pi\circ T^{-1}\le c^{-k}\nu$ and~$\mathcal E_{\mathrm{Gl}}^\nu\le c^{-k}\mathcal E_{\mathrm{Gl}}^\pi$. To bound~$f-f\circ T$, pin spins in $A$ successively. Before the~$j$th change, the law has Radon--Nikodym derivative at most~$c^{-(j-1)}$ to~$\pi$, and forcing the prescribed spin costs at most~$c^{-1}$ times the Dirichlet form for a single-spin resampling at that vertex. Telescoping and applying Cauchy--Schwarz give
\begin{equation*}
 \pi[(f-f\circ T)^2]\le k c^{-k}\sum_{v\in A}\mathcal E_{\mathrm{Gl},v}^\pi(f,f)\le k c^{-k} \mathcal E_{\rm Gl}^\pi(f,f)\,,
\end{equation*}
where $\mathcal E_{\mathrm{Gl},v}$ is the Dirichlet form associated to spin flips at $v$. 
Using~$\var_\pi(f)\le2\var_\pi(f\circ T)+2\pi[(f-f\circ T)^2]$ and~$\var_\pi(f\circ T)\le c^{-k}\var_\nu(f)$ we get~\eqref{eq:pin-comparison}. The WPI comparison follows by plugging in, and using~$\Osc(f|_{\{\sigma_A=\xi\}})\le\Osc(f)$.
\end{proof}

Given Lemma~\ref{lem:remove-pins-wpi}, we then can apply all our previous arguments to get the following analogue of Lemma~\ref{lem:diagonally-Dobrushin-WPI}, with the remaining adjustment in the proof being that finite-volume weights are with respect to the near-square rather than square. 

\begin{lemma}\label{lem:near-square-Dobrushin-WPI}
    Fix $\beta>\beta_c$. For every $C'$, there exists $C$ such that the following holds. Let $\Lambda^{(k^2+j)}$ be any near-square in $\Lambda_n$. Then Ising Glauber dynamics on $\Lambda^{(k^2 +j)}$ with the above induced diagonally Dobrushin boundary conditions satisfy an $(n^C, n^{-C'})$-WPI. 
\end{lemma}

\begin{proof}
    The proof is essentially identical to the proof of Lemma~\ref{lem:diagonally-Dobrushin-WPI}, except that we must check that the small kink in the geometry of the near-square as compared to the square does not affect the WPI constants by more than  a change of constants. By the cut-width bound of Lemma~\ref{lem:cutwidth-mixing-bound}, we can assume that $k \ge R \log n$ for a sufficiently large constant $R$ that it is much bigger than the $O(\log n)$-strip width and cutoffs used in the preceding arguments.

Consider $\Lambda^{(k^2 +j)}$ for some $1\le j\le 2k$, condition the spin at~$(1,k+1)$ to plus, and write~$\nu$ for the resulting law. Regard this vertex as part of the fixed boundary: this forces an initial south step in the interface, and the remainder of the interface joins source and sink points
\begin{equation*}
 a=(0,k),\qquad b=(k,0)\,,
\end{equation*}
as it would on $\Lambda_k$. Due to the cone constraint, every intervening bulk path connecting $a\to b$ lies in the true~$k\times k$-square $\Lambda^{(k^2)}$. This pinning of a vertex spin costs a constant $C_{\beta}$ in the WPI per Lemma~\ref{lem:remove-pins-wpi}. The Glauber WPI given the interface held for arbitrary domains with monochromatic boundary conditions and therefore is unchanged. 

When defining the interface dynamics (both finite-volume and infinite-volume) and spin block dynamics, we allow the strips to range one column further to the right, to include the partial column at $x$-coordinate $k+1$ so the strips still cover $\Lambda^{(k^2+j)}$. Then the analysis of the infinite-volume interface dynamics is identical to the case when we just were considering the domain $\Lambda_k$ from the start. The finite-volume comparison of the WPI using the Radon--Nikodym derivative $\psi_{\Lambda} = q_{\Lambda}(\Gamma)/q_\infty(\Gamma)$ extends to $\Lambda= \Lambda^{(k^2 + j)}$  via the near-squares generalization of Lemma~\ref{lem:OZ-finite-volume-comparison}.
The further comparisons to the Glauber Dirichlet form estimates apply with no trouble as they were built off of cut-width bounds that hold up to trivial changes of $k\mapsto k+1$. 
\end{proof}

\subsection{Censoring to get mixing time}

\begin{proof}[Proof of Theorem~\ref{thm:diagonally-Dobrushin-mixing}]
Recall the reduction of mixing time to mixing from the extremal initializations per~\eqref{eq:tv-controlled-by-extremal-inits} (taking $\Lambda$ to be $\Lambda_n$ and $\tau = {\rm Dob}$). By symmetry of the diagonal Dobrushin boundary conditions, it is sufficient to upper bound~\eqref{eq:tv-controlled-by-extremal-inits} with the  all-plus initialization. 

Let $\overline{X}_{+,t}$ be the censored version of $X_{+,t}$ with censoring scheme (recall Definition~\ref{def:censoring}) given by the growing domains $(\Lambda^{(i)})_{1\le i\le n^2}$ of Definition~\ref{def:near-squares} and times $t_0=0$ and $$t_{i}- t_{i-1} = C_0 n^{C_1}\qquad i\ge 1$$
where $C_1$ is bigger than the constant from Lemma~\ref{lem:near-square-Dobrushin-WPI}, and $C_0$ is a sufficiently large constant depending on~$\beta$. By the censoring Theorem 1.1 of~\cite{PWcensoring},  for every $t$ we have 
\begin{align}\label{eq:diagonal-Dobrushin-PI-almost-there}
   \sum_{v\in \Lambda_n} |\mathbb P(X_{+,t}(v) =1) - \pi_{\Lambda_n}^{\rm Dob}(\sigma_v =1)| \le n^2\|\mathbb P(X_{+,t}\in \cdot)  - \pi_{\Lambda_n}^{\rm Dob}\|_\tv \le n^2 \|\mathbb P(\overline{X}_{+,t}\in \cdot) - \pi_{\Lambda_n}^{\rm Dob}\|_{\tv} 
\end{align}
We next claim inductively, that for $0\le s\le C_0 n^{C_1}$, the law of $\overline{X}_{+,t_{i-1} +s}$ is within $in^{-10}$ total-variation distance of a chain on $\Lambda^{(i)}$, $\mathbb P_{\pi_{\Lambda^{(i-1)}}^{\rm Dob}\sqcup+}(Y_s^{(i)}\in \cdot)$, where $Y^{(i)}$ is an Ising Glauber dynamics on $\Lambda^{(i)}$ with $\rm Dob$ boundary conditions, and the initialization is from the previous domain's stationary distribution $\pi_{\Lambda^{(i-1)}}^{\rm Dob}$ concatenated with a $+1$ spin at the single vertex $\Lambda^{(i)}\setminus \Lambda^{(i-1)}$. The base case is vacuously true because the initialization on the single vertex $\Lambda^{(1)}$ will just be $+$ as the initialization of $\overline X_+$ is.
Now, suppose that is the case for an $i$, so that by the inductive assumption, 
$$\|\mathbb P(\overline{X}_{+,t_{i}}(\Lambda^{(i)})\in \cdot) - \mathbb P_{\pi_{\Lambda^{(i-1)}}^{\rm{Dob}}\sqcup +}(Y_{C_0 n^{C_1}}^{(i)}\in \cdot)\|_{\tv} \le in^{-10}\,.$$ 
The initial distribution $\pi_{\Lambda^{(i-1)}}^{\rm{Dob}}$ concatenated with a $+1$ at $\Lambda^{(i)}\setminus \Lambda^{(i-1)}$ has Radon--Nikodym derivative to $\pi_{\Lambda^{(i)}}^{\rm Dob}$ that is trivially bounded by $B = 1+e^{ 8 \beta}$. Therefore, by Lemma~\ref{lem:near-square-Dobrushin-WPI} and~\eqref{eq:WPI-consequence}, 
\begin{align*}
    \|\mathbb P_{\pi_{\Lambda^{(i-1)}\sqcup +}^{\rm{Dob}}}(Y_{C_0 n^{C_1}}^{(i)}\in \cdot) - \pi_{\Lambda^{(i)}}^{\rm Dob} \|_{\tv} \le n^{-10} \qquad\text{as long as  $t_i - t_{i-1} \ge n^{C} \log (e^{8\beta} n^{20} )$}
\end{align*}
(using that for $C'$ sufficiently large, $B n^{-C'/2} \le n^{-10}$). Adjusting the constants $C_0, C_1$, we evidently can ensure that for all $n$, we have $t_i - t_{i-1} \ge n^{C} \log ((1+e^{8\beta}) n^{20})$. Then take the optimal coupling of $\overline{X}_{+,t_{i}}(\Lambda^{(i+1)})$ with the sample $Y_0^{(i+1)}\sim \pi_{\Lambda^{(i)}}^{\rm Dob}\sqcup +$. Since $\Lambda^{(i+1)}\setminus \Lambda^{(i)}$ had been censored before time $t_{i}$, 
\begin{align*}
    \mathbb P(\overline{X}_{+,t_{i}} (\Lambda^{(i+1)}) \ne Y^{(i+1)}_0) \le (i+1)n^{-10}\,.
\end{align*}
Since the vertices of $\Lambda_n\setminus \Lambda^{(i+1)}$ are frozen to be plus for times $[t_i,t_{i+1})$ in the chain $\bar X_{+,t_{i}+s}$ they, together with the rest of the diagonally Dobrushin boundary conditions of $\Lambda_n$ impose diagonally Dobrushin boundary conditions on $\overline X_{+,t_i +s}$. Therefore using the same randomness to run the $\overline X_{+,t_i +s}$ and $Y^{(i+1)}_s$  Markov chains for times $0 \le s \le C_0 n^{C_1}$ retains 
\begin{align*}
    \mathbb P(\overline{X}_{+,t_{i}+s} (\Lambda^{(i+1)}) \ne Y^{(i+1)}_s) \le (i+1)n^{-10}\,,
\end{align*}
for all $0 \le s\le C_0 n^{C_1}$ as desired. In particular, $\|\mathbb P( \bar X_{+,t_{n^2}}\in \cdot) - \pi_{\Lambda_n}^{\rm Dob}\|_{\tv} \le n^{-8}$ which plugged into~\eqref{eq:diagonal-Dobrushin-PI-almost-there} and then~\eqref{eq:tv-controlled-by-extremal-inits} implies the mixing time is at most $t_{n^2} = C_0 n^{C_1 +2}$. 
\end{proof}

\appendix 
\section{Proofs of Ornstein--Zernike inputs}\label{app:OZ-inputs}\label{appendix:OZ-estimates}

We give the proofs of the estimates used in Section~\ref{sec:renewal-construction}. These are modifications of by-now-standard estimates in Ornstein--Zernike theory, with the flavor of differences to our setting being e.g., in the literature they have only been written for slope-$0$ instead of slope-$-1$, or written for percolation or random-cluster OZ theory rather than the random-line representation of the Ising model. 

\subsection{Comparison of domains}\label{proof:domain-weight-comparison}

\begin{lemma}\label{lem:deletion}
Fix $\beta>\beta_c$. For a finite dual graph~$H$ and sets of edges~$A,B$, put
\begin{equation*}
 F_H(A,B)=\frac{Z_0(H\setminus(A\cup B))Z_0(H)}
                    {Z_0(H\setminus A)Z_0(H\setminus B)}\,.
\end{equation*}
There exist constants $c_\beta,C_\beta$ such that 
\begin{equation}\label{eq:deletion-influence}
 0\le\log F_H(A,B)\le C_\beta\sum_{e\in A,\,f\in B}e^{-c_\beta d(e,f)}\,.
\end{equation}
\end{lemma}

\begin{proof}
Iteratively restore the edges of~$A\cap E(H)$ in~$H\setminus A$, and as each edge is added, compare the factors in~\eqref{eq:edge-addition} with the same factors when adding those edges in~$H\setminus (A\cup B)$. For~$e=uv\notin B$, write~$K$ for the graph just before adding~$e$. By GKS each such edge has additive contribution to~$\log F_H(A,B)$ between zero and
\begin{equation*}
 e^{-2\beta}\bigl(\langle\sigma_u\sigma_v\rangle_{K,\beta^*}
       -\langle\sigma_u\sigma_v\rangle_{K\setminus B,\beta^*}\bigr).
\end{equation*}
To control this quantity, turn on the couplings of the edges~$f=xy$ in~$B\cap E(K)$, tuning their strength from zero to~$\beta^*$. The derivative of~$\langle\sigma_u\sigma_v\rangle$ in the coupling of~$f$ is~$\cov(\sigma_u\sigma_v,\sigma_x\sigma_y)$. By the four-point inequality~\cite[p.~91, remark (i)]{Lebowitz74} and GKS, it is at most
\begin{equation*}
 G_\beta(u-x)G_\beta(v-y)+G_\beta(u-y)G_\beta(v-x)\le C_\beta e^{-c_\beta d(e,f)}\,.
\end{equation*}
Integrating this over $\beta'\in (0,\beta^*)$, then summing over~$e,f$ proves the lemma. (Note that for $e=f$, the quantity is bounded by $\log (1+e^{-2\beta})$ which is also bounded by its summand in ~\eqref{eq:deletion-influence}. 
\end{proof}

\begin{proof}[Proof of Lemma~\ref{lem:domain-weight-comparison} ]
For~$\Lambda_1\subset\Lambda_2$, \eqref{eq:q-Lambda} gives~$q_{\Lambda_1}(\eta)/q_{\Lambda_2}(\eta)=F_{\Lambda_2}(Q_\eta,B)$. Lemma~\ref{lem:deletion} proves the bound and the sign. For general domains, compare each to~$\Lambda_1\cap\Lambda_2$. The two sets of deleted edges partition~$B$, giving~\eqref{eq:domain-weight-comparison}. 
\end{proof}

\begin{proof}[Proof of Lemma~\ref{lem:domain-law-comparison}] On the common family~$\mathcal A$, the density is~$\psi/\nu_2(\psi)$ for~$\psi=q_{\Lambda_1}/q_{\Lambda_2}$, which gives~\eqref{eq:domain-radon-nikodym}. Summing the exponential over all edges at distance at least~$r$ from a fixed edge costs~$O(e^{-cr})$, since there are~$O(\ell+1)$ edges at distance~$\ell$. This gives both the bound for changing~$k$ edges and the bound for a change at distance~$r$.
\end{proof}

\subsection{Number of cone-points}\label{proof:cone-time}

Let~$\phi_H$ be the free FK--Ising measure on the dual graph~$H$ at edge parameter~$p=1-e^{-2\beta^*}$. Writing~$Z_{a,b}(H)$ for the sum of~$e^{-2\beta|F|}$ over subgraphs with odd-degree vertices~$a,b$, the high-temperature expansion and the Edwards--Sokal coupling give
\begin{equation}\label{eq:odd-fk-normalization}
 \frac{Z_{a,b}(H)}{Z_0(H)}
 =\langle\sigma_a\sigma_b\rangle_{H,\beta^*}
 =\phi_H(a\leftrightarrow b)\,.
\end{equation}
For an FK-cluster $K$, define cone-points as in Definition~\ref{def:cone-point} and denote their number by~$N_{\mathcal C}(K;a,b)$.

\begin{lemma}\label{lem:cone-time}
Fix~$\beta>\beta_c$ and a compact subinterval of~$(-2,0)$. There are~$\rho,c,C>0$, depending only on~$\beta$ and this subinterval, such that, whenever~$x=b-a$ has positive horizontal coordinate and slope in this subinterval, every finite lattice subgraph~$H$ satisfies
\begin{equation}\label{eq:no-cone-tail}
 \phi_H\bigl(a\leftrightarrow b,\ N_{\mathcal C}(C_a;a,b)<\rho\|x\|_1\bigr)
 \le C e^{-c\|x\|_1}G_\beta(x)\,.
\end{equation}
The same estimate holds on the full lattice.
\end{lemma}

Every  path from~$a$ to~$b$ in the cluster visits each separating cone-point. The coupling of~\eqref{eq:odd-fk-normalization} thus gives the corresponding random-line estimate.

\begin{lemma}\label{lem:contour-cone-deficiency}
Under the same restrictions on the endpoints $a,b$ as in Lemma~\ref{lem:cone-time} and for any domain $\Lambda$ including the full lattice,
\begin{equation}\label{eq:physical-cone-deficiency}
 \sum_{\substack{\Gamma:a\to b\\ \#\{\text{interior cone-points of }\Gamma\}<\rho\|x\|_1}}
 q_\Lambda(\Gamma)
 \le C e^{-\xi(x)-c\|x\|_1}\,.
\end{equation}
\end{lemma}

\begin{proof}[Proof of Lemmas~\ref{lem:cone-time} and~\ref{lem:contour-cone-deficiency}]
The paper~\cite{CIV08} gave a decomposition of the subcritical FK cluster connecting $a\to b$ with constant density of ``sparse points". By doing some finite-size modifications on a random subset of sparse points, we will show a constant density of cone-points. 

Write~$n=\|b-a\|_1$ and  work in the full lattice $\mathbb Z^2$. At the dual temperature~$\beta^*<\beta_c$, the FK--Ising measure is unique and its connection probabilities decay exponentially. By e.g.,~\cite[Proposition~A.1]{CIV08}, in a square of radius~$R$ the probability that the center connects to the boundary is at most~$Ce^{-cR}$ even when all boundary vertices are identified. This is assumption~(1.2) of~\cite{CIV08}.

On~$a\leftrightarrow b$, represent the cluster~$C_a$ by the tree constructed in~\cite[Section~2.2]{CIV08} at scale~$K$. Its vertices are centers of balls of radius~$O(K)$ covering~$C_a$, and its \emph{trunk} is the tree path from~$a$ to a center whose ball covers~$b$. Choose a closed cone~$\mathcal C'$ whose interior contains every allowed endpoint direction and whose nonzero directions lie strictly inside~$\mathcal C$. For an endpoint slope~$s$, the forward cone in~\cite[(2.9)]{CIV08} consists of the directions~$y$ satisfying~$h(s)\cdot y>(1-2\delta)\xi(y)$. By strict convexity, $h(s)\cdot y=\xi_t(y)$ only on the forward ray through~$(1,s)$. Compactness of the slope interval therefore lets us choose~$\delta>0$ so that all these forward cones lie in~$\mathcal C'$.

Fix~$K$ sufficiently large for this~$\delta$. By~\cite[Theorem~2.1]{CIV08}, conditional on~$a\leftrightarrow b$, except with probability~$Ce^{-cn}$ there are at least~$c n/K$ trunk vertices~$v$ at which the earlier trunk vertices lie in~$v-\mathcal C'$, the later ones lie in~$v+\mathcal C'$, and every branch lies in the union of the two cones. Since the balls centered at the tree vertices cover~$C_a$, each such~$v$ satisfies
\begin{equation}\label{eq:coarse-cone-opportunity}
 C_a\subseteq\bigl((v-\mathcal C')\cup(v+\mathcal C')\bigr)+[-h_0,h_0]^2
\end{equation}
with~$h_0=O(K)$. We now modify the edges near~$v$ to make it a separating cone-point of the cluster for our cone~$\mathcal C$ of Definition~\ref{def:cone-point}, adapting~\cite[Section~2.9]{CIV08}.

Take fixed~$r$ large enough that outside the $\ell^\infty$ ball~$B_{r-1}(v)$ the two sets in~\eqref{eq:coarse-cone-opportunity} lie in~$v\pm\mathcal C$, respectively. Via $O(r^2)$ many forced modifications on the states of the edges in $B_{r-1}(v)$, we can make it so that the two paths stay in their respective cones everywhere and meet exactly at the vertex $v$, while keeping the cluster outside $B_{r-1}(v)$ fixed.  

By the strict inclusion of the cone directions, for fixed~$L$ large, the square~$B_r(v)$ lies entirely in one of~$z\pm\mathcal C$ for every cluster vertex~$z$ outside~$B_{Lr}(v)$. Thus the modification preserves whether such a vertex is a cone-point and creates at most a fixed number~$M= O(r)$ many cone-points inside $B_r(v)$. Choose~$\alpha n$ centers away from~$a,b$ with disjoint squares~$B_{Lr}(v)$ for some fixed~$\alpha>0$. Starting from a configuration with fewer than~$\rho n$ separating cone-points, modify any~$m=\delta n$ of these centers. The resulting cluster still connects~$a\to b$ and has at most~$\rho n+Mm$ separating cone-points, among them the modified centers. Specifying these centers and the old states of the edges incident to their squares recovers the original configuration. By enumerating the preimages of the modifications on the balls and accepting a cost of $C^m$ for the change in probability for these modifications, we find that the probability of such a configuration is bounded by 
\begin{equation*}
 C^m\frac{\binom{\rho n+Mm}{m}}{\binom{\alpha n}{m}}
 \le \left[\frac{C'(\rho+M\delta)}{\alpha}\right]^m \,.
\end{equation*}
Taking~$\delta,\rho$ sufficiently small gives an $e^{-\Omega(n)}$ upper bound on the probability. Adding the $e^{-\Omega(n)}$ probability that the tree in the random-cluster side was atypical gives the claim on $\mathbb Z^2$. 

For a subgraph~$H$ of $\mathbb Z^2$, use the monotone coupling of~\cite[proof of Lemma~2.2]{IOVW-invarianceprinciple}. A separating cone-point of the larger cluster is its only vertex in that column, so every $a$-to-$b$ path in the smaller cluster passes through it. Thus the event in~\eqref{eq:no-cone-tail}, including the connection, is increasing, and the probability becomes smaller when shrinking the domain. Finally, recall the coupling of high-temperature FK connections and low-temperature OZ interfaces; namely, from $F(\sigma)$ independently opening each dual edge outside~$F$
with probability~$\tanh(\beta^*)$ gives a sample from the 
FK--Ising measure at~$\beta^*$ conditioned on~$a\leftrightarrow b$. Under this coupling, each separating random-cluster cone-point is also a cone-point of~$\Gamma$. Applying~\eqref{eq:odd-fk-normalization} and~$G_\beta(b-a)\le e^{-\xi(b-a)}$ gives~\eqref{eq:physical-cone-deficiency}.
\end{proof}

\subsection{Renewal components}\label{app:renewal}\label{proof:renewal-expansion}\label{app:boundary-modifications}
\begin{proof}[Proof of Propositions~\ref{prop:domain-renewal} and~\ref{prop:domain-component-bounds}]
We first prove both propositions for the internal renewal components, working in the full lattice $\mathbb Z^2$. At the end of the proof, we describe the modifications for handling the boundary renewal components.  

The interval expansion of~\cite[Theorem~5.9 and Appendix~A.2]{IOSV-prewetting} gives that the weight of a decorated OZ path $(\Gamma,\mathfrak{n})$ is exactly the product of~\eqref{eq:decorated-component-weight} over its renewal components. Indeed the factors~$1+f_I$ telescope to~$q_\infty(\Gamma)/\prod_i q_\infty(\eta_i)$, and multiplication by~$\mathcal P_\Gamma(\mathfrak n)$ leaves one factor~$f_I$ for each open interval. Thus the path weights tensorize across the open components of $\mathfrak{n}$. 

\smallskip
\noindent \emph{Exponential moments.} Fix~$h=h(s)$, $s\in\mathfrak S$. Consider a path~$w$ of length~$m$ formed by concatenating~$r$ constituent paths, including any boundary paths. For fixed displacement~$x$, its total tilted weight is at most~$C e^{h\cdot x}G_\beta(x)\le C$ by~\eqref{eq:raw-source-sum}. There are~$O(m^2)$ possible displacements. If~$r\ge\varepsilon m$, the probability that the decoration connects all~$r$ indices is at most~$Ce^{-c\varepsilon m}$ by~\cite[Appendix~A.3.1, (A.83)]{IOSV-prewetting}. This estimate uses the independent interval law~\eqref{eq:conditional-mark-law} and the exponential bound on~$f_I$.

It remains to bound the tilted weight of paths with~$r<\varepsilon m$. If~$\|x\|_1\le\delta m$, the length bound~\eqref{eq:raw-path-length-tail} absorbs~$e^{h\cdot x}$ for~$\delta$ sufficiently small. For~$\|x\|_1>\delta m$, fix a compact interval~$J\subset(-2,0)$ containing~$\mathfrak S$ in its interior. If~$x$ points in a direction~$(1,s')$ with~$s'\in J$, Lemma~\ref{lem:contour-cone-deficiency} gives a bound~$Ce^{-c\|x\|_1}$ for the tilted sum, since~$h\cdot x\le\xi(x)$ and~$\varepsilon m<\rho\|x\|_1$ for~$\varepsilon$ sufficiently small. In all other directions, strict convexity and compactness give~$\xi(x)-h\cdot x\ge c\|x\|_1$, and the two-point bound gives the same $Ce^{ - c\|x\|_1}$ bound. Summing over possible $x$, the total tilted weight of a renewal component of length~$m$ is at most~$Ce^{-cm}$. 

\smallskip
\noindent 
\emph{Normalization and mean slope.} To see why the bulk weights have total weight one, let~$F(z)$ be their sum at tilt~$z$ and let~$\chi(z)$ be the sum of all path weights at this tilt, summed over displacements and boundary directions. The random-line representation gives~$\chi(z)\asymp\sum_x e^{z\cdot x}G_\beta(x)$, where $\asymp$ denotes proportionality up to constant. Summing the representation~\eqref{eq:decorated-component-weight} over the number of bulk components gives
\begin{equation*}
 \chi(z)=b_L(z)b_R(z)\sum_{j\ge0}F(z)^j+R(z)\,,
\end{equation*}
where~$b_L,b_R$ are the two boundary masses and~$R$ is the mass of paths consisting of a single component. The boundary masses are positive, and the exponential moments just proved make all four functions~$F,b_L,b_R,R$ analytic near~$h$. Thus~$F(z)<1$ for~$z\in\operatorname{int}\mathcal K_{\beta}$ near~$h$. If~$F(h)<1$, then the geometric series would converge throughout a neighborhood of $h$, so $\sum_x e^{z\cdot x}G_\beta(x)$ would be finite throughout that neighborhood. This contradicts $h\in\partial\mathcal K_\beta$ from the definition~\ref{eq:dual-convex-set}. By continuity, $F(h)=1$.

 Since~$\partial_1F(h)=\E T>0$, write the level curve~$F=1$ near~$h$ as~$(\lambda(\theta),\theta)$. Differentiating twice, with~$\rho=\E Y/\E T$, gives
\begin{equation}\label{eq:physical-slope-identity}
 \lambda'(\theta)=-\rho,\qquad
 \rho'(\theta)=\frac{\E[(Y-\rho T)^2]}{\E T}>0\,.
\end{equation}
At~$h=h(s)$, $(1,s)$ is normal to the level curve, so~$\rho=s$.

\smallskip
\noindent \emph{Restriction to the square.}
Once the two boundary components are chosen, the intervening internal renewal components stay between their endpoint columns and heights by cone confinement. These cones are necessarily confined to the square, and therefore the laws of these internal components do not depend on the underlying square. To get~\eqref{eq:decorated-interface-factorization}, collect all configurations with fewer than three renewal components in the exceptional weight $R_{\Lambda,h}$ and otherwise, count the weights of the two boundary components, and then draw the internal components given those.  

To see~\eqref{eq:domain-renewal}, given all the renewal points under~$\widehat\pi^\infty$, the displacements of the individual components are fixed, so the product of their weights normalizes separately on each segment. The same argument applies in the presence of the boundary measures, whenever every intervening bulk concatenation between them is admissible. 
Note that this is also true in the near-square with its top-left vertex pinned to $+1$ from Section~\ref{sec:near-squares}, as the cones of the boundary components will lie entirely in $\Lambda_k$.

For the exponential moments of the boundary measures in Proposition~\ref{prop:domain-component-bounds}, note that their exponential moment only decreases from the restriction that the boundary imposes. In the square or near-square, each boundary measure retains a bounded path of $\Omega(1)$ weight. Finally, the contribution to the exceptional path weight~$R_{\Lambda,h}$ coming from having exactly the two boundary components and no internal components is bounded by the product of their exponential moments. This proves~\eqref{eq:domain-component-moments}.
\end{proof}

\subsection{OZ bridges under a length cutoff}\label{proof:renewal-estimates}

To construct the tilt $h_Q$ that normalizes the truncated weights on renewal components in~\eqref{eq:tilted-p-distribution}, write~$h=(\lambda(\theta),\theta)$ and let~$\theta$ range over a compact interval whose mean slopes contain~$\mathfrak S$ in their interior. Choose~$\delta_Q=\delta_Q(\theta)$ so that
\begin{equation}\label{eq:cutoff-normalization}
 \sum_{|\omega|\le Q} \mathsf p_h(\omega)e^{\delta_Q T(\omega)}=1\,.
\end{equation}
The discarded mass and its first two derivatives in~$\theta$ are~$O(e^{-cQ})$ by Proposition~\ref{prop:domain-component-bounds}. The left side is strictly increasing in~$\delta_Q$, with derivative bounded below since~$T\ge1$. Thus~$\delta_Q$ and its first two derivatives are~$O(e^{-cQ})$. Differentiating~\eqref{eq:cutoff-normalization} gives mean slope~$-(\lambda+\delta_Q)'(\theta)$. By~\eqref{eq:physical-slope-identity} this slope is strictly increasing and its range still contains~$\mathfrak S$ for~$Q$ large. Choosing~$\theta$ to give slope~$s$ and putting~$h_Q(s)=(\lambda(\theta)+\delta_Q(\theta),\theta)$ gives the claimed law. The tilted boundary measures retain their exponential moments and their options of bounded paths of positive weight.

\begin{proof}[Proof of Lemma~\ref{lem:OZ-local-limit}]
The local limit theorem type results of Lemma~\ref{lem:OZ-local-limit} already appear for fixed cutoffs $Q$ in~\cite{Ioffe14} and fixed slope. In the below we justify that the constants therein are uniform over all $Q \ge Q_0$ and over all slopes $s\in \mathfrak{S}$. 

Applying~\cite[Section~3.1, Proposition~3.2, (83) and (89)--(96)]{Ioffe14} with time~$T$ and spatial displacement~$Y$, Fourier inversion gives the uniform upper bound and a lower bound of order~$x^{-1/2}$ exactly at the lattice point $y$ nearest the chord of slope $s$. When considering general~$|y-sx|\le A\sqrt x$, we  add a tilt to make the mean slope~$y/x$. The change~$\Delta\theta$ in the vertical tilt is~$O(A/\sqrt x)$ by~\eqref{eq:physical-slope-identity} and~\eqref{eq:cutoff-normalization}. Since~$(\lambda+\delta_Q)'=-s$ at the original tilt, the log change of weight at~$(x,y)$ is~$(y-sx)\Delta\theta+O(x|\Delta\theta|^2)=O(A^2)$. This proves the lower bound throughout the stated window.

The Fourier estimates in~\cite{Ioffe14} are uniform in~$s,Q$. Indeed the step laws, including the untruncated law, form a compact family in an exponentially weighted~$\ell^1$ norm, and the three internal renewal components ($E, ESE, ESSE$) have uniformly lower  bounded weights and the differences of their displacement vectors generate $\mathbb Z^2$. This bounds the covariance away from degeneracy and the modulus of the characteristic function away from $1$ away from the origin.

To incorporate the convolutions with the two boundary renewal components, we simply condition on the two boundary components. For the lower bound, choose a bounded path of positive weight at each end. For the upper bound, if their combined length is at most~$x/2$, the remaining width is at least~$x/2$ and the bulk upper bound applies. Here the Fourier bound for the bulk holds for all sufficiently large widths, uniformly in~$Q$. The complementary event has probability~$O(e^{-cx})$, and~$U_{s,Q}\le1$ since a bulk walk visits any given column at most once.
\end{proof}

\begin{proof}[Proof of Lemma~\ref{lem:OZ-bridge-concentration}]
Translate~$a$ to the origin and write~$b=(\ell,s\ell)$. For~$s\in\mathfrak S$, sample the left boundary component and then independent bulk components of mean slope~$s$, up to the first renewal point at or beyond~$\ell/2$. By~\eqref{eq:domain-renewal}, relative to this law, the density of the sampled portion under the bridge is the renewal mass of the remaining displacement, convolved with the right boundary law, divided by the renewal mass for the whole bridge with both boundary laws. By~\eqref{eq:OZ-local-limit} this density is bounded. 

By Proposition~\ref{prop:domain-component-bounds} the total centered bulk displacement~$Y-sT$ (evaluated at the full bulk renewal portion $\omega_1 \circ \cdots \circ \omega_k$) has mean zero and a uniform exponential moment, so~$\E e^{\vartheta(Y-sT)}\le e^{C\vartheta^2}$ for~$|\vartheta|$ small. The maximal inequality for its exponential martingale bounds the probability of a deviation at least~$u$ in the first~$\ell$ steps by~$C\exp\{-c\min(u^2/\ell,u)\}$. There are at most~$\ell$ steps before the crossing, and the initial boundary displacement independently has an exponential tail. The $O(1)$ density bound described in the first paragraph of this proof therefore gives~\eqref{eq:OZ-bridge-concentration} up to the crossing. Reading the bridge backwards from~$b$ gives the other half and we can union bound over the two.

For arbitrary feasible slopes $\rho\in [-2,0]$, condition on the boundary components and on the entire set of bulk renewal components, but forget their order. The bulk components therefore occur in a uniformly random order. If their total displacement has slope~$\rho$, put~$b_i=Y_i-\rho T_i$ for the $i$'th centered displacement. Then~$\sum_i b_i=0$ and, by cone confinement and the length cutoff,
\begin{equation*}
 \sum_i b_i^2\le C\sum_i T_i^2\le CQ\ell\,,
\end{equation*}
since $Y_i/T_i$ and the internal slope both lie in $[-2,0]$. 
For~$k$ numbers with zero sum, the partial sums~$S_j$ in uniformly random order satisfy~$\E\max_j|S_j|^2\le C\sum_i b_i^2$. Indeed~$S_j/(k-j)$ is a martingale, and~$\E S_j^2=j(k-j)\sum_i b_i^2/[k(k-1)]$, so this follows from Doob's inequality applied to the first half (and in reverse order for the second). Passing from the bulk chord to the chord from~$a$ to~$b$ changes the deviations by at most~$CQ$, since the two end components have length at most~$Q$ and the bulk slopes are bounded. As~$Q^2\le Q\ell$, this proves~\eqref{eq:OZ-bridge-rough-moment}.
\end{proof}

\begin{proof}[Proof of Lemma~\ref{lem:OZ-exponential-tails}]
There are at most~$n$ many renewal components since each advances at least one column. Before conditioning, the exponential renewal component tail bound of Proposition~\ref{prop:domain-component-bounds}, together with a union bound over these $n$ steps and the two boundary components costs~$Cne^{-cQ}$. Dividing by the normalizing constant of~\eqref{eq:physical-reference-normalizer} costs an extra factor of~$C\sqrt n$. Finally, the exceptional probability is at most~$C\sqrt n e^{-cn}$ by~\eqref{eq:exceptional-paths-unlikely}, proving~\eqref{eq:OZ-cutoff-error}.
\end{proof}

\subsection{Finite-volume weights}\label{proof:boundary-influence}

\begin{proof}[Proof of Lemma~\ref{lem:OZ-finite-volume-comparison}]
We write the proof in a unified manner for both the square $\Lambda_n$ and for a near-square $\Lambda^{(n^2+j)}$. Denote the endpoints of the interface as~$a=(0,u)$ and~$b=(n,0)$, noting that $u\in \{n,n+1\}$. The slope will be $s=(v-u)/n$, and the tilt will be~$h=h(s)$. The domain contains $\Lambda_n$ and is contained in $\Lambda_{n+1}$; the chord from $a\to b$ differs from~$y=n-x$ by at most one. Lemma~\ref{lem:domain-weight-comparison} gives~$\psi_\Lambda=q_\Lambda/q_\infty\ge1$, so normalizing the two weights gives~\eqref{eq:OZ-finite-volume-density}. 

It remains to bound the upper tail of~$\psi_\Lambda$. Put~$d(x)=\min\{|x|,|n-x|\}$ and fix~$\eta>0$ small. Let~$E_R$ be the event that every~$z=(x,y)\in\Gamma$ satisfies
\begin{equation*}
 |y-u-sx|\le R+\eta d(x)\,.
\end{equation*}
For a vertex violating this bound, choose~$w$ from~$z-a,b-z$ with~$|w_1|=d(x)$. Then~$|w_2-sw_1|>R+\eta|w_1|$. Strict convexity of~$\xi$ gives~$\xi(w)-h\cdot w\ge c|w|$ in these directions: on the unit circle they are separated from the ray~$(1,s)$ where equality holds in~\eqref{eq:dual-vector}. The constant is uniform over our slopes by compactness. The other displacement has nonnegative excess as well, whence
\begin{equation*}
 \xi(z-a)+\xi(b-z)-h\cdot(b-a)\ge cR\,.
\end{equation*}
By~\eqref{eq:through-point}, the tilted weight of all paths visiting such a vertex is at most~$Ce^{-cR}$. There are~$O(n^2)$ vertices. The same upper bound applies after restricting admissible paths, and also with~$q_\infty\le q_\Lambda$. Dividing by the tilted normalizing constants, both at least~$cn^{-1/2}$ by~\eqref{eq:physical-reference-normalizer}, gives probability at most~$Cn^{5/2}e^{-cR}$ for~$E_R^c$ under either law.

On~$E_R$ apply~\eqref{eq:domain-weight-comparison} with~$B$ the edges missing from~$\Lambda$. For each~$e\in Q_\Gamma$, summing~$e^{-c d(e,f)}$ over~$f\in B$ costs at most~$Ce^{-c' d(e,B)}$, since the number of edges at distance~$r$ from~$e$ is~$O(r+1)$. For an edge inside the contained square, this is at most a constant times the sum of the four exponential kernels in its distances to the sides. The same bound holds in the extra layer, where one of the four distances is bounded. In column~$x$ the envelope allows~$O(R+d(x))$ reserved edges. The contribution of the vertical sides is thus bounded by
\begin{equation*}
 C\sum_{d\ge0}(R+d)e^{-cd}\le CR\,.
\end{equation*}
For either horizontal side, the exponential summed over one column is bounded by a constant. The columns with~$d(x)\le R/\eta$ therefore contribute~$O(R)$. In the remaining columns the chord is at distance at least~$d(x)-1$ from both horizontal sides, and the envelope has width at most~$2\eta d(x)$. Their contribution obeys the same exponential sum. Thus~$\log\psi_\Lambda\le CR$ on~$E_R$. Together with the estimate on~$E_R^c$ this proves~\eqref{eq:OZ-likelihood-tail}.
\end{proof}

\bibliographystyle{alpha}
\small
\bibliography{references}

\end{document}